\documentclass[12pt]{article}
\usepackage{amssymb}
\usepackage{latexsym} 
\usepackage{amsmath}

\def\Xint#1{\mathchoice {\XXint\displaystyle\textstyle{#1}} 
{\XXint\textstyle\scriptstyle{#1}}
{\XXint\scriptstyle\scriptscriptstyle{#1}} 
{\XXint\scriptscriptstyle\scriptscriptstyle{#1}}\!\int} 
\def\XXint#1#2#3{{\setbox0=\hbox{$#1{#2#3}{\int}$} 
\vcenter{\hbox{$#2#3$}}\kern-.5\wd0}}

\makeatletter
\let\@xp=\expandafter

\let\@nx=\noexpand
 
\makeatother
\numberwithin{equation}{section}
\newtheorem{thm}[equation]{Theorem}
\newtheorem{prop}[equation]{Proposition}
\newtheorem{defn}[equation]{Definition}
\newtheorem{rem}[equation]{Remark}
\newtheorem{lem}[equation]{Lemma}
\newtheorem{corol}[equation]{Corollary}
\title{Introduction to Morrey spaces}
\author{Massimo Lanza de Cristoforis \\
Dipartimento di Matematica, Tullio Levi-Civita\\
Universit\`{a} degli Studi di Padova,\\
Via Trieste 63,\\
35121 Padova, Italia. \\
Fax: ++39 049 827 1204\\
e-mail: mldc@math.unipd.it\\
Draft of November 21,  2020. Updated December 28, 2023.}

\date{\ }

\begin{document}

\maketitle

\vfill
\begin{center}
Southern Federal University, Rostov-on-Don, Russia, 2020
\end{center}

\newpage

\section{Introduction.} The present informal set of notes covers the material that has been presented by the author in a series of lectures for the Doctoral School 
in Mathematics of the
Southern Federal State University of Rostov-on-Don in the Fall of 2020 and that develops from the first part of the notes of \cite{La12}  that collects the material of the lectures of the author at the Eurasian National University, Astana, Kazakhstan in the Spring of 2012. The aim is to present 
 some elementary classical properties of Morrey spaces and some corresponding approximation results by smooth functions. 


\section{Technical preliminaries} 
Let ${\mathbb{N}}$ denote the set of natural numbers including $0$. Let ${\mathbb{R}}$ and ${\mathbb{C}}$ denote the sets of real and complex numbers, respectively. Let 
\[
n\in {\mathbb{N}}\setminus\{0\}\,.
\]
Let 
${\mathbb{D}}\subseteq {\mathbb {R}}^{n}$. Then $\overline{{\mathbb{D}}}$ 
denotes the 
closure of ${\mathbb{D}}$ and $\partial {\mathbb{D}}$ denotes the boundary of ${\mathbb{D}}$. For all $R>0$, $ x\in{\mathbb{R}}^{n}$, 
$x_{j}$ denotes the $j$-th coordinate of $x$, 
$| x|$ denotes the Euclidean modulus of $ x$ in
${\mathbb{R}}^{n}$, and ${\mathbb{B}}_{n}( x,R)$ denotes the ball $\{
y\in{\mathbb{R}}^{n}:\, | x- y|<R\}$. 

Let ${\mathbb{D}}\subseteq {\mathbb {R}}^{n}$. Then $B( {\mathbb{D}} )$ denotes the space of bounded functions from ${\mathbb{D}}$ to ${\mathbb{R}}$.
Then 
$C^{0}_{u}( {\mathbb{D}} )$ denotes the subspace of those $f\in C^{0} ( {\mathbb{D}}  )$ which are uniformly continuous, and
$C^{0}_{b}( {\mathbb{D}} )$ denotes $B( {\mathbb{D}} )
\cap C^{0} ( {\mathbb{D}}  )$.  Then we set $C^{0}_{ub}(  {\mathbb{D}}  )=C^{0}_{u}(  {\mathbb{D}}  )\cap C^{0}_{b}(  {\mathbb{D}}  )$.
It is well known that $B( {\mathbb{D}} )$ with the sup-norm is a Banach space, and that $C^{0}_{b}( {\mathbb{D}} )$ is a closed subspace of $B( {\mathbb{D}} )$. 
Let $\Omega$ be an open subset of ${\mathbb{R}}^{n}$. The space of $m$ times continuously 
differentiable real-valued functions on $\Omega$ is denoted by 
$C^{m}(\Omega)$. 
If ${\mathbb{D}}\subseteq {\mathbb {R}}^{n}$, then 
\begin{equation}\label{eq:ctildom}
\tilde{C}^{m}({\mathbb{D}})
\end{equation}
 denotes the set of functions $f$ from ${\mathbb{D}}$ to ${\mathbb{R}}$ such that there exist an open subset $W_{f}$ of ${\mathbb {R}}^{n}$ containing
${\mathbb{D}}$ and a function $F_{f}\in C^{m}(W_{f})$ such that $f=F_{f}$ in ${\mathbb{D}}$. 

$C^{\infty}_{c}(\Omega)$ denotes the space of functions of  $C^{\infty}(\Omega)$
with compact support.

Let  $\eta\equiv
(\eta_{1},\dots ,\eta_{n})\in{\mathbb{N}}^{n}$, $|\eta |\equiv
\eta_{1}+\dots +\eta_{n}  $. Then $D^{\eta} f$ denotes
$\frac{\partial^{|\eta|}f}{\partial
x_{1}^{\eta_{1}}\dots\partial x_{n}^{\eta_{n}}}$.    The
subspace of $C^{m}(\Omega )$ of those functions $f$ whose derivatives $D^{\eta }f$ of
order $|\eta |\leq m$ can be extended with continuity to 
$\overline{\Omega}$  is  denoted $C^{m}(\overline{\Omega} )$. Obviously, we have $\tilde{C}^{m}(\overline{\Omega})\subseteq
C^{m}(\overline{\Omega} )$.

If $\Omega$ is an open subset of ${\mathbb{R}}^{n}$, $m\in {\mathbb{N}}$, we set
\[
C^{m}_{b}(\overline{\Omega})\equiv
\{
u\in C^{m}(\overline{\Omega}):\,
D^{\gamma}u\ {\mathrm{is\ bounded}}\ \forall\gamma\in {\mathbb{N}}^{n}\
{\mathrm{such\ that}}\ |\gamma|\leq m
\}\,,
\]
and we endow $C^{m}_{b}(\overline{\Omega})$ with its usual  norm
\[
\|u\|_{ C^{m}_{b}(\overline{\Omega}) }\equiv
\sum_{|\gamma|\leq m}\sup_{x\in \overline{\Omega} }|D^{\gamma}u(x)|\qquad\forall u\in C^{m}_{b}(\overline{\Omega})\,. 
\]

We denote by $m_{n}$ the $n$-dimensional Lebesgue measure
and by ${\mathcal{L}}_{n}$  the set of Lebesgue measurable subsets of
${\mathbb{R}}^{n}$.  We
denote by  $\omega_{n}$ the  $n$-dimensional
measure of the ball ${\mathbb{B}}_{n}(0,1)$. 

Let $(X,{\mathcal{M}},\mu)$ be a measured space. If
$p\in[1,+\infty[$, we denote by ${\mathcal{L}}^{p}_{\mu}(X)$ the
space of measurable functions from $X$ to ${\mathbb{R}}$ such that
\[
\left( \int_{X}|f|^{p}\,d\mu \right)^{1/p}<+\infty\,.
\]
Then we set
\[
\|f\|_{ {\mathcal{L}}^{p}_{\mu}(X) }\equiv \left(
\int_{X}|f|^{p}\,d\mu \right)^{1/p}\qquad\forall f\in
{\mathcal{L}}^{p}_{\mu}(X)\,.
\]
We denote by ${\mathcal{L}}^{\infty}_{\mu}(X) $  the space of
measurable functions from $X$ to ${\mathbb{R}}$ such that
\[
{\mathrm{ess\,sup}}_{X}|f|<+\infty\,.
\]
Then we set
\[
\|f\|_{ {\mathcal{L}}^{\infty}_{\mu}(X) }\equiv
{\mathrm{ess\,sup}}_{X}|f| \qquad\forall f\in
{\mathcal{L}}^{\infty}_{\mu}(X)\,.
\]
For each $p\in[1,+\infty]$, we denote by $L^{p}_{\mu}(X)$ the
quotient space of ${\mathcal{L}}^{p}_{\mu}(X)$ with respect to the
equivalence relation of equality almost everywhere. If $[f]\in
L^{p}_{\mu}(X)$, then we set
\[
\|[f]\|_{ L^{p}_{\mu}(X) }\equiv \|f\|_{
{\mathcal{L}}^{p}_{\mu}(X) }\,.
\]
Here $[f]$ denotes the equivalence class of $f$. As customary, we
shall simply write $f$ instead of $[f]$. As is well known,
$(L^{p}_{\mu}(X) ,\|\cdot \|_{ L^{p}_{\mu}(X) })$ is a Banach
space for all $p\in[1,+\infty]$.

\section{Definition of Morrey spaces} 
The definition of Morrey space may sound a little unnatural. Thus we now start with a few introductory words. Let $\Omega$ be an open subset of ${\mathbb{R}}^{n}$. Let $p\in[1,+\infty]$. Then we can define the extension operator $E_{\Omega}$ from ${\mathcal{L}}^{p}(\Omega)$ to ${\mathcal{L}}^{p}({\mathbb{R}}^{n})$ by setting
\[
E_{\Omega}f(x)\equiv
\left\{
\begin{array}{ll}
f(x)   &  \forall x\in \Omega\,,
 \\
0  &   \forall x\in  {\mathbb{R}}^{n}\setminus\Omega\,. 
\end{array}
\right.
\]
Now let $p<+\infty$. Since $|f|^{p}\in {\mathcal{L}}^{1}(\Omega)$, we have $|E_{\Omega}f|^{p}\in {\mathcal{L}}^{1}({\mathbb{R}}^{n})$, and accordingly, the Lebesgue Differentiation Theorem  implies that
\[
\lim_{r\to 0} \mbox{$ ^{\_}\!\!\!\!\int$}_{ {\mathbb{B}}_{n}(x,r)}
|E_{\Omega}f(y)|^{p}\,dy=|E_{\Omega}f(x)|^{p}
\qquad{\mathrm{a.a.}}\ x\in {\mathbb{R}}^{n}\,,
\]
and accordingly that
\[
\lim_{r\to 0} 
\left(\frac{1}{m_n({\mathbb{B}}_{n}(x,r))}\int_{ {\mathbb{B}}_{n}(x,r)\cap\Omega}
| f(y)|^{p}\,dy
\right)^{1/p}=| f(x)| 
\qquad{\mathrm{a.a.}}\ x\in \Omega\,.
\]
In particular, for almost all $x\in \Omega$, there exist $\rho_{x}>0$, $c_{x}>0$ such that
\[
\left(\int_{ {\mathbb{B}}_{n}(x,r)\cap\Omega}
| f(y)|^{p}\,dy
\right)^{1/p}\leq c_{x}r^{n/p}\qquad\forall r\in]0,\rho_{x}[\,.
\]
In general, we cannot expect that such an inequality holds for a constant $c_x$ that is independent of  $x\in \Omega$ for an arbitrary $f\in {\mathcal{L}}^{p}({\mathbb{R}}^{n})$. In other words, we cannot expect that there exist  $\rho >0$, $c >0$ such that
\begin{equation}
\label{introd1}
\left(\int_{ {\mathbb{B}}_{n}(x,r)\cap\Omega}
| f(y)|^{p}\,dy
\right)^{1/p}\leq c r^{n/p}\qquad\forall r\in]0,\rho [\,.
\end{equation}
Indeed, the existence of such $\rho$ and  $c$ together with the above limiting relation would imply that
\[
|f(x)|\leq c\omega_{n}^{-1/p} \qquad{\mathrm{a.a.}}\ x\in \Omega\,,
\]
and thus that $f\in {\mathcal{L}}^{\infty}(\Omega)$.\par

The idea of Morrey was to   consider condition  (\ref{introd1}) with a  choice of $c$ and $\rho$
that is independent of  $x\in\Omega$, but with $r^{n/p}$ replaced by a power $r^{\lambda}$ for some $\lambda\in]0,n/p[$. In other words, Morrey considered uniformity in $x\in\Omega$ at the cost of replacing $r^{n/p}$ with a power with a lower exponent. 

Namely, Morrey considered the case in which  there exist $\rho$ and  $c$ such that
\begin{equation}
\label{introd2}
\left(\int_{ {\mathbb{B}}_{n}(x,r)\cap\Omega}
| f(y)|^{p}\,dy
\right)^{1/p}\leq c r^{\lambda}\qquad\forall r\in]0,\rho [\,.
\end{equation}
In case $\Omega$ has a finite measure, and if we require $f$ to satisfy a stronger summability assumption, \textit{i.e.}, if we require 
 $f\in {\mathcal{L}}^{q}(\Omega)$ 
for some $q>p$, then the H\"{o}lder inequality implies that
\begin{eqnarray}\label{introd3}
\lefteqn{
\left(\int_{ {\mathbb{B}}_{n}(x,r)\cap\Omega}
| f(y)|^{p}\,dy
\right)^{1/p}
}
\\ \nonumber
&&\qquad
\leq m_{n}( {\mathbb{B}}_{n}(x,r)\cap\Omega )^{
\frac{1}{p}-\frac{1}{q}}\|f\|_{ {\mathcal{L}}^{q}(\Omega)}
\leq \omega_{n}^{
\frac{1}{p}-\frac{1}{q}}r^{^{
\frac{n}{p}-\frac{n}{q}}}\|f\|_{ {\mathcal{L}}^{q}(\Omega)}\,,
\end{eqnarray}
for all $r\in]0,+\infty[$. In other words, a stronger summability property of $f$ implies the existence of $c$ as in (\ref{introd2}) with $\lambda =(\frac{n}{p}-\frac{n}{q}) >0$ for all $ r\in]0,+\infty [$. 

Now condition (\ref{introd2})  has revealed to be important in a huge number of applications,  especially in the regularity theory for elliptic equations, and has thus motivated the definition of Morrey space on a bounded domain as the set of (equivalence classes of) functions of $L^{p}$ for which (\ref{introd2}) holds for all $x\in\Omega$. 

Condition (\ref{introd2}) has been first  considered in a paper of Morrey \cite{Mo38}   on quasi-linear elliptic equations, although the analysis Morrey spaces has developed in the 1960’s with the work of  Y. Brudnyi, S. Campanato, J. Peetre.

We now look closely at the `Morrey' condition (\ref{introd2}) and we find convenient to set.
\[
|f|_{\rho,r^{-\lambda},p,\Omega}\equiv\sup_{(x,r)\in\Omega\times ]0,\rho[}
r^{-\lambda}\|f\|_{  L^{p}(  {\mathbb{B}}_{n}(x,r)\cap\Omega)    }
\qquad\forall \rho\in]0,+\infty]\,,
\]
for all measurable functions $f$ from $\Omega$ to
${\mathbb{R}}$ and for  all $p\in [1,+\infty]$, $\lambda\in [0,+\infty[$. Clearly, $|f|_{\rho,r^{-\lambda},p,\Omega}\in[0,+\infty]$. Then we introduce the Morrey spaces
\[
{\mathcal{M}}_{p}^{r^{-\lambda},\rho }(\Omega)\equiv
\left\{
f\in  {\mathbb{R}}^{\Omega}:\, f\ 
{\text{is\ measurable}}\,,\ 
|f|_{\rho,r^{-\lambda},p,\Omega}<+\infty
\right\}\,,
\]
for all $p\in [1,+\infty]$, $\lambda\in [0,+\infty[$. One can easily check that
\[
\left({\mathcal{M}}_{p}^{r^{-\lambda},\rho }(\Omega),
\|\cdot\|_{{\mathcal{M}}_{p}^{r^{-\lambda},\rho }(\Omega)}
\equiv
|\cdot|_{\rho,r^{-\lambda},p,\Omega}\right)
\]
is a normed space.  If $\rho_1$, $\rho_2\in ]0,+\infty]$, $\rho_1\leq\rho_2$, then we obviously have
\begin{equation}\label{eq:rhomorincl}
{\mathcal{M}}_{p}^{r^{-\lambda},\rho_2 }(\Omega)
\leq
{\mathcal{M}}_{p}^{r^{-\lambda},\rho_1 }(\Omega)
\end{equation}
and the corresponding embedding is continuous. We now prove that as long as both $\rho_1$, $\rho_2$ are finite, the above inclusion is actually an equality. Indeed, the following holds.
\begin{prop}\label{prop:moud}
 Let $\Omega$ be an open subset of ${\mathbb{R}}^n$. Let $\lambda\in[0,+\infty[$, $p\in [1,+\infty]$. If 
$\rho_1$, $\rho_2\in ]0,+\infty[$, then
\[
{\mathcal{M}}^{r^{-\lambda},\rho_1}_p(\Omega)={\mathcal{M}}^{r^{-\lambda},\rho_2}_p(\Omega)\,,
\]
and the corresponding norms are equivalent. 
\end{prop}
{\bf Proof.} We can clearly assume that $\rho_1<\rho_2$. By the above inclusion (\ref{eq:rhomorincl}), it suffices to estimate 
$
|f|_{\rho_2,r^{-\lambda},p,\Omega}$
in terms of $ |f|_{\rho_1,r^{-\lambda},p,\Omega}$ for all
$ f\in {\mathcal{M}}^{r^{-\lambda},\rho_1}_p(\Omega)$.
Let $f\in {\mathcal{M}}^{r^{-\lambda},\rho_1}_p(\Omega)$.
\begin{eqnarray}\label{prop:moud1}
\lefteqn{
|f|_{\rho_2,r^{-\lambda},p,\Omega}=\sup_{(x,r)\in\Omega\times ]0,\rho_2[}
r^{-\lambda}\|f\|_{  L^{p}(  {\mathbb{B}}_{n}(x,r)\cap\Omega)    }
}
\\ \nonumber
&& 
\leq\max\left\{
\sup_{(x,r)\in\Omega\times ]0,\rho_1[}
r^{-\lambda}\|f\|_{  L^{p}(  {\mathbb{B}}_{n}(x,r)\cap\Omega)    }
,
\sup_{(x,r)\in\Omega\times [\rho_1,\rho_2[}
r^{-\lambda}\|f\|_{  L^{p}(  {\mathbb{B}}_{n}(x,r)\cap\Omega)    }
\right\}
\\ \nonumber
&& 
=\max\left\{
|f|_{\rho_1,r^{-\lambda},p,\Omega},
\sup_{(x,r)\in\Omega\times [\rho_1,\rho_2[}
r^{-\lambda}\|f\|_{  L^{p}(  {\mathbb{B}}_{n}(x,r)\cap\Omega)    }
\right\}\,.
\end{eqnarray}
Next we note that
\begin{eqnarray}\label{prop:moud2}
\lefteqn{
\sup_{(x,r)\in\Omega\times [\rho_1,\rho_2[}
r^{-\lambda}\|f\|_{  L^{p}(  {\mathbb{B}}_{n}(x,r)\cap\Omega)    }
\leq\rho_1^{-\lambda}
\sup_{(x,r)\in\Omega\times [\rho_1,\rho_2[}\|f\|_{  L^{p}(  {\mathbb{B}}_{n}(x,r)\cap\Omega)    }
}
\\ \nonumber
&&\qquad\qquad\qquad\qquad
\leq\rho_1^{-\lambda}
\sup_{ x \in\Omega }\|f\|_{  L^{p}(  {\mathbb{B}}_{n}(x,\rho_2)\cap\Omega)    }
=\rho_1^{-\lambda}
\sup_{ x \in\Omega }\|E_\Omega f\|_{  L^{p}(  {\mathbb{B}}_{n}(x,\rho_2) )    }\,.
\end{eqnarray}
Next we choose an arbitrary number in $]0,\rho_1[$, for example $\rho_1/2$, and we try to estimate
\[
\sup_{ x \in\Omega }\|E_\Omega f\|_{  L^{p}(  {\mathbb{B}}_{n}(x,\rho_2) )    }
\]
in terms of $\sup_{ z \in\Omega }\|E_\Omega f\|_{  L^{p}(  {\mathbb{B}}_{n}(z,\rho_1/2) )    }$. Indeed $\sup_{ z \in\Omega }\|E_\Omega f\|_{  L^{p}(  {\mathbb{B}}_{n}(z,\rho_1/2) )    }$ can be estimated as follows in terms of $ |f|_{\rho_1,r^{-\lambda},p,\Omega}$
\begin{eqnarray}\label{prop:moud2a}
\lefteqn{
\sup_{ z \in\Omega }\|E_\Omega f\|_{  L^{p}(  {\mathbb{B}}_{n}(z,\rho_1/2) )    }
=\sup_{ z \in\Omega }\|  f\|_{  L^{p}(  {\mathbb{B}}_{n}(z,\rho_1/2) \cap\Omega)    }
}
\\ \nonumber
&&\qquad
\leq
\sup_{ (z,r) \in\Omega\times]0,\rho_1[ }\|  f\|_{  L^{p}(  {\mathbb{B}}_{n}(z,r) \cap\Omega)    }
=
\sup_{ (z,r) \in\Omega\times]0,\rho_1[ }r^\lambda r^{-\lambda}\|  f\|_{  L^{p}(  {\mathbb{B}}_{n}(z,r) )    }
\\ \nonumber
&&\qquad
\leq
\rho_1^\lambda\sup_{ (z,r) \in\Omega\times]0,\rho_1[ }  r^{-\lambda}\|  f\|_{  L^{p}(  {\mathbb{B}}_{n}(z,r) )    }=\rho_1^\lambda  |f|_{\rho_1,r^{-\lambda},p,\Omega}\,.
\end{eqnarray}
To do so, we observe that there exist finitely many points
$\xi_1$, \dots, $\xi_{m(\rho_2)}$ of $ \overline{  {\mathbb{B}}_{n}(0,\rho_2) }$ such that
\[
\bigcup_{j=1}^{ m(\rho_2) }  {\mathbb{B}}_{n}(\xi_j,\rho_1/4)\supseteq \overline{ {\mathbb{B}}_{n}(0,\rho_2)}\,,
\]
and that accordingly
\[
\bigcup_{j=1}^{ m(\rho_2) }  {\mathbb{B}}_{n}(x+\xi_j,\rho_1/4)\supseteq \overline{ {\mathbb{B}}_{n}(x,\rho_2)}
\qquad\forall x\in {\mathbb{R}}^n\,. 
\]
If   $x\in \Omega$, we have
\begin{eqnarray}\label{prop:moud3}
\lefteqn{
\|E_\Omega f\|_{  L^{p}(  {\mathbb{B}}_{n}(x,\rho_2) )    }
}
\\ \nonumber
&&\qquad
\leq m(\rho_2)^{1/p}\sup_{j\in\{1,\dots,m(\rho_2)\}}\|E_\Omega f\|_{  L^{p}( {\mathbb{B}}_{n}(x+\xi_j,\rho_1/4) )}
\\ \nonumber
&&\qquad
= m(\rho_2)^{1/p}\sup_{j\in \{1,\dots,m(\rho_2)\},\,
{\mathbb{B}}_{n}(x+\xi_j,\rho_1/4)\cap\Omega\neq \emptyset
}\|E_\Omega f\|_{  L^{p}( {\mathbb{B}}_{n}(x+\xi_j,\rho_1/4) )}\,.
\end{eqnarray}
Next we choose a point $y_{j,x}\in {\mathbb{B}}_{n}(x+\xi_j,\rho_1/4)\cap\Omega$ for each $j\in \{1,\dots,m(\rho_2)\}$ such that 
${\mathbb{B}}_{n}(x+\xi_j,\rho_1/4)\cap\Omega\neq \emptyset$. Then we have
\[
{\mathbb{B}}_{n}(x+\xi_j,\rho_1/4)\cap\Omega\subseteq {\mathbb{B}}_{n}(y_{j,x}, \rho_1/2)
\]
for each $j\in \{1,\dots,m(\rho_2)\}$ such that 
${\mathbb{B}}_{n}(x+\xi_j,\rho_1/4)\cap\Omega\neq \emptyset$. Hence, 
\begin{eqnarray*}
\lefteqn{
m(\rho_2)^{1/p}\sup_{j\in \{1,\dots,m(\rho_2)\},\,
{\mathbb{B}}_{n}(x+\xi_j,\rho_1/4)\cap\Omega\neq \emptyset
}\|E_\Omega f\|_{  L^{p}( {\mathbb{B}}_{n}(x+\xi_j,\rho_1/4)) }
}
\\ \nonumber
&&\qquad
\leq
m(\rho_2)^{1/p}\sup_{j\in \{1,\dots,m(\rho_2)\},\,
{\mathbb{B}}_{n}(x+\xi_j,\rho_1/4)\cap\Omega\neq \emptyset
}\|E_\Omega f\|_{  L^{p}( {\mathbb{B}}_{n}(y_{j,x},\rho_1/2)) }
\\ \nonumber
&&\qquad
\leq
m(\rho_2)^{1/p}\sup_{z\in \Omega }\|E_\Omega f\|_{  L^{p}( {\mathbb{B}}_{n}(z,\rho_1/2)) }
\end{eqnarray*}
Then inequalities (\ref{prop:moud1})--(\ref{prop:moud3}) imply that
\begin{eqnarray*}
\lefteqn{
|f|_{\rho_2,r^{-\lambda},p,\Omega}
\leq
\max\left\{
|f|_{\rho_1,r^{-\lambda},p,\Omega},
\sup_{(x,r)\in\Omega\times [\rho_1,\rho_2[}
r^{-\lambda}\|f\|_{  L^{p}(  {\mathbb{B}}_{n}(x,r)\cap\Omega)    }
\right\}
}
\\ \nonumber
&&\qquad\qquad\qquad 
\leq\max\left\{
|f|_{\rho_1,r^{-\lambda},p,\Omega},
\rho_1^{-\lambda}
\sup_{ x \in\Omega }\|E_\Omega f\|_{  L^{p}(  {\mathbb{B}}_{n}(x,\rho_2) )    }
\right\}
\\ \nonumber
&&\qquad\qquad\qquad 
\leq\max\left\{
|f|_{\rho_1,r^{-\lambda},p,\Omega},
\rho_1^{-\lambda}
m(\rho_2)^{1/p}\sup_{z\in \Omega }\|E_\Omega f\|_{  L^{p}( {\mathbb{B}}_{n}(z,\rho_1/2)) }
\right\}
\\ \nonumber
&&\qquad\qquad\qquad 
\leq
\max\left\{
|f|_{\rho_1,r^{-\lambda},p,\Omega},
\rho_1^{-\lambda}m(\rho_2)^{1/p}\rho_1^{\lambda}|f|_{\rho_1,r^{-\lambda},p,\Omega}
\right\}
 \\ \nonumber
&&\qquad\qquad\qquad 
=\max\left\{
1,
 m(\rho_2)^{1/p} 
\right\}|f|_{\rho_1,r^{-\lambda},p,\Omega}
\,
\end{eqnarray*}
and thus the proof is complete.\hfill  $\Box$ 

\vspace{\baselineskip}

The previous proposition shows that the spaces ${\mathcal{M}}^{r^{-\lambda},\rho}_p(\Omega)$ coincide for all finite values of $\rho$.

\vspace{\baselineskip}

If $\rho$ is finite, the definition of $|f|_{\rho,r^{-\lambda},p,\Omega}$ requires information on $f$ on the balls $ {\mathbb{B}}_{n}(x,r)\cap\Omega$ for $r$ only in the right neighborhood $]0,\rho[$, and we can take $\rho$ as small as we wish,  and thus the membership of $f$ in ${\mathcal{M}}^{r^{-\lambda},\rho}_p(\Omega)$ can be interpreted as a regularity condition of $f$.

\vspace{\baselineskip}

If instead $\rho=+\infty$, the definition of $|f|_{+\infty,r^{-\lambda},p,\Omega}$ involves information on $f$ on the balls $ {\mathbb{B}}_{n}(x,r)\cap\Omega$   both  when $r$ is small and when $r$ is arbitrarily large and thus the membership of $f$ in ${\mathcal{M}}^{r^{-\lambda},+\infty}_p(\Omega)$ can be interpreted as a condition that is at one hand a regularity condition and on the other hand also a condition on the behavior of $f$ at infinity.

\vspace{\baselineskip}

In order to treat  ${\mathcal{M}}^{r^{-\lambda},\rho}_p(\Omega)$ both in case $\rho<+\infty$ and $\rho=+\infty$ at the same time, we resort to the generalized Morrey spaces, that we now introduce by means of the following. 
\begin{defn}
\label{morrey}
Let $\Omega$ be an open subset of ${\mathbb{R}}^{n}$. Let $p\in [1,+\infty]$. Let $w$ be a function from $]0,+\infty[$ to $[0,+\infty[$. Let
\[
|f|_{\rho,w,p,\Omega}\equiv\sup_{(x,r)\in\Omega\times ]0,\rho[}
w (r)\|f\|_{  L^{p}(  {\mathbb{B}}_{n}(x,r)\cap\Omega)    }
\qquad\forall \rho\in]0,+\infty]\,,
\]
for all measurable functions $f$ from $\Omega$ to
${\mathbb{R}}$.

Assume that there exists $r_0\in]0,+\infty[$ such that $w(r_0)\neq 0$. Then we define as generalized Morrey space with weight $w $ and exponent $p$ the set
\[
{\mathcal{M}}_{p}^{w }(\Omega)\equiv
\left\{
f\in  {\mathbb{R}}^{\Omega}:\, f\ 
{\text{is\ measurable}}\,,\ 
|f|_{+\infty,w,p,\Omega}<+\infty
\right\}\,.
\]
Then we set
\[
\|f\|_{  {\mathcal{M}}_{p}^{w }(\Omega)  }\equiv |f|_{+\infty,w,p,\Omega}
\qquad\forall f\in {\mathcal{M}}_{p}^{w }(\Omega)\,.
\]
\end{defn}
One can easily verify that 
$({\mathcal{M}}_{p}^{w }(\Omega), \|\cdot\|_{  {\mathcal{M}}_{p}^{w }(\Omega)  })$ is a
 normed space.
 
 Here we mention that the definition of generalized Morrey spaces	 
   is not uniform in the literature.  So for example if $\Omega={\mathbb{R}}^n$, the definition here coincides with that of Gogatishvili and Mustafayev \cite{GoMu11}. Then   we can obtain  the definition of   Nakai \cite{Na00}  by taking $w(r)=\varphi(r)^{-1}m_n({\mathbb{B}}_n(0,r))^{-1/p}$ and that of
Sawano \cite{Sa19} by taking $w(r)=\varphi(r) m_n({\mathbb{B}}_n(0,r))^{-1/p}$ for some real valued function $\varphi$. We also note that here we say that $w$ is a `weight' (as in Gogatishvili and Mustafayev \cite{GoMu11}), while other authors reserve the word weight for a weight put on the measure in $\Omega$ as in Samko \cite{Sam08} (a case that we do not discuss in these notes).

One may wonder whether   $|f|_{\rho,w,p,\Omega}$ would change if we replace the supremum in $x\in \Omega$ with $x\in \overline{\Omega}$. The answer is no as the following lemma shows.
 \begin{lem}
\label{lem:mocls}
Let $\Omega$ be an open subset of ${\mathbb{R}}^{n}$. Let $p\in [1,+\infty]$. Let $w$ be a   function from $]0,+\infty[$ to $[0,+\infty[$.    If $f$ is a measurable function  from $\Omega$ to ${\mathbb{R}}$ such that  $f\in L^{p}({\mathbb{B}}_{n}(x,r)\cap\Omega)$ for all $(x,r)\in \Omega\times ]0,+\infty[$, then 
\[
\sup_{(x,r)\in\Omega\times ]0,\rho[}
w (r)\|f\|_{  L^{p}(  {\mathbb{B}}_{n}(x,r)\cap\Omega)    }
=
\sup_{(x,r)\in\overline{\Omega}\times ]0,\rho[}
w (r)\|f\|_{  L^{p}(  {\mathbb{B}}_{n}(x,r)\cap\Omega)    }
\,,
\]
for all $\rho\in]0,+\infty]$.
\end{lem}
{\bf Proof.} It clearly suffices to prove that
\[
\sup_{(x,r)\in\partial{\Omega}\times ]0,\rho[}
w (r)\|f\|_{  L^{p}(  {\mathbb{B}}_{n}(x,r)\cap\Omega)    }
\leq
\sup_{(x,r)\in\Omega\times ]0,\rho[}
w (r)\|f\|_{  L^{p}(  {\mathbb{B}}_{n}(x,r)\cap\Omega)    }\,,
\]
for all $\rho\in]0,+\infty]$. Let $\rho\in]0,+\infty]$. Let $(\tilde{x},r)\in (\partial{\Omega})\times]0,\rho[$. Let $\{x_{j}\}_{j\in {\mathbb{N}} }$
be a sequence in $\Omega$ such that $\lim_{j\to\infty}x_{j}=\tilde{x}$.
It clearly suffices to show that
\begin{equation}
\label{prelprgm6}
w (r)\|f\|_{ L^{p}({\mathbb{B}}_{n}(\tilde{x},r)\cap\Omega)}
\leq \limsup_{j\to\infty}w (r)\|f\|_{ L^{p}({\mathbb{B}}_{n}(x_{j},r)\cap\Omega)}\,.
\end{equation}
Indeed, 
\[
\limsup_{j\to\infty}w (r)\|f\|_{ L^{p}({\mathbb{B}}_{n}(x_{j},r)\cap\Omega)}
\leq
\sup_{(x,r)\in\Omega\times ]0,\rho[}
w (r)\|f\|_{  L^{p}(  {\mathbb{B}}_{n}(x,r)\cap\Omega)    }\,.
\]
Possibly neglecting a finite number of terms, we can clearly assume that 
\[
|x_{j}-\tilde{x}|\leq 1\,,
\]
 and accordingly that
\[
{\mathbb{B}}_{n}(x_{j},r)\subseteq {\mathbb{B}}_{n}(\tilde{x},r+1)
\qquad\forall j\in {\mathbb{N}}\,.
\]
Obviously,
\begin{equation}
\label{prelprgm7}
\lim_{j\to\infty}\chi_{ {\mathbb{B}}_{n}(x_{j},r)\cap\Omega}(x)f (x)
=
\chi_{ {\mathbb{B}}_{n}(\tilde{x},r)\cap\Omega}(x)f (x)
\qquad\forall x\in\Omega
\setminus
\partial{\mathbb{B}}_{n}(\tilde{x},r)
\,,
\end{equation}
and accordingly for almost all $x\in \Omega$. Indeed, $m_{n}(
\partial{\mathbb{B}}_{n}(\tilde{x},r))=0$.
Next we consider separately cases $p<\infty$ and case $p=\infty$.

Let $p<\infty$. Since
\[
{\mathbb{B}}_{n}(x_{j},r)\subseteq
 {\mathbb{B}}_{n}(\tilde{x},r+1)\subseteq
{\mathbb{B}}_{n}(x_{0}, r+2)
\qquad\forall j\in {\mathbb{N}}\,,
\]
we have
\[
| \chi_{ {\mathbb{B}}_{n}(x_{j},r)\cap\Omega} f  |^{p}
\leq
|f |^{p}\chi_{ {\mathbb{B}}_{n}(x_{0}, r+2)\cap\Omega} 
\in L^{1}(\Omega)\,,
\]
for all $j\in {\mathbb{N}}$, and the convergence of (\ref{prelprgm7}) and the Dominated Convergence Theorem imply that
\[
\lim_{j\to\infty}\int_{\Omega}| \chi_{ {\mathbb{B}}_{n}(x_{j},r)\cap\Omega} f  |^{p}\,dx
=
\int_{\Omega}| \chi_{ {\mathbb{B}}_{n}(\tilde{x},r)\cap\Omega} f  |^{p}\,dx\,,
\]
and thus (\ref{prelprgm6}) holds. 

We now assume that $p=\infty$. Let $N$ be a subset of $\Omega$ of measure zero such that
\begin{eqnarray*}
&&\partial{\mathbb{B}}_{n}(\tilde{x},r)\subseteq N\,,
\\
&&|f(x)|\leq\|f\|_{ L^{\infty}( {\mathbb{B}}_{n}(x_{0}, r+2)\cap\Omega) }\qquad\forall x\in ({\mathbb{B}}_{n}(
x_{0}, r+2)\cap\Omega)\setminus N\,,
\\
&&|f(x)|\leq\|f\|_{ L^{\infty}( {\mathbb{B}}_{n}(x_{l},r )\cap\Omega) }\qquad\forall x\in ({\mathbb{B}}_{n}(
x_{l},r)\cap\Omega)\setminus N\,,
\qquad\forall l\in {\mathbb{N}}\,,
\end{eqnarray*}
(\textit{cf.} (\ref{prelprgm7})). Then we have 
\begin{eqnarray*}
\lefteqn{
w(r)|f(x) \chi_{ {\mathbb{B}}_{n}(\tilde{x},r)\cap\Omega}(x)|
\leq
w(r)\limsup_{j\to\infty}|f(x) \chi_{ {\mathbb{B}}_{n}(x_{j},r)\cap\Omega}(x)|
}
\\
\nonumber
&& \qquad 
\leq w(r) \inf_{j\in{\mathbb{N}}}\sup_{l\geq j}
|f(x) \chi_{ {\mathbb{B}}_{n}(x_{l},r)\cap\Omega}(x)|
\\
\nonumber
&&\qquad
\leq \inf_{j\in{\mathbb{N}}}\sup_{l\geq j}w(r)\|f\|_{L^{\infty}(
{\mathbb{B}}_{n}(x_{l},r)\cap\Omega)}
\\
\nonumber
&& \qquad 
\leq  
\limsup_{j\to\infty}w(r)\|f\|_{L^{\infty}(
{\mathbb{B}}_{n}(x_{l},r)\cap\Omega)}
\qquad\forall x\in ({\mathbb{B}}_{n}(x_{0},r+2)\cap\Omega)\setminus N\,.
\end{eqnarray*}
Since $({\mathbb{B}}_{n}(\tilde{x},r)\cap\Omega)\setminus N \subseteq
({\mathbb{B}}_{n}(x_{0}, r+2)\cap\Omega)\setminus N $, 
such an inequality holds in particular for all 
$x\in ({\mathbb{B}}_{n}(\tilde{x},r)\cap\Omega)\setminus N$
and accordingly
\[
w(r)\|f\|_{L^{\infty}(
{\mathbb{B}}_{n}(\tilde{x},r)\cap\Omega)}
\leq \limsup_{j\to\infty}w(r)\|f\|_{L^{\infty}(
{\mathbb{B}}_{n}(x_{l},r)\cap\Omega)}
\,,
\]
and thus (\ref{prelprgm6}) holds.\hfill  $\Box$ 

\vspace{\baselineskip}

Next we introduce the  choices of $w$ that we consider most in the present notes. If we choose the function $w(r)\equiv r^{-\lambda}$ for all $r\in]0,+\infty[$, we have
\[
{\mathcal{M}}_{p}^{r^{-\lambda} }(\Omega)={\mathcal{M}}_{p}^{r^{-\lambda,+\infty} }(\Omega)\,.
\]
If instead we take $\rho\in]0,+\infty[$ and we consider the function
\[
w_{\lambda,\rho}(r)\equiv\left\{
\begin{array}{ll}
r^{-\lambda} & \text{if}\ r\in]0,\rho[\,,
\\
0& \text{if}\ r\in [\rho,+\infty[\,,
 \end{array}
\right.
\]
we have the space
\[
{\mathcal{M}}_{p}^{w_{\lambda,\rho} }(\Omega)={\mathcal{M}}_{p}^{r^{-\lambda},\rho }(\Omega)\,.
\]
In other words,  the space ${\mathcal{M}}^{r^{-\lambda},\rho}_p(\Omega)$   is a specific 
generalized Morrey space both in case $\rho<+\infty$ and $\rho=+\infty$.

Also important  for the analysis of the regularity of functions are the spaces
\[
{\mathcal{M}}_{p}^{r^{-\lambda},\rho }(\Omega)\cap L^p(\Omega)
\]
with the norm
\[
\|f\|_{{\mathcal{M}}_{p}^{r^{-\lambda},\rho }(\Omega)\cap L^p(\Omega)}
=
\max\left\{
\|f\|_{{\mathcal{M}}_{p}^{r^{-\lambda},\rho }(\Omega) },
\|f\|_{L^p(\Omega)}
\right\}
\]
when $\rho\in]0,+\infty[$. As the following lemma shows, it is just another case of generalized Morrey space. The following statement is a variant of a consequence of a more general statement of Sawano \cite[Example 3.4, third point]{Sa19} for generalized Morrey spaces in ${\mathbb{R}}^n$. Since here we are dealing  with generalized Morrey spaces in a domain $\Omega$  of ${\mathbb{R}}^n$, to invoke Sawano's result we would need to  prove  the existence of a linear and continuous extension operator from the generalized Morrey space  in $\Omega$ to the corresponding generalized Morrey space  in ${\mathbb{R}}^n$.   Here instead we include  a direct proof.

\begin{lem}
\label{prem}
Let $\Omega$ be an open subset of ${\mathbb{R}}^{n}$. Let $p\in [1,+\infty]$, $\lambda\in [0,+\infty[$. Let $w_{\lambda}$ be the function  from $]0,+\infty[$ to itself defined by 
\begin{equation}
\label{prem1}
w_{\lambda}(r)\equiv
 \left\{
 \begin{array}{ll}
r^{-\lambda} & {\mathrm{if}}\ r\in]0,1]
\\
1 & {\mathrm{if}}\ r\in [1,+\infty[\,.
\end{array}
 \right.
\end{equation}
If $\rho\in]0,+\infty[$, then 
\begin{equation}\label{prem1a}
{\mathcal{M}}_{p}^{r^{-\lambda},\rho }(\Omega)\cap L^p(\Omega)
={\mathcal{M}}_{p}^{w_\lambda }(\Omega)
\end{equation}
and the corresponding norms are equivalent.
 \end{lem}
{\bf Proof.} We can clearly assume that $\Omega$ is not empty. We first assume that $f\in {\mathcal{M}}_{p}^{r^{-\lambda},\rho }(\Omega)\cap L^p(\Omega)$ and we show that $f\in {\mathcal{M}}_{p}^{w_\lambda }(\Omega)$. To do so, we must estimate 
\[
\sup_{(x,r)\in\Omega\times]0,+\infty[}w_{\lambda}(r)\|f\|_{  L^{p}(  {\mathbb{B}}_{n}(x,r)\cap\Omega)    }\,.
\]
 We do so by splitting case $r\in ]0,\min\{\rho,1\}[$ where $w_\lambda(r)=r^{-\lambda}$ and case $r\in [\min\{\rho,1\},+\infty[$, where 
$w_\lambda(r)\leq w_{\lambda}(\min\{\rho,1\})$. Thus we have
\begin{eqnarray*}
\lefteqn{
\sup_{(x,r)\in\Omega\times]0,\min\{\rho,1\}[}w_{\lambda}(r)\|f\|_{  L^{p}(  {\mathbb{B}}_{n}(x,r)\cap\Omega)    }
=
\sup_{(x,r)\in\Omega\times]0,\min\{\rho,1\}[}r^{-\lambda}\|f\|_{  L^{p}(  {\mathbb{B}}_{n}(x,r)\cap\Omega)    }
}
\\ \nonumber
&&\qquad\qquad\qquad\qquad\qquad
\leq 
\sup_{(x,r)\in\Omega\times]0,\rho[}r^{-\lambda}\|f\|_{  L^{p}(  {\mathbb{B}}_{n}(x,r)\cap\Omega)    }
\leq  \|f\|_{{\mathcal{M}}_{p}^{r^{-\lambda},\rho }(\Omega) }
\,,
\end{eqnarray*}
and
\[
\sup_{(x,r)\in\Omega\times [\min\{\rho,1\},+\infty[}w_{\lambda}(r)\|f\|_{  L^{p}(  {\mathbb{B}}_{n}(x,r)\cap\Omega)    }\leq
w_{\lambda}(\min\{\rho,1\})\|f\|_{  L^{p}(   \Omega)    }
\,.
\]
Hence, $f\in {\mathcal{M}}_{p}^{w_\lambda }(\Omega)$ and 
\[
\|f\|_{{\mathcal{M}}_{p}^{w_\lambda }(\Omega)}\leq
w_{\lambda}(\min\{\rho,1\})
 \|f\|_{{\mathcal{M}}_{p}^{r^{-\lambda},\rho }(\Omega)\cap L^p(\Omega)}
\]
 Conversely, we now assume that  $f\in {\mathcal{M}}_{p}^{w_\lambda }(\Omega)$. We first show that $f\in {\mathcal{M}}_{p}^{r^{-\lambda},\rho }(\Omega) $. Since $r^{-\lambda}\leq w_\lambda(r)$ for all $r\in]0,+\infty[$, we have
\begin{eqnarray*}
\lefteqn{
 \sup_{(x,r)\in\Omega\times]0,\rho[}r^{-\lambda}\|f\|_{  L^{p}(  {\mathbb{B}}_{n}(x,r)\cap\Omega)    }
 \leq
 \sup_{(x,r)\in\Omega\times]0,\rho[}w_\lambda(r)\|f\|_{  L^{p}(  {\mathbb{B}}_{n}(x,r)\cap\Omega)    }
 }
 \\ \nonumber
&&\qquad\qquad\qquad\qquad\qquad
 \leq 
 \sup_{(x,r)\in\Omega\times]0,+\infty[}w_\lambda(r)\|f\|_{  L^{p}(  {\mathbb{B}}_{n}(x,r)\cap\Omega)    }
 =
 \|f\|_{{\mathcal{M}}_{p}^{w_\lambda }(\Omega)}\,.
\end{eqnarray*}
 Hence, $f\in {\mathcal{M}}_{p}^{r^{-\lambda},\rho }(\Omega)$.  We now show that $f\in L^{p}(\Omega)$. 
If $p<\infty$, then  we fix $\xi\in \Omega$.  Then inequality
\[
 \int_{ {\mathbb{B}}_{n}(\xi,s)\cap\Omega}
| f(y)|^{p}\,dy
 \leq w_{\lambda}(s)^{-p}\left(\sup_{(x,r)\in\Omega\times ]0,+\infty[}w_{\lambda}(r)\|f\|_{  L^{p}(  {\mathbb{B}}_{n}(x,r)\cap\Omega)    }\right)^{p}\,,
\]
for all $s\in [1,+\infty[$ 
together with the Monotone Convergence Theorem imply that
\begin{equation}
\label{prem2}
\int_{  \Omega}
| f(y)|^{p}\,dy\leq
 \left(\|f\|_{ {\mathcal{M}}_{p}^{w_\lambda }(\Omega)}\right)^{p}\,.
\end{equation}
Indeed, $w_{\lambda}(s)^{-p}=1$ for all $s\in [1,+\infty[$.
Hence,
$\int_{  \Omega}
| f(y)|^{p}\,dy<\infty$ and  accordingly $f\in L^{p}(\Omega)$.

If $p=+\infty$, then  we fix $\xi\in \Omega$. Then the inequality
\[
{\mathrm{ess\,sup}}_{{\mathbb{B}}_{n}(\xi,s)\cap\Omega}|f|
\leq w_{\lambda}(s)^{-1} \sup_{(x,r)\in\Omega\times ]0,+\infty[}w_{\lambda}(r)\|f\|_{  L^{\infty}(  {\mathbb{B}}_{n}(x,r)\cap\Omega)    } 
\quad\forall s\in [1,+\infty[\,,
\]
implies that 
\begin{equation}
\label{prem3}
{\mathrm{ess\,sup}}_{ \Omega}|f|
\leq  \|f\|_{ {\mathcal{M}}_{\infty}^{w_\lambda }(\Omega)} \,.
\end{equation}
Indeed, $w_{\lambda}(s)^{-p}=1$ for all $s\in [1,+\infty[$.
Hence, $f\in L^{\infty}(\Omega)$. Moreover, the above inequalities imply that
\[
 \|f\|_{{\mathcal{M}}_{p}^{r^{-\lambda},\rho }(\Omega)\cap L^p(\Omega)}
 \leq
 \|f\|_{ {\mathcal{M}}_{p}^{w_\lambda }(\Omega)}\,,
\]
both in case $p<+\infty$ and in case $p=+\infty$. 
\hfill  $\Box$

\vspace{\baselineskip}

If  $p\in [1,+\infty]$ and $\lambda\in[0,+\infty[$, we find convenient to set
\begin{equation}\label{morreym}
M_{p}^{\lambda }(\Omega)\equiv  {\mathcal{M}}_{p}^{w_{\lambda} }(\Omega)\,,
\end{equation}
and 
\[
\|f\|_{  M_{p}^{\lambda }(\Omega)  }\equiv
\|f\|_{  {\mathcal{M}}_{p}^{w_{\lambda} }(\Omega)  }
\qquad\forall f\in M_{p}^{\lambda }(\Omega)\,.
\]
 
\begin{corol}\label{corol:prembd}
 Let $\Omega$ be a bounded open subset of ${\mathbb{R}}^{n}$. Let $p\in [1,+\infty]$, $\lambda\in [0,+\infty[$. If $\rho\in]0,+\infty[$, then 
 \begin{equation}\label{corol:prembd1}
{\mathcal{M}}_{p}^{r^{-\lambda},\rho }(\Omega) 
={\mathcal{M}}_{p}^{w_\lambda }(\Omega)
\end{equation}
 and the corresponding norms are equivalent.
 \end{corol}
{\bf Proof.} Let $\rho_1>{\mathrm{diam}}(\Omega)$, $r\in]0,{\mathrm{diam}}(\Omega)[$. Then we have
\begin{eqnarray*}
\lefteqn{
\|g\|_{L^p(\Omega)}=\|g\|_{L^p(\Omega\cap{\mathbb{B}}_n(x,r))}
\leq w_{\lambda,\rho_1}(r)^{-1}w_{\lambda,\rho_1}(r)\|g\|_{L^p(\Omega\cap{\mathbb{B}}_n(x,r))}
}
\\ \nonumber
&&\quad
\leq w_{\lambda,\rho_1}(r)^{-1}\|g\|_{{\mathcal{M}}_{p}^{r^{-\lambda},\rho_1 }(\Omega) }
=(r^{-\lambda})^{-1}\|g\|_{{\mathcal{M}}_{p}^{r^{-\lambda},\rho_1 }(\Omega) }
\qquad\forall g\in {\mathcal{M}}_{p}^{r^{-\lambda},\rho_1 }(\Omega) 
\end{eqnarray*}
and accordingly
\[
{\mathcal{M}}_{p}^{r^{-\lambda},\rho_1 }(\Omega)\cap L^p(\Omega)={\mathcal{M}}_{p}^{r^{-\lambda},\rho_1 }(\Omega)
\]
and the corresponding norms are equivalent. Hence, Proposition \ref{prop:moud} and Lemma \ref{prem} imply that
\[
{\mathcal{M}}_{p}^{r^{-\lambda},\rho  }(\Omega)={\mathcal{M}}_{p}^{r^{-\lambda},\rho_1 }(\Omega) 
={\mathcal{M}}_{p}^{r^{-\lambda},\rho_1 }(\Omega)\cap L^p(\Omega)={\mathcal{M}}_{p}^{w_\lambda }(\Omega)
\]
and the corresponding norms are equivalent.
\hfill  $\Box$ 

\vspace{\baselineskip}

\begin{rem}{\em  In case $\Omega={\mathbb{R}}^n$, the norm of ${\mathcal{M}}_{p}^{r^{-\lambda}}({\mathbb{R}}^n)$ satisfies an important homogeneity property. Indeed, an elementary computation shows that
  \[
  \|f(\alpha \cdot)\|_{{\mathcal{M}}_{p}^{r^{-\lambda}}({\mathbb{R}}^n)}
  =\alpha^{\lambda-(n/p)}
  \|f(  \cdot)\|_{{\mathcal{M}}_{p}^{r^{-\lambda}}({\mathbb{R}}^n)}
  \]
  for all $f\in {\mathcal{M}}_{p}^{r^{-\lambda}}({\mathbb{R}}^n)$ and $\alpha\in]0,+\infty[$. Hence, ${\mathcal{M}}_{p}^{r^{-\lambda}}({\mathbb{R}}^n)$ is said to be the `homogeneous' Morrey space of exponents $\lambda$, $p$.
 }
\end{rem}

\begin{rem}{\em
N.~Samko \cite{Sa09} has introduced a generalization of Morrey spaces by replacing the Lebesgue measure by a weighted measure, that we do not consider here.}
\end{rem}

\section{Vanishing generalized Morrey spaces}
 We now define an important subspace of a (generalized) Morrey space: the `vanishing (generalized) Morrey space' or `little (generalized) Morrey space'. 
\begin{defn}
\label{lmorrey}
Let $\Omega$ be an open subset of ${\mathbb{R}}^{n}$. Let $p\in [1,+\infty]$. 
  Let $w$ be a function from $[0,+\infty[$ to $]0,+\infty[$.  Assume that there exists $r_0\in]0,+\infty[$ such that $w(r_0)\neq 0$.

Then we define as generalized vanishing (or little)  Morrey space with weight $w $ and exponent $p$ the subspace
\[
{\mathcal{M}}_{p}^{w ,0}(\Omega)\equiv
\left\{
f\in {\mathcal{M}}_{p}^{w }(\Omega):\,\lim_{\rho\to 0}
|f|_{\rho,w,p,\Omega}=0
\right\}
\]
of ${\mathcal{M}}_{p}^{w }(\Omega)$.
 
\end{defn}
We always think the space   ${\mathcal{M}}_{p}^{w ,0}(\Omega)$ as endowed of the norm  of
${\mathcal{M}}_{p}^{w }(\Omega)$. In particular, we have the corresponding `vanishing' variants 
${\mathcal{M}}_{p}^{r^{-\lambda},0}(\Omega)$, ${\mathcal{M}}_{p}^{r^{-\lambda} ,\rho,0}(\Omega)$ for $\rho\in]0,+\infty[$ and $M_p^{\lambda,0}(\Omega)$
for the specific cases  ${\mathcal{M}}_{p}^{r^{-\lambda}}(\Omega)$, ${\mathcal{M}}_{p}^{r^{-\lambda} ,\rho}(\Omega)$ for $\rho\in]0,+\infty[$ and $M_p^\lambda(\Omega)$
that are associated to the weights $r^{-\lambda}$, $w_{\lambda,\rho}$  for $\rho\in]0,+\infty[$, 
and $w_{\lambda}$.

\section{Trivial and nontrivial generalized Morrey spaces}

In this section we briefly discuss the issue of triviality of the generalized Morrey spaces and we first prove the following triviality condition for a generalized Morrey space (and for the corresponding vanishing space).
\begin{prop}\label{prop:motri}
 Let $\Omega$ be an open subset of ${\mathbb{R}}^{n}$. Let $p\in [1,+\infty]$. 
  Let $w$ be a function from $]0,+\infty[$ to $[0,+\infty[$. If there exists $r_0\in]0,+\infty[$ such that 
  $w(r)\neq 0$ for all $r\in ]0,r_0[$, and if 
\begin{equation}\label{prop:motri1}
 \limsup_{r\to 0}w(r)r^{n/p}=+\infty\,,
\end{equation}
  then ${\mathcal{M}}_{p}^{w}(\Omega)=\{0\}$. 
\end{prop}
{\bf Proof.} Let $f\in  {\mathcal{M}}_{p}^{w}(\Omega)$. Since ${\mathcal{M}}_{p}^{w}(\Omega)\subseteq  L^p({\mathbb{B}}_n(0,r)\cap\Omega)$ for all $r\in]0,+\infty[$, we have $f\in 
 L^p_{{\mathrm{loc}}}(\Omega)\subseteq L^1_{{\mathrm{loc}}}(\Omega)$ and  the Lebesgue Dfferentiation Theorem implies that
\[
\lim_{r\to 0} \Xint-_{
{\mathbb{B}}_n(x,r)
}|f(y)|\,dy=|f(x)|\qquad {\mathrm{a.a.}}\ x\in \Omega\,.
\]
So we now try to show that our assumptions imply that the limit in the right hand side equals zero. To do so,  we note that the H\"{o}lder inequality implies that
\begin{eqnarray*}
\lefteqn{
\Xint-_{
{\mathbb{B}}_n(x,r)
}|f(y)|\,dy=m_n({\mathbb{B}}_n(x,r))^{-1}\int_{ {\mathbb{B}}_n(x,r) }|f(y)|\,dy
}
\\ \nonumber
&&\qquad
\leq m_n({\mathbb{B}}_n(x,r))^{-1}m_n({\mathbb{B}}_n(x,r))^{\frac{1}{1}-\frac{1}{p}}
\|f\|_{L^p({\mathbb{B}}_n(x,r) )}
\\ \nonumber
&&\qquad
=m_n({\mathbb{B}}_n(x,r))^{ -\frac{1}{p}}w(r)^{-1}w(r) 
\|f\|_{L^p({\mathbb{B}}_n(x,r)\cap\Omega)}
\\ \nonumber
&&\qquad
\leq\frac{1}{
\omega_n^{ 1/p}r^{n/p}w(r)}\|f\|_{{\mathcal{M}}_{p}^{w}(\Omega)}\,,
\end{eqnarray*}
for all $x\in \Omega$ and $r\in]0,r_0[$ such that ${\mathbb{B}}_n(x,r)\subseteq \Omega$.
By assumption (\ref{prop:motri1}), we have
\[
 \liminf_{r\to 0}\frac{1}{
 r^{n/p}w(r)}=0
\]
and thus we have
\begin{eqnarray*}
\lefteqn{
|f(x)|=\lim_{r\to 0} \Xint-_{
{\mathbb{B}}_n(x,r)
}|f(y)|\,dy = \liminf_{r\to 0} \Xint-_{
{\mathbb{B}}_n(x,r)
}|f(y)|\,dy 
}
\\ \nonumber
&&\qquad\qquad
\leq
\liminf_{r\to 0}
\frac{1}{
\omega_n^{ 1/p}r^{n/p}w(r)}\|f\|_{{\mathcal{M}}_{p}^{w}(\Omega)}=0
\qquad {\mathrm{a.a.}}\ x\in \Omega\,,
\end{eqnarray*}
and thus the proof is complete. \hfill  $\Box$ 

\vspace{\baselineskip}

By applying the previous proposition to the special weights $w_{\lambda,\rho}$, $r^{-\lambda}$, $w_{\lambda}$, we deduce the validity of the following  corollary. 
\begin{corol}\label{corol:motri}
 Let $\Omega$ be an open subset of ${\mathbb{R}}^{n}$. Let $p\in [1,+\infty]$. 
 Let $\lambda\in [0,+\infty[$, $\lambda>n/p$, then  ${\mathcal{M}}_{p}^{r^{-\lambda},\rho}(\Omega)=\{0\}$ for all $\rho\in]0,+\infty[$, ${\mathcal{M}}_{p}^{r^{-\lambda}}(\Omega)=\{0\}$, $M^{\lambda}_p(\Omega)=\{0\}$. 
\end{corol}
Corollary \ref{prop:motri} says in particular that the Morrey spaces ${\mathcal{M}}_{p}^{r^{-\lambda},\rho}(\Omega)$ for all $\rho\in]0,+\infty[$, ${\mathcal{M}}_{p}^{r^{-\lambda}}(\Omega)$, $M^{\lambda}_p(\Omega)$ or the corresponding vanishing variants can be nontrivial only if $\lambda\in [0,n/p]$.

Then one can wonder whether there are sufficient conditions for nontriviality of a generalized Morrey space. Thus we introduce the following  as in Burenkov V.I., Guliyev H.V.~\cite{BuGu}, Burenkov V.I., Jain P., Tararykova T.V.~\cite{BuJaTa}, Burenkov V.I.,  Liflyand E.~\cite{BuLi20}.
 \begin{defn}
\label{opinf}
Let $p\in [1,+\infty]$. We denote by $W_{p,\infty}$ the set of measurable functions $w$ from $]0,+\infty[$ to itself such that there exist $t_{1},t_{2}\in ]0,+\infty[$ such that
\[
{\mathrm{ess\,sup}}_{r\in]t_{1},+\infty[}w(r)<+\infty\,,
\qquad
{\mathrm{ess\,sup}}_{r\in]0, t_{2} [}w(r)r^{n/p}<+\infty\,.
\]
\end{defn}
Then we have the following result of Burenkov V.I., Guliyev H.V.~\cite{BuGu}, Burenkov V.I., Jain P., Tararykova T.V.~\cite{BuJaTa}, which we do not prove here.
\begin{thm}
Let $p\in [1,+\infty]$. Let $\Omega$ be an open subset of ${\mathbb{R}}^{n}$.  If $w \in W_{p,\infty}$, then 
${\mathcal{M}}_{p}^{w }(\Omega)\neq \{0\}$. 
\end{thm}

\section{Elementary embeddings  of generalized Morrey spaces into Lebesgue spaces}
We first  prove that the functions of ${\mathcal{M}}_{p}^{w}(\Omega)$ are locally $p$-summable.
\begin{prop}
\label{lpim}
Let $\Omega$ be an open subset of ${\mathbb{R}}^{n}$. Let $p\in [1,+\infty]$.   Let $w$ be a   function from $]0,+\infty[$ to $[0,+\infty[$.  
Then the following two statements hold.
\begin{enumerate}
\item[(i)] Let $\xi\in \Omega$, $R\in]0,+\infty[$, $w(R)\neq 0$. Then the map from 
\[
{\mathcal{M}}_{p}^{w}(\Omega)\qquad\text{to}\qquad
 L^{p}({\mathbb{B}}_{n}(\xi,R)\cap \Omega)
 \]
 that takes $f\in 
{\mathcal{M}}_{p}^{w}(\Omega)$ to the restriction $f_{|{\mathbb{B}}_{n}(\xi,R)\cap \Omega}$  is linear and continuous. In particular,
\[
\|f\|_{L^{p}({\mathbb{B}}_{n}(\xi,R)\cap \Omega)}
\leq
w(R)^{-1}\|f\|_{{\mathcal{M}}_{p}^{w}(\Omega)}
\qquad\forall f\in {\mathcal{M}}_{p}^{w}(\Omega)\,.
\]
\item[(ii)] Assume that there exists $r_0\in]0,+\infty[$ such that $w(r_0)\neq 0$.
Then ${\mathcal{M}}_{p}^{w}(\Omega)$ is continuously embedded into $L^{p}_{{\mathrm{loc}} }(\Omega)$.
\end{enumerate}
\end{prop}
{\bf Proof.} We can clearly assume that $\Omega$ is not empty.
Statement (i) is an immediate consequence of the inequality
\begin{eqnarray*}
\lefteqn{ 
\|f\|_{L^{p}({\mathbb{B}}_{n}(\xi,R)\cap \Omega)}=
w(R)^{-1}w(R)\|f\|_{L^{p}({\mathbb{B}}_{n}(\xi,R)\cap \Omega)}
}
\\
\nonumber
&& \qquad
\leq w(R)^{-1}
\sup_{(x,r)\in \Omega\times]0,+\infty[}
w(r)\|f\|_{L^{p}({\mathbb{B}}_{n}(x,r)\cap \Omega)}
=w(R)^{-1}\|f\|_{{\mathcal{M}}_{p}^{w}(\Omega)}\,,
\end{eqnarray*}
for all $f\in {\mathcal{M}}_{p}^{w}(\Omega)$. 

We now consider statement (ii).  Let $K$ be a compact subset of $\Omega$. Since $K$ is compact,  there exist $\xi_1$, \dots, $\xi_m$ such that
\[
K\subseteq \bigcup_{j=1}^m{\mathbb{B}}_{n}(\xi_j,r_0)\cap \Omega
\]
Then the inequality of (i) implies that
\begin{eqnarray*}
\lefteqn{
\|f\|_{L^{p}(K)}
\leq 
m^{1/p}\sup_{j=1,\dots,m}\|f\|_{L^{p}({\mathbb{B}}_{n}(\xi_j,r_0)\cap \Omega)}
}
\\ \nonumber
&&\qquad\qquad
\leq 
m^{1/p}w(r_0)^{-1}\|f\|_{{\mathcal{M}}_{p}^{w}(\Omega)}
\qquad\forall  f\in {\mathcal{M}}_{p}^{w}(\Omega)\,.
\end{eqnarray*}
\hfill  $\Box$

\vspace{\baselineskip}

If instead we want the functions of ${\mathcal{M}}_{p}^{w}(\Omega)$  to be $p$ summable, we need to introduce some restriction on $w$. Indeed, the following holds.
\begin{prop}\label{prop:molp}
 Let $\Omega$ be an open subset of ${\mathbb{R}}^{n}$. Let $p\in [1,+\infty]$. Let $w$ be a   function from $]0,+\infty[$ to itself.    Let
\begin{equation}
\label{prelprgm1_}
\eta_{w} \equiv\inf_{r\in ]0,+\infty[}w(r)>0\,.
\end{equation}
Then the space ${\mathcal{M}}_{p}^{w}(\Omega)$ is continuously embedded into $L^{p} (\Omega)$, and
\begin{equation}
\label{prelprgm1a}
\|f\|_{  L^{p} (\Omega)  }\leq \eta_{w}^{-1}
\|f\|_{  {\mathcal{M}}_{p}^{w}(\Omega)  }\qquad\forall f\in  {\mathcal{M}}_{p}^{w}(\Omega)\,.
\end{equation}
\end{prop}
{\bf Proof.} Let $f\in {\mathcal{M}}^{w}_{p}(\Omega)$.

If $p<\infty$, then   the inequality
\begin{eqnarray*}
 \lefteqn{
 \int_{ {\mathbb{B}}_{n}(x,r)\cap\Omega}
| f(y)|^{p}\,dy
 \leq w (r)^{-p}\left(\sup_{r\in]0,+\infty[}w (r)\|f\|_{  L^{p}(  {\mathbb{B}}_{n}(x,r)\cap\Omega)    }\right)^{p}
}
\\
\nonumber
&&\qquad\qquad
\leq
\eta_{w} ^{-p}\left(\sup_{r\in]0,+\infty[}w (r)\|f\|_{  L^{p}(  {\mathbb{B}}_{n}(x,r)\cap\Omega)    }\right)^{p}
 \ \ \forall r\in ]0,+\infty[\,,
\end{eqnarray*}
together with the Monotone Convergence Theorem imply that
\begin{equation}
\label{prelprgm4}
\int_{  \Omega}
| f(y)|^{p}\,dy\leq
\eta_{w} ^{-p}\left(\sup_{r\in]0,+\infty[}w (r)\|f\|_{  L^{p}(  {\mathbb{B}}_{n}(x,r)\cap\Omega)    }\right)^{p}\,,
\end{equation}
and thus $f\in L^{p}(\Omega)$.
On the other hand, if $p=+\infty$, then the inequality
\[
{\mathrm{ess\,sup}}_{{\mathbb{B}}_{n}(x,r)\cap\Omega}|f|
\leq \eta_{w} ^{-1} \sup_{r\in]0,+\infty[}w(r)\|f\|_{  L^{\infty}(  {\mathbb{B}}_{n}(x,r)\cap\Omega)    } 
\quad\forall r\in ]0,+\infty[\,,
\]
implies that 
\begin{equation}
\label{prelprgm5}
{\mathrm{ess\,sup}}_{ \Omega}|f|
\leq \eta_{w} ^{-1} \sup_{r\in]0,+\infty[}w (r)\|f\|_{  L^{\infty}(  {\mathbb{B}}_{n}(x,r)\cap\Omega)    } \,,
\end{equation}
and thus $f\in L^{\infty}(\Omega)$.
Inequalities (\ref{prelprgm4}) and (\ref{prelprgm5}) imply that
both in case $p<\infty$ and $p=\infty$, we have 
\[
\|f\|_{L^{p} (\Omega)}\leq \eta_{w} ^{-1}|f|_{+\infty,w,p,\Omega}
=\eta_{w} ^{-1}\|f\|_{{\mathcal{M}}_{p}^{w }(\Omega)}\,.
\]
Hence, the statement follows.\hfill  $\Box$ 

\vspace{\baselineskip}

By applying the previous proposition to the special weight   $w_{\lambda}$, we deduce the validity of the following  corollary. 
\begin{corol}\label{corol:molp}
 Let $\Omega$ be an open subset of ${\mathbb{R}}^{n}$. Let $p\in [1,+\infty]$. 
 Let $\lambda\in [0,+\infty[$. Then the space $M_{p}^{\lambda}(\Omega)$ is continuously embedded into $L^{p} (\Omega)$, and $\|f\|_{L^{p} (\Omega)}\leq  \|f\|_{M_{p}^{\lambda}(\Omega)}$
for all $ f\in M_{p}^{\lambda}(\Omega)$. 
 \end{corol}
If instead we want   ${\mathcal{M}}_{p}^{w}(\Omega)$  to contain the space of $p$ summable functions $L^p(\Omega)$, we need to introduce some other restriction on $w$. 
\begin{prop}\label{prop:mocolp}
 Let $\Omega$ be an open subset of ${\mathbb{R}}^{n}$. Let $p\in [1,+\infty]$. 
  Let $w$ be a function from $]0,+\infty[$ to $[0,+\infty[$. Assume that there exists $r_0\in]0,+\infty[$ such that $w(r_0)\neq 0$. If 
  \begin{equation}\label{prop:mocolp1}
\varsigma_w\equiv\sup_{r\in]0,+\infty[}w(r)<+\infty\,,
\end{equation}
then $L^p(\Omega)$ is continuously embedded into ${\mathcal{M}}_{p}^{w}(\Omega)$ and
\begin{equation}\label{prop:mocolp2}
\|f\|_{ {\mathcal{M}}_{p}^{w}(\Omega)   }\leq \varsigma_w \|f\|_{  L^{p}(  \Omega)    }\qquad\forall f\in  L^{p}(  \Omega) \,.
\end{equation}
\end{prop}
{\bf Proof.} By assumption, we have
\begin{eqnarray}
\label{prop:mocolp3}
\lefteqn{
 w(r)\|f\|_{  L^{p}(  {\mathbb{B}}_{n}(x,\rho)\cap\Omega)    }
 \leq\varsigma_w \|f\|_{  L^{p}(  {\mathbb{B}}_{n}(x,\rho)\cap\Omega)    }
 }
 \\ \nonumber
 &&\qquad\leq\varsigma_w \|f\|_{  L^{p}(  \Omega)    }
\qquad\forall (x,r)\in \Omega\times]0,+\infty[ \,,
\end{eqnarray}
and inequality (\ref{prop:mocolp2}) holds true.\hfill  $\Box$ 

\vspace{\baselineskip}

However condition (\ref{prop:mocolp1}) is quite restrictive and is satisfied only in certain limiting cases as that of the  Morrey spaces ${\mathcal{M}}_{p}^{r^{-0},\rho}(\Omega)$ for all $\rho\in]0,+\infty[$, ${\mathcal{M}}_{p}^{r^{-0}}(\Omega)$, $M^{0}_p(\Omega)$.
\begin{corol}\label{corol:mo0colp}
Let $\Omega$ be an open subset of ${\mathbb{R}}^{n}$. Let $p\in [1,+\infty]$. 
 Then  $L^p(\Omega)$ is continuously embedded into 
 the Morrey spaces ${\mathcal{M}}_{p}^{r^{-0},\rho}(\Omega)$ for all $\rho\in]0,+\infty[$, ${\mathcal{M}}_{p}^{r^{-0}}(\Omega)$, $M^{0}_p(\Omega)$. Moreover, the norm of the corresponding inclusion is less or equal to $1$. 
\end{corol}

By combining Propositions \ref{prop:molp} and \ref{prop:mocolp}, we deduce the validity of the following, that says that if the weight $w$ is both bounded and bounded away from $0$, then the corresponding generalized Morrey space coincides with a Lebesgue space. 
\begin{prop}\label{prop:mo=lp}
 Let $\Omega$ be an open subset of ${\mathbb{R}}^{n}$. Let $p\in [1,+\infty]$. Let $w$ be a   function from $]0,+\infty[$ to itself.    Let
\begin{equation}
\label{prelprgm1}
 \eta_{w} \equiv\inf_{r\in ]0,+\infty[}w(r)>0\,,\qquad \varsigma_w\equiv\sup_{r\in]0,+\infty[}w(r)<+\infty\,.
\end{equation}
Then the space ${\mathcal{M}}_{p}^{w}(\Omega)$ equals $L^{p} (\Omega)$. Moreover,
\[
 \eta_{w}\|f\|_{  L^{p} (\Omega)  }\leq  
 \|f\|_{ {\mathcal{M}}_{p}^{w}(\Omega)   }\leq \varsigma_w \|f\|_{  L^{p}(  \Omega)    }\qquad\forall f\in  L^{p}(  \Omega) \,.
\]
\end{prop}
By applying the previous proposition to the special weights   $r^{-0}=1=w_{0}$, we deduce the validity of the following  corollary. 
\begin{corol}\label{corol:mo=lp}
 Let $\Omega$ be an open subset of ${\mathbb{R}}^{n}$. Let $p\in [1,+\infty]$. 
 
 Then the Morrey spaces  ${\mathcal{M}}_{p}^{r^{-0}}(\Omega)$, $M^{0}_p(\Omega)$ coincide with
 $L^p(\Omega)$  and the corresponding norms are equal.
 \end{corol}
We now look at what happens for the vanishing Morrey spaces under the assumptions of condition (\ref{prelprgm1}). 
\begin{prop}\label{prop:vmo=lp}
 Let $\Omega$ be an open subset of ${\mathbb{R}}^{n}$. Let $p\in [1,+\infty]$. Let $w$ be a   function from $]0,+\infty[$ to itself.    Let condition 
\begin{equation}
\label{prop:vmo=lp_}
 \eta_{w} \equiv\inf_{r\in ]0,+\infty[}w(r)>0\,,\qquad \varsigma_w\equiv\sup_{r\in]0,+\infty[}w(r)<+\infty\,.
\end{equation}
If $p<+\infty$, then 
\[
 {\mathcal{M}}_{p}^{w,0}(\Omega)=\left\{f\in L^{p} (\Omega):\,
 \lim_{\rho\to 0}\left(\sup_{x\in\Omega}\|f\|_{  L^{p}(  {\mathbb{B}}_{n}(x,\rho)\cap\Omega)    }\right)=0\right\}\,. 
 \]
If $p=+\infty$, then ${\mathcal{M}}_{p}^{w,0}(\Omega)=\{0\}$.
\end{prop}
{\bf Proof.} Let $p<+\infty$. By Proposition \ref{prop:mo=lp}, we have ${\mathcal{M}}_{p}^{w}(\Omega)=L^p(\Omega)$ and accordingly
\[
 {\mathcal{M}}_{p}^{w,0}(\Omega)= \left\{f\in L^{p} (\Omega):\,
 \lim_{\rho\to 0}  |f|_{\rho,w,p,\Omega} =0\right\}\,. 
\]
Thus it suffices to show that if  $f\in L^{p} (\Omega)$, then
\begin{equation}
\label{prop:vmo=lp0}
  \eta_{w}\sup_{x\in\Omega}\|f\|_{  L^{p}(  {\mathbb{B}}_{n}(x,\rho)\cap\Omega)    }
 \leq |f|_{\rho,w,p,\Omega}
 \leq\varsigma_w \sup_{x\in\Omega}\|f\|_{  L^{p}(  {\mathbb{B}}_{n}(x,\rho)\cap\Omega)    }
  \qquad\forall \rho\in ]0,+\infty[\,,
\end{equation}
an inequality that we now turn to prove.  Let $f\in L^{p} (\Omega)$, $\rho\in ]0,+\infty[$. 
By the Monotone Convergence Theorem, we know that
\begin{equation}
\label{prop:vmo=lp0a}
\sup_{r\in]0,\rho[}\|f\|_{  L^{p}(  {\mathbb{B}}_{n}(x,r)\cap\Omega)    }
=
\|f\|_{  L^{p}(  {\mathbb{B}}_{n}(x,\rho)\cap\Omega)    }
\qquad\forall x\in \Omega
\,.
\end{equation}
Hence, 
\begin{eqnarray*}
\lefteqn{
 \eta_{w}\sup_{x\in\Omega}\|f\|_{  L^{p}(  {\mathbb{B}}_{n}(x,\rho)\cap\Omega)    }
 =\eta_{w}\sup_{x\in\Omega}\sup_{r\in]0,\rho[}\|f\|_{  L^{p}(  {\mathbb{B}}_{n}(x,r)\cap\Omega)    }
}
\\ \nonumber
&&\qquad
\leq 
\sup_{x\in\Omega}\sup_{r\in]0,\rho[}w(r)\|f\|_{  L^{p}(  {\mathbb{B}}_{n}(x,r)\cap\Omega)    }
\\ \nonumber
&&\qquad
\leq
\sup_{(x,r)\in \Omega\times]0,\rho[} w(r)\|f\|_{  L^{p}(  {\mathbb{B}}_{n}(x,r)\cap\Omega)    }
=
|f|_{\rho,w,p,\Omega}
\\ \nonumber
&&\qquad
\leq 
\varsigma_w\sup_{(x,r)\in \Omega\times]0,\rho[}  \|f\|_{  L^{p}(  {\mathbb{B}}_{n}(x,r)\cap\Omega)    }
\leq 
\varsigma_w\sup_{(x,r)\in \Omega\times]0,\rho[}  \|f\|_{  L^{p}(  {\mathbb{B}}_{n}(x,\rho)\cap\Omega)    }
\\ \nonumber
&&\qquad
=
\varsigma_w\sup_{x\in \Omega\ }  \|f\|_{  L^{p}(  {\mathbb{B}}_{n}(x,\rho)\cap\Omega)    }
\end{eqnarray*}
and thus inequality (\ref{prop:vmo=lp0}) holds true. 

We now  consider case $p=+\infty$. So we now take  $f\in {\mathcal{M}}_{\infty}^{w,0}(\Omega)$ and we show that $f=0$. Since $\eta_w>0$,  Proposition \ref{prop:molp} implies that ${\mathcal{M}}_{\infty}^{w,0}(\Omega)\subseteq L^\infty(\Omega)$ and accordingly $f\in L^{1}_{{\mathrm{loc}}}(\Omega)$.  Then  
 the Lebesgue Differentiation Theorem, implies that
\begin{equation}
\label{prop:vmo=lp1}
\lim_{r\to 0} \Xint-_{ {\mathbb{B}}_{n}(x,r)\cap\Omega}
| f(y)|\,dy
  =| f(x)| 
\qquad{\mathrm{a.a.}}\ x\in \Omega\,.
\end{equation}
So we now try to show that the limit in the right hand side equals $0$ for almost all points $x\in\Omega$. To do so, we observe that
\[
\|f\|_{  L^{\infty}(  {\mathbb{B}}_{n}(x,\rho/2)\cap\Omega)    }
 \leq w(\rho/2)^{-1}
|f|_{\rho,w,\infty,\Omega}
 \leq \eta_w^{-1}
|f|_{\rho,w,\infty,\Omega}\quad\forall (x,\rho)\in\Omega\times]0,+\infty[\,. 
\]
Since $f\in {\mathcal{M}}_{\infty}^{w,0}(\Omega)$, we have
$\lim_{\rho\to 0}|f|_{\rho,w,\infty,\Omega}=0$ and accordingly, 
\[
\lim_{\rho\to 0}
\|f\|_{  L^{\infty}(  {\mathbb{B}}_{n}(x,\rho/2)\cap\Omega)    }=0
\qquad\forall x\in\Omega
\,.
\]
Since
\[\Xint-_{ {\mathbb{B}}_{n}(x,\rho/2)\cap\Omega}
| f(y)|\,dy\leq \|f\|_{  L^{\infty}(  {\mathbb{B}}_{n}(x,\rho/2)\cap\Omega)    }
\]
 for all $(x,\rho)\in\Omega\times ]0,+\infty[$, the above limiting relations imply that
$| f(x)|=0$ for almost all $x\in\Omega$. Hence, the statement holds for $p=+\infty$. \hfill  $\Box$ 

\vspace{\baselineskip}

By applying the previous proposition to the special weights   $r^{-0}=1=w_{0}$, we deduce the validity of the following  corollary. 
\begin{corol}\label{corol:vmo=lp}
 Let $\Omega$ be an open subset of ${\mathbb{R}}^{n}$.   
 If $p\in [1,+\infty[$, then 
\[
 {\mathcal{M}}_{p}^{r^{-0},0}(\Omega)=M_{p}^{0,0}(\Omega)=\left\{f\in L^{p} (\Omega):\,
 \lim_{\rho\to 0}\left(\sup_{x\in\Omega}\|f\|_{  L^{p}(  {\mathbb{B}}_{n}(x,\rho)\cap\Omega)    }\right)=0\right\}\,. 
 \]
If $p=+\infty$, then  ${\mathcal{M}}_{p}^{r^{-0},0}(\Omega)=M_{p}^{0,0}(\Omega)=\{0\}$.
\end{corol}
We now wonder when the functions of ${\mathcal{M}}_{p}^{w}(\Omega)$ are essentially bounded and we prove the following statement. 
\begin{prop}\label{prop:mobd}
 Let $\Omega$ be an open subset of ${\mathbb{R}}^{n}$. Let $p\in [1,+\infty]$. 
  Let $w$ be a function from $]0,+\infty[$ to $[0,+\infty[$. If  
\begin{equation}\label{prop:mobd1}
 \limsup_{r\to 0}w(r)r^{n/p}>0\,,
\end{equation}
  then ${\mathcal{M}}_{p}^{w}(\Omega)$  is continuously embedded into  $L^\infty(\Omega)$.
  \end{prop}
  {\bf Proof.} Let $f\in  {\mathcal{M}}_{p}^{w}(\Omega)$. Since ${\mathcal{M}}_{p}^{w}(\Omega)\subseteq L^p_{{\mathrm{loc}}}(\Omega)\subseteq L^1_{{\mathrm{loc}}}(\Omega)$, the Lebesgue Differentiation Theorem implies that
\[
\lim_{r\to 0} \Xint-_{
{\mathbb{B}}_n(x,r)
}|f(y)|\,dy=|f(x)|\qquad {\mathrm{a.a.}}\ x\in \Omega\,.
\]
So we now try to estimate the limit in the left hand side. To do so, we note that the H\"{o}lder inequality implies that
\begin{eqnarray*}
\lefteqn{
\Xint-_{
{\mathbb{B}}_n(x,r)
}|f(y)|\,dy=m_n({\mathbb{B}}_n(x,r))^{-1}\int_{ {\mathbb{B}}_n(x,r) }|f(y)|\,dy
}
\\ \nonumber
&&\qquad
\leq m_n({\mathbb{B}}_n(x,r))^{-1}m_n({\mathbb{B}}_n(x,r))^{\frac{1}{1}-\frac{1}{p}}
\|f\|_{L^p({\mathbb{B}}_n(x,r) )}
\\ \nonumber
&&\qquad
=m_n({\mathbb{B}}_n(x,r))^{ -\frac{1}{p}}w(r)^{-1}w(r) 
\|f\|_{L^p({\mathbb{B}}_n(x,r)\cap\Omega)}
\\ \nonumber
&&\qquad
\leq\frac{1}{
\omega_n^{ 1/p}r^{n/p}w(r)}\|f\|_{{\mathcal{M}}_{p}^{w}(\Omega)}\,,
\end{eqnarray*}
for all $x\in \Omega$ and $r\in]0,+\infty[$ such that $w(r)\neq 0$ and ${\mathbb{B}}_n(x,r)\subseteq \Omega$.
By assumption (\ref{prop:mobd1}), there exists a sequence $\{r_j\}_{j\in {\mathbb{N}}}$ in $]0,+\infty[$ such that
\[
\lim_{j\to \infty}r_j=0\,,\qquad  \lim_{j\to \infty}w(r_j)r_j^{n/p} =\limsup_{r\to 0}w(r)r^{n/p}>0\,.
\]
Possibly neglecting a finite number of indexes, we can assume that $w(r_j)\neq 0$ for all $j\in{\mathbb{N}}$. 
Hence, we have
\begin{eqnarray*}
\lefteqn{
|f(x)|=\lim_{r\to 0} \Xint-_{
{\mathbb{B}}_n(x,r)
}|f(y)|\,dy 
}
\\ \nonumber
&&\qquad\qquad
= \lim_{j\to \infty}\Xint-_{
{\mathbb{B}}_n(x,r_j)
}|f(y)|\,dy
=
\liminf_{j\to\infty} \Xint-_{
{\mathbb{B}}_n(x,r_j)
}|f(y)|\,dy 
\\ \nonumber
&&\qquad\qquad
\leq
\liminf_{j\to\infty} 
\frac{1}{
\omega_n^{ 1/p}r_j^{n/p}w(r_j)}\|f\|_{{\mathcal{M}}_{p}^{w}(\Omega)}
=
\lim_{j\to\infty} 
\frac{1}{
\omega_n^{ 1/p}r_j^{n/p}w(r_j)}\|f\|_{{\mathcal{M}}_{p}^{w}(\Omega)}
\\ \nonumber
&&\qquad\qquad
=\frac{1}{
\omega_n^{ 1/p}\lim_{j\to\infty}r_j^{n/p}w(r_j)}\|f\|_{{\mathcal{M}}_{p}^{w}(\Omega)}
\\ \nonumber
&&\qquad\qquad
=\frac{1}{
\omega_n^{ 1/p}\limsup_{r\to 0}w(r)r^{n/p}}\|f\|_{{\mathcal{M}}_{p}^{w}(\Omega)}
\qquad {\mathrm{a.a.}}\ x\in \Omega\,,
\end{eqnarray*}
where we understand that the right hand side equals $0$ if $\limsup_{r\to 0}w(r)r^{n/p}$ equals $+\infty$. 
Hence, the proof is complete. \hfill  $\Box$ 

\vspace{\baselineskip}
  
 Then the following immediate corollary holds. 
 \begin{corol}\label{corol:mobd}
 Let $\Omega$ be an open subset of ${\mathbb{R}}^{n}$. Let $p\in [1,+\infty]$. 
  Let $\lambda=n/p$. Then  the  Morrey spaces ${\mathcal{M}}_{p}^{r^{-n/p},\rho}(\Omega)$ for all $\rho\in]0,+\infty[$, ${\mathcal{M}}_{p}^{r^{-n/p}}(\Omega)$, $M^{n/p}_p(\Omega)$  are all  continuously embedded into $L^\infty(\Omega)$.
  \end{corol}
Conversely, we may ask under what assumptions  the space $L^\infty(\Omega)$ is included into
${\mathcal{M}}_{p}^{w}(\Omega)$ and we prove the following statement.
\begin{prop}\label{prop:bdmo}
 Let $\Omega$ be an open subset of ${\mathbb{R}}^{n}$. Let $p\in [1,+\infty]$. 
  Let $w$ be a function from $]0,+\infty[$ to $[0,+\infty[$. Assume that there exists $r_0\in]0,+\infty[$ such that $w(r_0)\neq 0$. If 
\begin{equation}\label{prop:mobd1a}
 \sup_{r\in]0,+\infty[ }w(r)r^{n/p}<+\infty\,,
\end{equation}
  then $L^\infty(\Omega)$   is continuously embedded into  ${\mathcal{M}}_{p}^{w}(\Omega)$.
\end{prop}
{\bf Proof.} If $f\in L^\infty(\Omega)$, then the H\"{o}lder inequality implies that
\begin{eqnarray}\label{prop:mobd2}
\lefteqn{
w (r)\|f\|_{ L^{p}(  {\mathbb{B}}_{n}(x,r)\cap\Omega)}
}
\\
\nonumber
&&\qquad
\leq w (r)m_{n}( {\mathbb{B}}_{n}(x,r)\cap\Omega)^{1/p}\|f\|_{L^{\infty}(\Omega)}
\\
\nonumber
&&\qquad
\leq w (r)\omega_n^{1/p}r^{n/p}\|f\|_{L^{\infty}(\Omega)}
\qquad\forall r\in ]0,+\infty[\,.
\end{eqnarray}
and thus the statement holds true.\hfill  $\Box$ 

\vspace{\baselineskip}

Then the Corollary \ref{corol:mobd} and the previous proposition imply the validity of the 
following corollary. 
 \begin{corol}\label{corol:bdmo}
 Let $\Omega$ be an open subset of ${\mathbb{R}}^{n}$. Let $p\in [1,+\infty]$. 
  Let $\lambda=n/p$. Then 
  \begin{equation}\label{corol:bdmo1}
{\mathcal{M}}_{p}^{r^{-n/p}}(\Omega)=L^\infty(\Omega)\,,\qquad
{\mathcal{M}}_{p}^{r^{-n/p},\rho}(\Omega)=L^\infty(\Omega)\qquad\forall \rho\in]0,+\infty[\,,
\end{equation}
and the corresponding norms are equivalent. 
  \end{corol}
  Instead for the weight $w_{n/p}$ we prove the following statement. 
  \begin{corol}\label{corol:bdMo}
 Let $\Omega$ be an open subset of ${\mathbb{R}}^{n}$ of finite measure. Let $p\in [1,+\infty]$. 
  Let $\lambda=n/p$. Then 
  \begin{equation}\label{corol:bdMo1}
M_{p}^{n/p}(\Omega) =L^\infty(\Omega)\,,
\end{equation}
and the corresponding norms are equivalent. 
\end{corol}
{\bf Proof.} If $f\in L^\infty(\Omega)$, then the H\"{o}lder inequality implies that
\begin{eqnarray*}
\lefteqn{
w_{n/p}(r)\|f\|_{ L^{p}(  {\mathbb{B}}_{n}(x,r)\cap\Omega)}
}
\\
\nonumber
&&\qquad
\leq w_{n/p}(r)m_{n}( {\mathbb{B}}_{n}(x,r)\cap\Omega)^{1/p}\|f\|_{L^{\infty}(\Omega)}\qquad\forall r\in ]0,+\infty[\,.
\end{eqnarray*}
Since
\begin{eqnarray*}
& w_{n/p}(r)m_{n}( {\mathbb{B}}_{n}(x,r)\cap\Omega)^{1/p}
\leq \omega_{n}^{1/p} \qquad\qquad&{\mathrm{if}}\ r\in]0,1]\,,
\\
& w_{n/p}(r)m_{n}( {\mathbb{B}}_{n}(x,r)\cap\Omega)^{1/p}
\leq m_{n}(\Omega)^{1/p} \qquad&{\mathrm{if}}\
 r\in]1,+\infty[\,,
\end{eqnarray*}
the above inequality implies that $L^{\infty}(\Omega)
 \subseteq
 M_{p}^{\lambda}(\Omega)$ and that
\begin{eqnarray*}
\lefteqn{
\|f\|_{ M_{p}^{\lambda}(\Omega) } 
}
\\ \nonumber
&&\quad
\equiv
\sup_{(x,r)\in]0,+\infty[}
 w_{n/p}(r)\|f\|_{ L^{p}(  {\mathbb{B}}_{n}(x,r)\cap\Omega)}
\leq \max\{\omega_{n}^{1/p},
 m_{n}(\Omega)^{1/p}\}\|f\|_{L^{\infty}(\Omega)}\,.
\end{eqnarray*}\hfill  $\Box$ 

\vspace{\baselineskip}
  
  \begin{rem}{\em Let $\Omega$ be an open subset of ${\mathbb{R}}^{n}$. Let $p\in [1,+\infty[$. The equality $M^{n/p}(\Omega)= L^\infty(\Omega)$ can hold only if 
 $L^\infty(\Omega)\subseteq L^p(\Omega)$. Indeed, $M^{n/p}_p(\Omega)\subseteq L^p(\Omega)$. 
 Now inclusion $L^\infty(\Omega)\subseteq L^p(\Omega)$ implies that $\chi_\Omega\in L^p(\Omega)$, and accordingly that $m_n(\Omega)<+\infty$. Thus   we can have $M^{n/p}_p(\Omega)= L^\infty(\Omega)$ only if $m_n(\Omega)<+\infty$. 
 
 If instead $p=+\infty$, we have already proved that $M^{n/p}_p(\Omega)=M^{0}_p(\Omega)= L^\infty(\Omega)$, no matter whether $m_n(\Omega)$ is finite or not. 
 }\end{rem}
 
 Then we have the following statement. 
\begin{prop}
\label{bgm}
Let $\Omega$ be an open subset of ${\mathbb{R}}^{n}$. 
 Let $p\in[1,+\infty]$.  Let $w$ be a   function from $]0,+\infty[$ to $[0,+\infty[$.
 Assume that there exists $r_0\in]0,+\infty[$ such that $w(r_0)\neq 0$. If
  \begin{equation}\label{bgm1}
\lim_{\rho\to 0}\sup_{r\in]0,\rho]} w(r)r^{n/p} =0\,,
\end{equation}
then
\[
{\mathcal{M}}^{w}_{p}(\Omega)\cap L^{\infty}(\Omega)\subseteq 
{\mathcal{M}}^{w,0}_{p}(\Omega)\,.
\]
\end{prop}
{\bf Proof.} Obviously, ${\mathcal{M}}^{w}_{p}(\Omega)\cap L^{\infty}(\Omega)\subseteq 
{\mathcal{M}}^{w}_{p}(\Omega)$. By the H\"{o}lder inequality, we have
\begin{eqnarray*}
\lefteqn{\sup_{(x,r)\in\Omega\times ]0,\rho[}
w (r)\|f\|_{  L^{p}(  {\mathbb{B}}_{n}(x,r)\cap\Omega)    }
}
\\
\nonumber
&&\qquad
\leq
 \sup_{(x,r)\in\Omega\times ]0,\rho[}
w (r)\omega_{n}^{1/p}r^{n/p}
\|f\|_{  L^{\infty}(  {\mathbb{B}}_{n}(x,r)\cap\Omega)    }
\\
\nonumber
&&\qquad
\leq
\|f\|_{  L^{\infty}(   \Omega)    }
\omega_{n}^{1/p}
\sup_{(x,r)\in\Omega\times ]0,\rho[}
w (r)r^{n/p}
\qquad\forall \rho\in]0,+\infty[\,.
\end{eqnarray*}
Hence, condition (\ref{bgm1}) implies the validity of the statement.
\hfill  $\Box$

\vspace{\baselineskip}

Then we have the following immediate corollary for the weights $w_{\lambda,\rho}$, with $\rho\in]0,+\infty[$, $r^{-\lambda}$, $w_\lambda$.
\begin{corol}
\label{bm}
Let $\Omega$ be an open subset of ${\mathbb{R}}^{n}$. 
 Let $p\in[1,+\infty[$.  If $\lambda\in[0,n/p[$, then 
 \begin{enumerate}
\item[(i)] ${\mathcal{M}}^{r^{-\lambda},\rho,0}_{p}(\Omega)\cap L^{\infty}(\Omega)\subseteq {\mathcal{M}}^{r^{-\lambda},\rho,0}_{p}(\Omega)$ for all $\rho\in]0,+\infty[$. 
\item[(ii)] ${\mathcal{M}}^{r^{-\lambda}, 0}_{p}(\Omega)\cap L^{\infty}(\Omega)\subseteq {\mathcal{M}}^{r^{-\lambda}, 0}_{p}(\Omega)$. 
\item[(iii)] $M^{\lambda}_{p}(\Omega)\cap L^{\infty}(\Omega)\subseteq M^{\lambda,0}_{p}(\Omega)$. 
\end{enumerate}
\end{corol}

 \begin{rem}{\em Let $\Omega$ be an open subset of ${\mathbb{R}}^{n}$ of finite measure. Let $p\in [1,+\infty[$. Let $\lambda\in ]0,n/p[$. 
Let $p_\lambda\in]p,+\infty[$ be such that
\[
\lambda=\frac{n}{p}-\frac{n}{p_\lambda}\,.
\]
Then inequality (\ref{introd3}) implies that $L^{p_\lambda}(\Omega)$ is continuously embedded into $M^{\lambda}_p(\Omega)$.}
 \end{rem}

\section{Restrictions and extensions of functions in generalized Morrey spaces}

By the definition of norm in a  generalized Morrey space, it follows that if the absolute value of  a measurable function can be estimated in terms of a function in a generalized Morrey space, than the function belongs to the generalized Morrey space too. More precisely, the following holds.
\begin{lem}
\label{lem:molat}
Let $\Omega$ be an open subset of ${\mathbb{R}}^{n}$. Let $p\in [1,+\infty]$. Let $w$ be a   function from $]0,+\infty[$ to $[0,+\infty[$. Assume that there exists $r_0\in]0,+\infty[$ such that $w(r_0)\neq 0$. 

Let $f,g$ be measurable functions from $\Omega$ to ${\mathbb{R}}$ such that $|f|\leq |g|$ almost everywhere in $\Omega$. Then we have
\[
|f|_{\rho,w,p,\Omega}\leq |g|_{\rho,w,p,\Omega}\qquad
\forall \rho\in]0,+\infty]\,.
\]
In particular, if $g\in {\mathcal{M}}_{p}^{w}(\Omega)$ and if
$g\in {\mathcal{M}}_{p}^{w,0}(\Omega)$, then 
$f\in {\mathcal{M}}_{p}^{w}(\Omega)$ and  
$f\in {\mathcal{M}}_{p}^{w,0}(\Omega)$, respectively.
\end{lem}
Next we collect in the following statement some properties of restrictions and extensions in generalized Morrey spaces. 
\begin{prop}
\label{prelprgm}
Let $\Omega$ be an open subset of ${\mathbb{R}}^{n}$. Let $p\in [1,+\infty]$. Let $w$ be a   function from $]0,+\infty[$ to $[0,+\infty[$.  Assume that there exists $r_0\in]0,+\infty[$ such that $w(r_0)\neq 0$. Then the following statements hold.
\begin{enumerate}
  \item[(i)] Let $V$ be an open subset of $\Omega$. Then
\[
|f|_{\rho,w,p,V}\leq |f|_{\rho,w,p,\Omega}\qquad
\forall \rho\in]0,+\infty]\,,
\]
for all measurable functions $f$ from $\Omega$ to ${\mathbb{R}}$. In particular, the restriction map is linear and continuous from
$ {\mathcal{M}}_{p}^{w}(\Omega)$ to $ {\mathcal{M}}_{p}^{w}(V)$ and maps
$ {\mathcal{M}}_{p}^{w,0}(\Omega)$ to $ {\mathcal{M}}_{p}^{w,0}(V)$.
\item[(ii)] Let $w(r)>0$ for all $r\in]0,+\infty[$. Let $w$ satisfy the following `doubling' condition  
\begin{equation}
\label{prelprgm2}
\sigma_{w} \equiv\sup_{r\in ]0,+\infty[}w(r)w^{-1}(2r)<+\infty \,.
\end{equation}
Then  
\[
|E_{\Omega}f|_{\rho,w,p,{\mathbb{R}}^{n}}
\leq \sigma_{w}|f|_{2\rho,w,p,\Omega}\qquad
\forall \rho\in]0,+\infty]\,,
\]
for all measurable functions from $\Omega$ to ${\mathbb{R}}$ such that  $f\in L^{p}({\mathbb{B}}_{n}(x,r)\cap\Omega)$ for all $(x,r)\in \Omega\times ]0,+\infty[$. (Here we understand that 
$2\rho=+\infty$ if $\rho=+\infty$ and that the right hand side equals $+\infty$ if $|f|_{2\rho,w,p,\Omega}=+\infty$.) In particular, the extension  map $E_{\Omega}$ is linear and continuous from
$ {\mathcal{M}}_{p}^{w }(\Omega)$ to $ {\mathcal{M}}_{p}^{w }({\mathbb{R}}^{n})$ and maps
$ {\mathcal{M}}_{p}^{w ,0}(\Omega)$ to $ {\mathcal{M}}_{p}^{w ,0}({\mathbb{R}}^{n})$.
\end{enumerate}
\end{prop}
{\bf Proof.}  The inequality in statement   (i) is obvious and implies the validity of the corresponding statement.\par

We now prove statement (ii).   Let $f$ be a measurable function  from $\Omega$ to ${\mathbb{R}}$ such that  $f\in L^{p}({\mathbb{B}}_{n}(x,r)\cap\Omega)$ for all $r\in]0,+\infty[$.  Let $\rho\in]0,+\infty]$. It suffices to prove that if 
 $(x,r)\in {\mathbb{R}}^{n}\times ]0,\rho[$,
 then 
 \[
 w(r)\|E_\Omega f\|_{L^{p}( {\mathbb{B}}_{n}(x,r)\cap {\mathbb{R}}^n) )}=w(r)\|f\|_{L^{p}( {\mathbb{B}}_{n}(x,r)\cap\Omega) )}
 \leq \sigma_{w}|f|_{2\rho,w,p,\Omega}\,.
 \]
 Let ${\mathrm{dist}}\,(x,\Omega)$ denote the distance between $x$ and $\Omega$, where we understand that such a distance equals $+\infty$ if $\Omega$ is empty.  We distinguish two cases.
 
  If $r\leq {\mathrm{dist}}\,(x,\Omega)$, then   ${\mathbb{B}}_{n}(x,r)\cap\Omega$ is empty and 
  \[
  w(r)\|E_\Omega f\|_{L^{p}( {\mathbb{B}}_{n}(x,r)\cap {\mathbb{R}}^n) )}=
  w(r)\|f\|_{L^{p}( {\mathbb{B}}_{n}(x,r)\cap\Omega) )}=0\leq \sigma_{w}|f|_{2\rho,w,p,\Omega}\,.
  \]
  If $r> {\mathrm{dist}}\,(x,\Omega)$, then there exists $\xi_{x}\in 
    {\mathbb{B}}_{n}(x,r)\cap\Omega$. By the triangular inequality, we have   
$ {\mathbb{B}}_{n}(x,r)\subseteq  {\mathbb{B}}_{n}(\xi_{x},2r)    $
 and accordingly
\begin{eqnarray*}
\lefteqn{
w(r)\|E_\Omega f\|_{L^{p}( {\mathbb{B}}_{n}(x,r)\cap {\mathbb{R}}^n) )}=
w(r)\|f\|_{L^{p}( {\mathbb{B}}_{n}(x,r)\cap\Omega) )}
}
\\ \nonumber
&&\qquad\qquad\qquad\qquad
\leq
 w(r)w(2r)^{-1} w(2r)\|f\|_{L^{p}( {\mathbb{B}}_{n}(\xi_{x},2r)\cap\Omega) )}
  \leq \sigma_{w}|f|_{2\rho,w,p,\Omega}\,.
\end{eqnarray*}
 Hence, (ii) holds true. \hfill  $\Box$

\vspace{\baselineskip}

\begin{corol}\label{prelprgmrr}
 Let $\Omega$ be an open subset of ${\mathbb{R}}^{n}$. Let $p\in [1,+\infty]$. 
 Let $\lambda\in [0,+\infty[$. Then the following statements hold.
 \begin{enumerate}
\item[(i)] 
The extension  operator $E_{\Omega}$ is linear and continuous from
$ {\mathcal{M}}_{p}^{r^{-\lambda} }(\Omega)$ to $ {\mathcal{M}}_{p}^{r^{-\lambda} }({\mathbb{R}}^{n})$ and maps
$ {\mathcal{M}}_{p}^{r^{-\lambda} ,0}(\Omega)$ to $ {\mathcal{M}}_{p}^{r^{-\lambda} ,0}({\mathbb{R}}^{n})$.
\item[(ii)] The extension  operator $E_{\Omega}$ is linear and continuous from
$ M_p^\lambda(\Omega)$ to $M_p^\lambda({\mathbb{R}}^{n})$ and maps
$ M_p^{\lambda,0}(\Omega)$ to $M_p^{\lambda,0}({\mathbb{R}}^{n})$.
\end{enumerate}
\end{corol}
{\bf Proof.}  Since   $r^{-\lambda} \left((2r)^{-\lambda}\right)^{-1}=2^{\lambda}$ for all $r\in]0,+\infty[$ and 
\[
w_{\lambda}(r)w^{-1}_{\lambda}(2r)=
\left\{
\begin{array}{ll}
r^{-\lambda} (2r)^{\lambda}=2^{\lambda}
&
{\mathrm{if}}\ 2r\in]0,1]\,,
\\
r^{-\lambda} & {\mathrm{if}}\ 2r\in]1,2]\,,
\\
1 & {\mathrm{if}}\ 2r\in]2,+\infty[\,,
\end{array}
\right.
\]
we have
\begin{equation}
\label{prelprgmr}
\sigma_{r^{-\lambda}}=2^{\lambda}\,,
\qquad
\sigma_{w_{\lambda}}=2^{\lambda}\,,
\end{equation}
and accordingly,   the statement of Proposition \ref{prelprgm} (ii) applies to the weights
$r^{-\lambda}$ and
 $w=w_{\lambda}$.  \hfill  $\Box$ 

\vspace{\baselineskip}

Proposition \ref{prelprgm} (ii) requires that $w$ does not vanish. In case $w$ vanishes we only prove the following for the specific weight $w_{\lambda, \rho}$ with $\rho\in ]0,+\infty[$. 

\begin{prop}\label{prop:mrhoext}
Let $\Omega$ be an open subset of ${\mathbb{R}}^{n}$. Let $p\in [1,+\infty]$. 
 Let $\lambda\in [0,+\infty[$. Then the following statements hold.
 \begin{enumerate}
\item[(i)] If $f$ is a measurable function from $\Omega$ to ${\mathbb{R}}$ 
such that  $f\in L^{p}({\mathbb{B}}_{n}(x,r)\cap\Omega)$ for all $(x,r)\in \Omega\times ]0,+\infty[$, then 
 \begin{equation}\label{prop:mrhoext1}
|E_\Omega f|_{\rho,r^{-\lambda},p,{\mathbb{R}}^n}
\leq 2^\lambda |  f|_{2\rho,r^{-\lambda},p,\Omega}\qquad\forall \rho\in]0,+\infty[\,.
\end{equation}
\item[(ii)] If $\rho\in ]0,+\infty[$, then the 
extension  map $E_{\Omega}$ is linear and continuous from
$ {\mathcal{M}}_{p}^{r^{-\lambda},\rho }(\Omega)$ to $ {\mathcal{M}}_{p}^{r^{-\lambda},\rho }({\mathbb{R}}^{n})$ and maps
$ {\mathcal{M}}_{p}^{r^{-\lambda} ,\rho,0}(\Omega)$ to $ {\mathcal{M}}_{p}^{r^{-\lambda},\rho ,0}({\mathbb{R}}^{n})$.
\end{enumerate}
 \end{prop}
{\bf Proof.}    Since  $r^{-\lambda} \left((2r)^{-\lambda}\right)^{-1}=2^{\lambda}$ for all $r\in]0,+\infty[$, Proposition \ref{prelprgm} (ii) with $w(r)=r^{-\lambda}$ implies the validity of the inequality (\ref{prop:mrhoext1}) in (i). 
  By inequality (\ref{prop:mrhoext1}), the map  $E_\Omega$ is linear and continuous from 
$ {\mathcal{M}}_{p}^{r^{-\lambda},2\rho }(\Omega)$ to $ {\mathcal{M}}_{p}^{r^{-\lambda},\rho }({\mathbb{R}}^{n})$.
By Proposition \ref{prop:moud}, we know that
 ${\mathcal{M}}_{p}^{r^{-\lambda},\rho }(\Omega)={\mathcal{M}}_{p}^{r^{-\lambda},2\rho }(\Omega)$ and that the corresponding norms are equivalent. Hence, $E_\Omega$ is linear and continuous from 
$ {\mathcal{M}}_{p}^{r^{-\lambda}, \rho }(\Omega)$ to $ {\mathcal{M}}_{p}^{r^{-\lambda},\rho }({\mathbb{R}}^{n})$.
 By inequality (\ref{prop:mrhoext1}),   $E_\Omega f$ belongs to $   {\mathcal{M}}_{p}^{r^{-\lambda},\rho,0 }({\mathbb{R}}^{n})$
 whenever $f\in {\mathcal{M}}_{p}^{r^{-\lambda}, \rho,0 }(\Omega)$.\hfill  $\Box$ 

\vspace{\baselineskip}

\section{Completeness of  generalized Morrey spaces} 

\begin{thm}
\label{gmco}
Let $\Omega$ be an open subset of ${\mathbb{R}}^{n}$. Let $p\in[1,+\infty]$. Let $w$ be a   function from $]0,+\infty[$ to $[0,+\infty[$. Assume that there exists $r_0\in]0,+\infty[$ such that $w(r_0)\neq 0$.

Then ${\mathcal{M}}^{w }_{p}(\Omega)$ is a Banach space. 
\end{thm}
{\bf Proof.} Let $\{u_{j}\}_{ j\in  {\mathbb{N}} }$ be a Cauchy sequence in ${\mathcal{M}}^{w }_{p}(\Omega)$. Since  
${\mathcal{M}}^{w }_{p}(\Omega)$ is continuously embedded into $L^{p}_{{\mathrm{loc}} }(\Omega)$ and $L^{p}_{{\mathrm{loc}} }(\Omega)$ is complete, there exist a 
$u\in L^{p}_{{\mathrm{loc}} }(\Omega)$ such that
\[
\lim_{j\to\infty}u_{j }=u\qquad{\mathrm{in}}\ L^{p}_{{\mathrm{loc}} }(\Omega)\,.
\]
Furthermore, there exists a subsequence $\{ u_{j_{k}}\}_{k\in {\mathbb{N}}}$ which converges  pointwise almost everywhere to $u$   in $\Omega$. 

We now show that $u\in {\mathcal{M}}^{w }_{p}(\Omega)$. 
Let $x\in \Omega$, $r\in]0,+\infty[$, $w(r)\neq 0$.
By Proposition \ref{lpim} (i), the restriction map is linear and continuous from 
${\mathcal{M}}^{w }_{p}(\Omega)$ to $L^{p}({\mathbb{B}}_{n}(x,r)\cap\Omega)$. Since $\{ u_{j_{k}}\}_{k\in {\mathbb{N}}}$ is a Cauchy sequence in ${\mathcal{M}}^{w }_{p}(\Omega)$, then the sequence of restrictions $\{
(u_{j})_{| {\mathbb{B}}_{n}(x,r)\cap\Omega)}
\}_{j  \in   {\mathbb{N}}   }$ is a Cauchy sequence in the complete space $L^{p}({\mathbb{B}}_{n}(x,r)\cap\Omega)$. Then there exists $v_{x,r}\in  L^{p}({\mathbb{B}}_{n}(x,r)\cap\Omega)$ such that
\[
\lim_{j\to\infty}(u_{j})_{| {\mathbb{B}}_{n}(x,r)\cap\Omega }
=
v_{x,r}\qquad{\mathrm{in}}\ L^{p}({\mathbb{B}}_{n}(x,r)\cap\Omega)\,.
\]
Accordingly, we have 
\[
\lim_{k\to\infty}(u_{j_{k}})_{| {\mathbb{B}}_{n}(x,r)\cap\Omega }
=
v_{x,r}\qquad{\mathrm{in}}\ L^{p}({\mathbb{B}}_{n}(x,r)\cap\Omega)\,,
\]
Now we know that the sequence in the left hand side converges pointwise almost everywhere in $\Omega$ and thus in  ${\mathbb{B}}_{n}(x,r)\cap\Omega$. Then we must necessarily have 
\[
u=
v_{x,r}\qquad \text{a.e.\ in}\ {\mathbb{B}}_{n}(x,r)\cap\Omega\,.
\]
In particular, 
\[
\lim_{j\to\infty}(u_{j})_{| {\mathbb{B}}_{n}(x,r)\cap\Omega }
=
u_{| {\mathbb{B}}_{n}(x,r)\cap\Omega }\qquad{\mathrm{in}}\ L^{p}({\mathbb{B}}_{n}(x,r)\cap\Omega)\,.
\]
Next we note that
\begin{eqnarray}
\label{gmco1}
\lefteqn{w(r)\|u\|_{L^{p}({\mathbb{B}}_{n}(x,r)\cap\Omega)}
}
\\
\nonumber
&& 
\leq
w(r)\|u-u_{j}\|_{L^{p}({\mathbb{B}}_{n}(x,r)\cap\Omega)}
+
w(r)\|u_{j}\|_{L^{p}({\mathbb{B}}_{n}(x,r)\cap\Omega)}
\\
\nonumber
&&
\leq
w(r)\|u-u_{j}\|_{L^{p}({\mathbb{B}}_{n}(x,r)\cap\Omega)}
+ \|u_{j}\|_{ {\mathcal{M}}_{p}^{w}(\Omega)}
\\
\nonumber
&&\leq
w(r)\|u-u_{j}\|_{L^{p}({\mathbb{B}}_{n}(x,r)\cap\Omega)}
+
\sup_{l\in {\mathbb{N}} }\|u_{l}\|_{ {\mathcal{M}}_{p}^{w}(\Omega)}
\qquad\forall j\in {\mathbb{N}}\,.
\end{eqnarray}
Since $\{u_{j}\}_{j\in {\mathbb{N}}}$ is a Cauchy sequence in ${\mathcal{M}}^{w }_{p}(\Omega)$, the supremum in the right hand side is finite. Then we take the limit as $j$ tends to infinity in the left and right hand sides of inequalities (\ref{gmco1}) and deduce that
\[
w(r)\|u\|_{L^{p}({\mathbb{B}}_{n}(x,r)\cap\Omega)}\leq 0
+\sup_{l\in {\mathbb{N}} }\|u_{l}\|_{ {\mathcal{M}}_{p}^{w}(\Omega)}\,.
\]
 Since $x$ and $r$ have been chosen arbitrarily with the restriction that $w(r)\neq 0$, we immediately deduce that $u\in {\mathcal{M}}_{p}^{w}(\Omega)$ and that
\begin{eqnarray*}
\lefteqn{
\|u\|_{{\mathcal{M}}_{p}^{w}(\Omega)}
 =\sup_{(x,r)\in\Omega\times]0,+\infty[}w(r)\|f\|_{  L^{p}(  {\mathbb{B}}_{n}(x,r)\cap\Omega)    }
}
\\ \nonumber
&&\qquad\qquad
=
\sup_{(x,r)\in\Omega\times]0,+\infty[, w(r)\neq 0}w(r)\|f\|_{  L^{p}(  {\mathbb{B}}_{n}(x,r)\cap\Omega)    }
\leq
\sup_{l\in {\mathbb{N}} }\|u_{l}\|_{ {\mathcal{M}}_{p}^{w}(\Omega)}\,.
\end{eqnarray*}
Next we show that $\{u_{j}\}_{j\in {\mathbb{N}}}$ converges to 
$u$ in ${\mathcal{M}}^{w }_{p}(\Omega)$. Let $\epsilon\in]0,+\infty[$. 
 Then there exists $j_{0}\in {\mathbb{N}}$ such that
 \[
 \|u_{j}-u_{l}\|_{{\mathcal{M}}^{w }_{p}(\Omega)}\leq \epsilon\qquad
 \forall l,j\geq j_{0}\,.
 \]
Now let $(x,r)\in \Omega\times]0,+\infty[$, $w(r)\neq 0$. Then we have 
\[
w(r)\|u_{j}-u_{l}\|_{L^{p}({\mathbb{B}}_{n}(x,r)\cap\Omega)}
\leq
\|u_{j}-u_{l}\|_{{\mathcal{M}}^{w }_{p}(\Omega)}\leq\epsilon\qquad
 \forall l,j\geq j_{0}\,.
\]
Then we can take the limit in the left and right hand side as $l$ tends to infinity, and obtain
\[
w(r)\|u_{j}-u \|_{L^{p}({\mathbb{B}}_{n}(x,r)\cap\Omega)}
\leq
\epsilon\qquad
 \forall j\geq j_{0}\,.
\]
Since $(x,r)$ is arbitrary with the restriction that $w(r)\neq 0$, we have
\begin{eqnarray*}
\lefteqn{
\|u_{j}-u \|_{{\mathcal{M}}^{w }_{p}(\Omega)}=\sup_{(x,r)\in\Omega\times]0,+\infty[}w(r)\|u_{j}-u \|_{L^{p}({\mathbb{B}}_{n}(x,r)\cap\Omega)}
}
\\ \nonumber
&&\qquad\qquad
=\sup_{(x,r)\in\Omega\times]0,+\infty[, w(r)\neq 0}w(r)\|u_{j}-u \|_{L^{p}({\mathbb{B}}_{n}(x,r)\cap\Omega)}
\leq
\epsilon\qquad
 \forall j\geq j_{0}\,.
\end{eqnarray*}
Hence, $\{u_{j}\}_{j\in {\mathbb{N}}}$ converges to 
$u$ in ${\mathcal{M}}^{w }_{p}(\Omega)$.\hfill  $\Box$

\vspace{\baselineskip}

Next we show that the  vanishing  generalized Morrey spaces are closed in the Morrey spaces. 
\begin{thm}
\label{lmc}
Let $\Omega$ be an open subset of ${\mathbb{R}}^{n}$. Let $p\in[1,+\infty]$. Let $w$ be a function from $]0,+\infty[$ to $[0,+\infty[$. 
Assume that there exists $r_0\in]0,+\infty[$ such that $w(r_0)\neq 0$.

Then ${\mathcal{M}}^{w ,0}_{p}(\Omega)$ is a closed subspace of ${\mathcal{M}}^{w }_{p}(\Omega)$. 
\end{thm}
{\bf Proof.} Let $f\in   {\mathcal{M}}^{w }_{p}(\Omega)$. Let $\{f_{j}\}_{
j\in {\mathbb{N}} }$ be a sequence in ${\mathcal{M}}^{w ,0}_{p}(\Omega)$
such that $f=\lim_{j\to\infty}f_{j}$ in ${\mathcal{M}}^{w }_{p}(\Omega)$. We must show that
$f\in   {\mathcal{M}}^{w,0 }_{p}(\Omega)$, \textit{i.e.}, that
\[
\lim_{\rho\to 0}| f |_{\rho,w,p,\Omega}=0\,.
\]
To do so, we note that
\begin{eqnarray*}
\lefteqn{
w(r)\|f\|_{L^{p}({\mathbb{B}}_{n}(x,r)\cap\Omega)}
}
\\
\nonumber
&&\qquad
\leq 
w(r)\|f_{j}\|_{L^{p}({\mathbb{B}}_{n}(x,r)\cap\Omega)}
+
w(r)\|f-f_{j}\|_{L^{p}({\mathbb{B}}_{n}(x,r)\cap\Omega)}
\\
\nonumber
&&\qquad
\leq | f_{j}|_{\rho,w,p,\Omega}+
|f- f_{j}|_{+\infty,w,p,\Omega}\,,
\end{eqnarray*}
for all $(x,r)\in \Omega\times]0,\rho[$, $\rho\in ]0,+\infty[$  and for all $j\in {\mathbb{N}}$.
Then we have 
\[
| f |_{\rho,w,p,\Omega}\leq | f_{j}|_{\rho,w,p,\Omega}+
\|f- f_{j}\|_{ {\mathcal{M}}^{w }_{p}(\Omega) }\qquad\forall \rho\in ]0,+\infty[\,,
\]
for all $j\in {\mathbb{N}}$. Now let $\eta>0$. By assumption, there exists $j_{0}\in {\mathbb{N}}$ such that
\[
 \|f- f_{j}\|_{ {\mathcal{M}}^{w }_{p}(\Omega) }
\leq \eta/2 \qquad\forall j\geq j_{0}\,.
\]
Since $f_{j_{0}}\in {\mathcal{M}}^{w ,0}_{p}(\Omega)$, then there exists $\rho_{\eta}>0$ such that 
\[
 | f_{j_{0}}|_{\rho,w,p,\Omega}\leq\eta/2\qquad\forall \rho\in]0,\rho_{\eta}]\,.
\]
Then we have
\[
| f |_{\rho,w,p,\Omega}\leq | f_{j_{0}}|_{\rho,w,p,\Omega}+
\|f- f_{j_{0}}\|_{ {\mathcal{M}}^{w }_{p}(\Omega) }
\leq (\eta/2)+ (\eta/2)=\eta\qquad\forall \rho\in]0,\rho_{\eta}]\,,
\]
and accordingly $f\in {\mathcal{M}}^{w ,0}_{p}(\Omega)$.
\hfill  $\Box$

\vspace{\baselineskip}

Then we have the following statement.
\begin{prop}
\label{intml}
Let $\Omega$ be an open subset of ${\mathbb{R}}^{n}$. 
Let $p\in [1,+\infty]$. Let $w$ be a   function from $]0,+\infty[$ to $[0,+\infty[$. Assume that there exists $r_0\in]0,+\infty[$ such that $w(r_0)\neq 0$. Then the following two statements hold.
\begin{enumerate}
\item[(i)] ${\mathcal{M}}_{p}^{w,0}(\Omega)$ is a Banach subspace of ${\mathcal{M}}_{p}^{w}(\Omega)$.
\item[(ii)] The function $\|\cdot\|_{
{\mathcal{M}}_{p}^{w}(\Omega)
\cap
L^{p}(\Omega)}$ from
 ${\mathcal{M}}_{p}^{w}(\Omega)
\cap
L^{p}(\Omega)$
  to $[0,+\infty[$ defined by
\[
\|f\|_{
{\mathcal{M}}_{p}^{w}(\Omega)
\cap
L^{p}(\Omega)}
\equiv
\max\{
\|f\|_{
{\mathcal{M}}_{p}^{w}(\Omega)}
,
\|f\|_{L^{p}(\Omega)}
\}
\qquad\forall f\in
{\mathcal{M}}_{p}^{w}(\Omega)
\cap
L^{p}(\Omega)\,,
\]
is a norm on ${\mathcal{M}}_{p}^{w}(\Omega)
\cap
L^{p}(\Omega)$ and 
$({\mathcal{M}}_{p}^{w}(\Omega)
\cap
L^{p}(\Omega),\|\cdot\|_{
{\mathcal{M}}_{p}^{w}(\Omega)
\cap
L^{p}(\Omega)})$ is a Banach space.
\item[(iii)] ${\mathcal{M}}_{p}^{w,0}(\Omega)
\cap
L^{p}(\Omega)$ is a Banach subspace of ${\mathcal{M}}_{p}^{w}(\Omega)
\cap
L^{p}(\Omega)$. 
\item[(iv)] Let $\eta_{w}>0$. Then ${\mathcal{M}}_{p}^{w}(\Omega)
\cap
L^{p}(\Omega)={\mathcal{M}}_{p}^{w}(\Omega)$ and ${\mathcal{M}}_{p}^{w,0}(\Omega)
\cap
L^{p}(\Omega)={\mathcal{M}}_{p}^{w,0}(\Omega)$ both algebraically and topologically.
\end{enumerate}
\end{prop}
{\bf Proof.} Since ${\mathcal{M}}_{p}^{w,0}(\Omega)$ is a closed subspace of ${\mathcal{M}}_{p}^{w}(\Omega)$, 
${\mathcal{M}}_{p}^{w,0}(\Omega)$  is a Banach space. Since the intersection of Banach spaces with the norm of the maximum is well known to be a Banach space, 
${\mathcal{M}}_{p}^{w,0}(\Omega)
\cap
L^{p}(\Omega)$ and ${\mathcal{M}}_{p}^{w}(\Omega)
\cap
L^{p}(\Omega)$ are Banach spaces. Hence, statements (ii) and (iii)  follow.

Next we consider statement (iv). If $\eta_{w}>0$, we already know that ${\mathcal{M}}_{p}^{w}(\Omega)
$ is continuously embedded into $L^{p}(\Omega)$, and we have
\[
\|f\|_{
{\mathcal{M}}_{p}^{w}(\Omega)}
\leq
\max\{
\|f\|_{
{\mathcal{M}}_{p}^{w}(\Omega)}
,
\|f\|_{L^{p}(\Omega)}
\}
\leq
\max\{1,\eta_{w}^{-1}\}
\|f\|_{
{\mathcal{M}}_{p}^{w}(\Omega)}
\]
for all $f\in {\mathcal{M}}_{p}^{w}(\Omega)
\cap
L^{p}(\Omega)$. Hence, ${\mathcal{M}}_{p}^{w}(\Omega)
\cap
L^{p}(\Omega)={\mathcal{M}}_{p}^{w}(\Omega)
$ both algebraically and topologically. Since the norm of 
${\mathcal{M}}_{p}^{w,0}(\Omega)$ is just the restriction of the norm of ${\mathcal{M}}_{p}^{w}(\Omega)$, the norm of ${\mathcal{M}}_{p}^{w,0}(\Omega)\cap
L^{p}(\Omega)$ is just the restriction of the norm of ${\mathcal{M}}_{p}^{w}(\Omega)\cap
L^{p}(\Omega)$, and   we have ${\mathcal{M}}_{p}^{w,0}(\Omega)
\cap
L^{p}(\Omega)={\mathcal{M}}_{p}^{w,0}(\Omega)
$ both algebraically and topologically.\hfill  $\Box$

\vspace{\baselineskip}

\section{A multiplication Theorem in Morrey spaces}

\begin{thm}
\label{mulgm}
Let $\Omega$ be an open subset of ${\mathbb{R}}^{n}$. Let $p_{1},p_{2}, p\in[1,+\infty]$ be such that $(1/p_{1})+(1/p_{2})=(1/p )$.
 Let $v_{1}$, $v_{2}$ be   functions from $]0,+\infty[$ to $[0,+\infty[$.   
Assume that there exist  $r_1$, $r_2\in]0,+\infty[$ such that $v_1(r_1)\neq 0\neq v_2(r_2)$. 
 
If $(f,g)\in {\mathcal{M}}^{v_{1}}_{p_{1}}(\Omega)\times
{\mathcal{M}}^{v_{2}}_{p_{2}}(\Omega)$, then $fg\in {\mathcal{M}}^{v_{1}v_{2}}_{p }(\Omega)$ and
\begin{eqnarray}
\label{mulgm1}
| fg |_{\rho,v_{1}v_{2},p ,\Omega}
&\leq&
| f |_{\rho,v_{1},p_{1},\Omega}
| g |_{\rho,v_{2},p_{2},\Omega}
\qquad\forall \rho\in]0,+\infty]\,,
\\
\label{mulgm2}
\|fg\|_{{\mathcal{M}}^{v_{1}v_{2} }_{p }(\Omega)} 
&\leq&
\|f\|_{{\mathcal{M}}^{v_{1}}_{p_{1}}(\Omega)} 
\|g\|_{{\mathcal{M}}^{v_{2}}_{p_{2}}(\Omega)} \,.
\end{eqnarray}
In particular, if either
\[
(f,g)\in {\mathcal{M}}^{v_{1},0}_{p_{1}}(\Omega)\times
{\mathcal{M}}^{v_{2} }_{p_{2}}(\Omega) 
\quad\text{or}\quad
 (f,g)\in {\mathcal{M}}^{v_{1} }_{p_{1}}(\Omega)\times
{\mathcal{M}}^{v_{2},0}_{p_{2}}(\Omega)\,,
\]
then $fg\in {\mathcal{M}}^{v_{1}v_{2},0}_{p }(\Omega)$.
\end{thm}
{\bf Proof.} Let $\rho\in]0,+\infty]$. If $(x,r)\in\Omega\times]0,\rho[$, then the H\"{o}lder inequality implies that
\begin{eqnarray*}
\lefteqn{v_{1}(r)v_{2}(r)\|fg\|_{L^{p}({\mathbb{B}}_{n}(x,r)\cap\Omega )}
}
\\
\nonumber
&&\qquad
\leq
v_{1}(r)\|f\|_{L^{p_{1}}({\mathbb{B}}_{n}(x,r)\cap\Omega )}
v_{2}(r)\|g\|_{L^{p_{2}}({\mathbb{B}}_{n}(x,r)\cap\Omega )}
\\
\nonumber
&&\qquad
\leq | f |_{\rho,v_{1},p_{1},\Omega}
| g |_{\rho,v_{2},p_{2},\Omega}\,.
\end{eqnarray*}
Then by taking the supremum in $(x,r)\in\Omega\times]0,\rho[$ in the left and right hand side, we obtain (\ref{mulgm1}). Inequality 
(\ref{mulgm2}) is an immediate consequence of inequality (\ref{mulgm1}). \hfill  $\Box$

\vspace{\baselineskip}

\begin{corol}
\label{mulgmi}
Let $\Omega$ be an open subset of ${\mathbb{R}}^{n}$. Let $p \in[1,+\infty]$. Let $w$ be  a  function  from $]0,+\infty[$ to $[0,+\infty[$. Assume that there exists  $r_0  \in]0,+\infty[$ such that $w(r_0)\neq 0$. 
 
If $(f,g)\in {\mathcal{M}}^{w }_{p }(\Omega)\times
L^{\infty}(\Omega)$, then $fg\in {\mathcal{M}}^{w }_{p }(\Omega)$ and
\begin{eqnarray}
\label{mulgmi1}
| fg |_{\rho,w  ,p  ,\Omega}
&\leq&
| f |_{\rho,w ,p ,\Omega}
\|g\|_{L^{\infty}(\Omega)}
\qquad\forall \rho\in]0,+\infty]\,,
\\
\label{mulgmi2}
\|fg\|_{{\mathcal{M}}^{w   }_{p }(\Omega)} 
&\leq&
\|f\|_{{\mathcal{M}}^{w }_{p }(\Omega)} 
\|g\|_{L^{\infty}(\Omega)} \,.
\end{eqnarray}
In particular, if $(f,g)\in {\mathcal{M}}^{w ,0}_{p }(\Omega)\times
L^{\infty}(\Omega)$, then $fg\in {\mathcal{M}}^{w ,0}_{p }(\Omega)$.
\end{corol}
{\bf Proof.} By Proposition \ref{prop:mo=lp}, we have ${\mathcal{M}}^{v_{2} ,0}_{p_{2}}(\Omega)=L^\infty(\Omega)$ when $p_{2}=+\infty$ and  $v_{2}=1$.
Then by taking $p=p_{1}$, $p_{2}=+\infty$, $v_{1}=w$,
 $v_{2}=1$ in the multiplication Theorem \ref{mulgm}, and by the obvious inequality
 \[
 | g |_{\rho,1,\infty,\Omega}\leq \|g\|_{L^{\infty}(\Omega)}
 \qquad\forall \rho\in]0,+\infty]\,,
 \]
we obtain inequalities (\ref{mulgmi1}) and (\ref{mulgmi2}), which imply the validity of the statement.
\hfill  $\Box$

\vspace{\baselineskip}

Then we can deduce the following immediate corollaries of the previous theorem and of its corollary.
\begin{corol}
\label{mulm}
Let $\Omega$ be an open subset of ${\mathbb{R}}^{n}$. Let $p_{1},p_{2},p\in[1,+\infty]$ be such that $(1/p_{1})+(1/p_{2})=(1/p )$. Let $\lambda_{1},\lambda_{2}\in[0,+\infty[$,   $\lambda=\lambda_{1}+\lambda_{2}$. Then the following statements hold.
\begin{enumerate}
\item[(i)] Let $\rho\in]0,+\infty[$. The pointwise multiplication is bilinear and continuous from ${\mathcal{M}}^{r^{-\lambda_{1}},\rho}_{p_{1}}(\Omega)\times
{\mathcal{M}}^{r^{-\lambda_{2}},\rho }_{p_{2}}(\Omega)$ to ${\mathcal{M}}^{r^{-\lambda},\rho }_{p }(\Omega)$ and maps
${\mathcal{M}}^{r^{-\lambda_{1}},\rho,0}_{p_{1}}(\Omega)\times
{\mathcal{M}}^{r^{-\lambda_{2}},\rho }_{p_{2}}(\Omega)$ to ${\mathcal{M}}^{r^{-\lambda},\rho,0 }_{p }(\Omega)$ and ${\mathcal{M}}^{r^{-\lambda_{1} },\rho}_{p_{1}}(\Omega)\times
{\mathcal{M}}^{r^{-\lambda_{2}},\rho,0}_{p_{2}}(\Omega)$ to ${\mathcal{M}}^{r^{-\lambda},\rho,0 }_{p }(\Omega)$.
\item[(ii)] The pointwise multiplication is bilinear and continuous from ${\mathcal{M}}^{r^{-\lambda_{1}} }_{p_{1}}(\Omega)\times
{\mathcal{M}}^{r^{-\lambda_{2}}  }_{p_{2}}(\Omega)$ to ${\mathcal{M}}^{r^{-\lambda}  }_{p }(\Omega)$ and maps
${\mathcal{M}}^{r^{-\lambda_{1}} ,0}_{p_{1}}(\Omega)\times
{\mathcal{M}}^{r^{-\lambda_{2}}  }_{p_{2}}(\Omega)$ to ${\mathcal{M}}^{r^{-\lambda} ,0 }_{p }(\Omega)$ and ${\mathcal{M}}^{r^{-\lambda_{1} } }_{p_{1}}(\Omega)\times
{\mathcal{M}}^{r^{-\lambda_{2}} ,0}_{p_{2}}(\Omega)$ to ${\mathcal{M}}^{r^{-\lambda} ,0 }_{p }(\Omega)$.
\item[(iii)] The pointwise multiplication is bilinear and continuous from $M^{\lambda_{1}}_{p_{1}}(\Omega)\times
M^{\lambda_{2}}_{p_{2}}(\Omega)$ to $M^{\lambda }_{p }(\Omega)$ and maps
$M^{\lambda_{1},0}_{p_{1}}(\Omega)\times
M^{\lambda_{2} }_{p_{2}}(\Omega)$ to $M^{\lambda,0 }_{p }(\Omega)$ and $M^{\lambda_{1} }_{p_{1}}(\Omega)\times
M^{\lambda_{2},0}_{p_{2}}(\Omega)$ to $M^{\lambda,0 }_{p }(\Omega)$.
\end{enumerate}
\end{corol}
{\bf Proof.} (i) It suffices to take $v_{1}=w_{\lambda_{1},\rho}$ and 
$v_{2}=w_{\lambda_{2},\rho}$ and apply the previous multiplication  theorem. (ii) It suffices to take $v_{1}=r^{-\lambda_{1}}$ and 
$v_{2}=r^{-\lambda_{2}}$  and apply the previous multiplication theorem. (iii) It suffices to take $v_{1}=w_{\lambda_{1}}$ and 
$v_{2}=w_{\lambda_{2}}$  and apply the previous multiplication theorem. \hfill  $\Box$ 

\vspace{\baselineskip}

\begin{corol}
\label{mulmi}
Let $\Omega$ be an open subset of ${\mathbb{R}}^{n}$. Let $p \in[1,+\infty]$. Let $\lambda\in[0,+\infty[$.  Then the following statements hold.
\begin{enumerate}
\item[(i)] Let $\rho\in]0,+\infty[$. The pointwise multiplication is bilinear and continuous from ${\mathcal{M}}^{r^{-\lambda},\rho }_{p }(\Omega)\times
L^{\infty}(\Omega)$ to ${\mathcal{M}}^{r^{-\lambda},\rho }_{p }(\Omega)$ and maps
${\mathcal{M}}^{r^{-\lambda},\rho,0}_{p}(\Omega)\times
L^{\infty}(\Omega) $ to ${\mathcal{M}}^{r^{-\lambda},\rho,0 }_{p }(\Omega)$.
\item[(ii)] The pointwise multiplication is bilinear and continuous from ${\mathcal{M}}^{r^{-\lambda}  }_{p }(\Omega)\times
L^{\infty}(\Omega)$ to ${\mathcal{M}}^{r^{-\lambda}  }_{p }(\Omega)$ and maps
${\mathcal{M}}^{r^{-\lambda} ,0}_{p}(\Omega)\times
L^{\infty}(\Omega) $ to ${\mathcal{M}}^{r^{-\lambda} ,0 }_{p }(\Omega)$.
\item[(iii)] The pointwise multiplication is bilinear and continuous from $M^{\lambda }_{p }(\Omega)\times
L^{\infty}(\Omega)$ to $M^{\lambda }_{p }(\Omega)$ and maps
$M^{\lambda,0}_{p}(\Omega)\times
L^{\infty}(\Omega) $ to $M^{\lambda,0 }_{p }(\Omega)$.
\end{enumerate}
\end{corol}
 
\section{An embedding theorem for generalized Morrey spaces}

We first prove the following elementary embedding theorem.
\begin{thm}
\label{thm:emgm}
Let $\Omega$ be an open subset of ${\mathbb{R}}^{n}$. Let $p,q\in [1,+\infty]$, $p\geq q$.  Let $v_{1}$, $v_{2}$ be   functions from $]0,+\infty[$ to $[0,+\infty[$. 
Assume that there exists  $r_0 \in]0,+\infty[$ such that $v_1(r_0)v_2(r_0)\neq 0$. Let
\begin{eqnarray*}
I (p,q,\rho)&\equiv&\left\{
c\in [0,+\infty[:\,   v_2(r)r^{n/q}
\leq c v_1(r)r^{n/p}  \ \forall r\in]0,\rho[
\right\}  \,,
\\
\iota  (p,q,\rho)&\equiv &\inf I (p,q,\rho)\qquad\forall \rho\in]0,+\infty]\,.
\end{eqnarray*}
If $\rho\in]0,+\infty]$ and 
\begin{equation}\label{thm:emgm1}
I (p,q,\rho)\neq\emptyset\,,
\end{equation}
then
\begin{equation}\label{thm:emgm1a}
|f|_{\rho,v_2,q,\Omega}
\leq
\omega_n^{\frac{1}{q}-\frac{1}{p}}\iota (p,q,\rho)
|f|_{\rho,v_1,p,\Omega}\qquad\forall f\in {\mathcal{M}}^{v_1}_p(\Omega)\,.
\end{equation}
In particular, if $I (p,q,+\infty)\neq\emptyset$, then
  ${\mathcal{M}}^{v_{1}}_{p }(\Omega)$ is continuously embedded into
${\mathcal{M}}^{v_{2}}_{q }(\Omega)$ and ${\mathcal{M}}^{v_{1},0}_{p }(\Omega)$ is continuously embedded into
${\mathcal{M}}^{v_{2},0}_{q }(\Omega)$. 
If we further assume that
\begin{equation}
\label{thm:emgm2}
\lim_{\rho\to 0} \iota (p,q,\rho)=0\,,
\end{equation}
then ${\mathcal{M}}^{v_{1}}_{p }(\Omega)$ is continuously embedded into
${\mathcal{M}}^{v_{2},0}_{q }(\Omega)$. 
\end{thm}
{\bf Proof.} Let $f\in {\mathcal{M}}^{v_{1}}_{p }(\Omega)$. If $\rho\in]0,+\infty]$, $I (p,q,\rho)\neq \emptyset$, $c\in I (p,q,\rho)$, then the H\"{o}lder inequality implies that
\begin{eqnarray}
\label{emgm4}
\lefteqn{
v_{2}(r)\|f\|_{L^{q}({\mathbb{B}}_{n}(x,r)\cap\Omega)}
}
\\
\nonumber
&& 
\leq 
v_{2}(r)  m_{n}({\mathbb{B}}_{n}(x,r)\cap\Omega)^{\frac{1}{q}-\frac{1}{p}}\|f\|_{L^{p}({\mathbb{B}}_{n}(x,r)\cap\Omega)}   
\\
\nonumber
&& 
\leq 
v_{2}(r)   \omega_n^{\frac{1}{q}-\frac{1}{p}}r^{\frac{n}{q}-\frac{n}{p}}\|f\|_{L^{p}({\mathbb{B}}_{n}(x,r)\cap\Omega)}   
\leq 
cv_{1}(r)   \omega_n^{\frac{1}{q}-\frac{1}{p}} \|f\|_{L^{p}({\mathbb{B}}_{n}(x,r)\cap\Omega)}   
\end{eqnarray}
for all $(x,r)\in\Omega\times]0,\rho[$. Then by taking first the supremum in $(x,r)\in\Omega\times]0,\rho[$, and then the infimum in $c\in I(p,q,\rho)$, we obtain inequality (\ref{thm:emgm1a}). 

In particular, if $I (p,q,+\infty)\neq\emptyset$, then
inequality (\ref{thm:emgm1a}) with $\rho=+\infty$ implies that ${\mathcal{M}}^{v_{1}}_{p }(\Omega)$ is continuously embedded into
${\mathcal{M}}^{v_{2}}_{q }(\Omega)$. Moreover, 
\[
I(p,q,\rho)\supseteq  I(p,q,+\infty)\qquad\forall \rho\in]0,+\infty]
\]
and accordingly
\[
\iota (p,q,\rho)\leq \iota (p,q,+\infty)\qquad\forall \rho\in]0,+\infty]\,.
\]
Hence, inequality  (\ref{thm:emgm1a}) implies that
\[
|f|_{\rho,v_2,q,\Omega}
\leq
\omega_n^{\frac{1}{q}-\frac{1}{p}}\iota (p,q,+\infty)
|f|_{\rho,v_1,p,\Omega}\qquad\forall f\in {\mathcal{M}}^{v_1}_p(\Omega)\,.
\]
for all $\rho\in]0,+\infty]$. Thus if $f\in{\mathcal{M}}^{v_{1},0}_{p }(\Omega)$, we have $f\in {\mathcal{M}}^{v_{2},0}_{q }(\Omega)$. 
 If we further assume that the limiting condition  (\ref{thm:emgm2}) holds and that $f\in {\mathcal{M}}^{v_{1}}_{p }(\Omega)$, then inequality (\ref{thm:emgm1a})  implies that 
\[
\lim_{\rho\to 0}|f|_{\rho,v_2,q,\Omega}=0
\]
 and accordingly that   $f\in {\mathcal{M}}^{v_{2},0}_{q }(\Omega)$.\hfill  $\Box$

\vspace{\baselineskip}

\begin{rem}\label{rem:emgm} 
{\em Under the assumptions of the previous statement, if we further assume that
\[
v_1(r)\neq 0\qquad\forall r\in]0,+\infty[\,,
\]
 then 
 \[
  \iota (p,q,\rho)=\sup_{r\in]0,\rho[}\frac{v_{2}(r)r^{n/q}}{v_{1}(r)r^{n/p}}\qquad\forall \rho\in]0,+\infty]\,.
\]
 In particular,  condition
 \begin{equation}
 \label{emgm1}
\sup_{r\in]0,+\infty[}\frac{v_{2}(r)r^{n/q}}{v_{1}(r)r^{n/p}}<+\infty\,,
\end{equation}
 implies that 
 \[
 \iota (p,q,+\infty)=\sup_{r\in]0,+\infty[}\frac{v_{2}(r)r^{n/q}}{v_{1}(r)r^{n/p}}\in I(p,q,+\infty)\neq\emptyset\,,
 \]
 and condition 
\begin{equation}
\label{emgm2}
 \lim_{\rho\to 0} \sup_{r\in]0,\rho[}\frac{v_{2}(r)r^{n/q}}{v_{1}(r)r^{n/p}}=0\,,
\end{equation}
implies that
\[
 \lim_{\rho\to 0} \iota (p,q,\rho)=0\,.
 \]
}\end{rem}
 Then we have the following corollary for the Morrey spaces with weights $w_{\lambda,\rho}$, $r^{-\lambda}$, $w_\lambda$.

\begin{corol}\label{corol:emgm}
 Let $\Omega$ be an open subset of ${\mathbb{R}}^{n}$. Let $p,q \in[1,+\infty]$. Let $\lambda\in[0,+\infty[$.  Then the following statements hold.
 \begin{enumerate}
\item[(i)] Let $\tilde{\rho}\in]0,+\infty[$. If
\begin{equation}\label{corol:emgm1}
p
\geq q\,,\qquad \left(\lambda- (n/p)\right)
\geq\left(\nu- (n/q)\right)
\end{equation}
then $ {\mathcal{M}}^{r^{-\lambda},\tilde{\rho}}_{p }(\Omega)$ is continuously embedded into
${\mathcal{M}}^{r^{-\nu},\tilde{\rho}}_{q }(\Omega)$ and $ {\mathcal{M}}^{r^{-\lambda},\tilde{\rho},0}_{p }(\Omega)$ is continuously embedded into
${\mathcal{M}}^{r^{-\nu},\tilde{\rho},0}_{q }(\Omega)$.  If we further assume that the second inequaluity of (\ref{corol:emgm1}) is strict, then 
$ {\mathcal{M}}^{r^{-\lambda},\tilde{\rho}}_{p }(\Omega)$ is continuously embedded into
${\mathcal{M}}^{r^{-\nu},\tilde{\rho},0}_{q }(\Omega)$.
\item[(ii)] If
\begin{equation}\label{corol:emgm2}
p
\geq q\,,\qquad \left(\lambda- (n/p)\right)
=\left(\nu- (n/q)\right)
\end{equation}
then $ {\mathcal{M}}^{r^{-\lambda}}_{p }(\Omega)$ is continuously embedded into
${\mathcal{M}}^{r^{-\nu} }_{q }(\Omega)$ and $ {\mathcal{M}}^{r^{-\lambda} ,0}_{p }(\Omega)$ is continuously embedded into
${\mathcal{M}}^{r^{-\nu} ,0}_{q }(\Omega)$.
\end{enumerate}
\item[(iii)] If 
\begin{equation}
\label{emm1}
p= q\,,\qquad \lambda>\nu\,,
\end{equation}
then $M^{\lambda}_{p }(\Omega)$ is continuously embedded into
$M^{\nu,0}_{q }(\Omega)$. 
 \end{corol}
 {\bf Proof.} (i) By definition of $w_{\nu,\tilde{\rho}}$, $w_{\lambda,\tilde{\rho}}$, we have
 \[
 w_{\nu,\tilde{\rho}}(r)r^{n/q}=r^{-\nu+(n/q)}\leq r^{-\lambda+(n/p)}=w_{\lambda,\tilde{\rho}}(r)r^{n/p}
 \]
for  $r\in]0,1]$, $r<\tilde{\rho}$,
\begin{eqnarray*}
\lefteqn{
w_{\nu,\tilde{\rho}}(r)r^{n/q}
=r^{-\nu+(n/q)}=r^{-\nu}r^{(n/q)}\leq \tilde{\rho}^{(n/q)}=(\tilde{\rho}^{(n/q)}\tilde{\rho}^\lambda)\tilde{\rho}^{-\lambda}
}
\\ \nonumber
&&\qquad\qquad\qquad
\leq  
(\tilde{\rho}^{n/q}\tilde{\rho}^\lambda)r^{-\lambda}r^{n/p}
= (\tilde{\rho}^{n/q}\tilde{\rho}^\lambda)w_{\lambda,\tilde{\rho}}(r)r^{n/p}
\end{eqnarray*}
for  $r\in [1,+\infty[$, $r<\tilde{\rho}$,
\[
w_{\nu,\tilde{\rho}}(r)r^{n/q}=0=w_{\lambda,\tilde{\rho}}(r)r^{n/p}
\]
for $r\geq\tilde{\rho}$. In conclusion,
\[
\max\{1, (\tilde{\rho}^{n/q}\tilde{\rho}^\lambda),0\}\in I(p,q,+\infty)
\]
and thus Theorem \ref{thm:emgm} implies the validity of the first part of the statement. 
Now let $f\in  {\mathcal{M}}^{r^{-\lambda},\tilde{\rho},0}_{p }(\Omega)$
If $\eta<1$, the first of the above computation shows that
\[
\iota(p,q,\eta)\leq\sup_{r\in]0,\eta[}r^{-\nu+(n/q) -(-\lambda+(n/p))}
=\eta^{-\nu+(n/q) -(-\lambda+(n/p))}
\,.
\]
Thus  if we further assume that the second inequality of (\ref{corol:emgm1}) is strict, then we have
$\lim_{\eta\to 0}\iota(p,q,\eta)=0$ and again 
Theorem \ref{thm:emgm} implies that $ {\mathcal{M}}^{r^{-\lambda},\tilde{\rho}}_{p }(\Omega)$ is continuously embedded into
${\mathcal{M}}^{r^{-\nu},\tilde{\rho},0}_{q }(\Omega)$.

(ii) Since
\[
r^{-\nu}r^{n/q}
=r^{-\lambda}r^{n/p}\qquad\forall r\in]0,+\infty[\,,
\]
we have $1\in \iota(p,q,+\infty)$ and thus Theorem \ref{thm:emgm} implies the validity of the statement.

(iii) By definition of $w_{\lambda}$ and $w_{\nu}$ and by condition  $p=q$, we have
\[
\frac{w_{\nu}(r)r^{n/q}}{w_{\lambda}(r)r^{n/p}}=
r^{ \lambda-\nu}\qquad\forall r\in ]0,1]\,,
\qquad
\frac{w_{\nu}(r)r^{n/q}}{w_{\lambda}(r)r^{n/p}}=1
\qquad\forall r\in ]1,+\infty[\,.
\]
Since $p=q$ and $\lambda>\nu$, Theorem \ref{thm:emgm} and Remark \ref{rem:emgm}
imply  that (iii) holds true.\hfill  $\Box$ 

\vspace{\baselineskip}

\subsection{Relation between $M^{\lambda}_{p}(\Omega)$ and ${\mathcal{M}}^{r^{-\lambda}}_{p }(\Omega)$}
The following corollary clarifies the relation between $M^{\lambda}_{p}(\Omega)$ and ${\mathcal{M}}^{r^{-\lambda}}_{p }(\Omega)$. 
\begin{corol}
\label{mlpwp}
Let $\Omega$ be an open subset of ${\mathbb{R}}^{n}$. Let $p\in [1,+\infty]$. 
Let $\lambda\in [0,+\infty[$. Then the following statements hold.
\begin{enumerate}
\item[(i)] $M^{\lambda}_{p}(\Omega)=
{\mathcal{M}}^{r^{-\lambda}}_{p }(\Omega)\cap L^{p}(\Omega)$, and the corresponding norms coincide. Namely,
\[
\|f\|_{ M^{\lambda}_{p}(\Omega) }
=
\max\{
\|f\|_{ {\mathcal{M}}^{r^{-\lambda}}_{p }(\Omega) },
\|f\|_{ L^{p}(\Omega) }
\}
\qquad\forall f\in M^{\lambda}_{p}(\Omega)\,.
\]
\item[(ii)] If $\Omega$ is bounded, then
\[
\|f\|_{ L^{p}(\Omega)  }
\leq
{\mathrm{diam}}\,^{\lambda}(\Omega)
\|f\|_{{\mathcal{M}}^{r^{-\lambda}}_{p }(\Omega) }
\qquad\forall f\in {\mathcal{M}}^{r^{-\lambda}}_{p }(\Omega) \,,
\]
and ${\mathcal{M}}^{r^{-\lambda}}_{p }(\Omega)$ is continuously embedded into $L^{p}(\Omega) $.
\item[(iii)]  If $\Omega$ is bounded, then
\[
M^{\lambda}_{p}(\Omega)={\mathcal{M}}^{r^{-\lambda}}_{p }(\Omega)\,,
\]
both algebraically and topologically. Moreover,
\[
\|f\|_{ {\mathcal{M}}^{r^{-\lambda}}_{p }(\Omega) }
\leq
\|f\|_{M^{\lambda}_{p}(\Omega) }
\leq
\max\{
1,{\mathrm{diam}}\,^{\lambda}(\Omega)
\}\|f\|_{  {\mathcal{M}}^{r^{-\lambda}}_{p }(\Omega) }
\,,
\]
for all  $f\in
M^{\lambda}_{p}(\Omega)={\mathcal{M}}^{
r^{-\lambda}}_{p }(\Omega)$. 
\end{enumerate}
\end{corol}
{\bf Proof.} We first consider statement (i). Since
\[
\sup_{r\in]0,+\infty[}\frac{r^{-\lambda}r^{n/p}}{w_{\lambda}(r)r^{n/p}}
=\sup_{r\in]0,+\infty[}\frac{r^{-\lambda}}{w_{\lambda}(r)}=1<+\infty\,,
\]
Theorem \ref{thm:emgm}  and Remark \ref{rem:emgm} imply that
$M^{\lambda}_{p}(\Omega)=
{\mathcal{M}}^{
w_{\lambda}}_{p }(\Omega)$ is continuously embedded into ${\mathcal{M}}^{
r^{-\lambda}}_{p }(\Omega)$ and that
\[
\|f\|_{
{\mathcal{M}}^{
r^{-\lambda}}_{p }(\Omega)
}
\leq
\|f\|_{M^{\lambda}_{p}(\Omega)}
\qquad\forall f\in M^{\lambda}_{p}(\Omega)\,.
\]
By Corollary \ref{corol:molp}, we already know that
$M^{\lambda}_{p}(\Omega)$ is continuously embedded into $L^{p}(\Omega)$ and that
\[
\|f\|_{ L^{p}(\Omega) }\leq
\|f\|_{ M^{\lambda}_{p}(\Omega) }
\qquad\forall f \in  M^{\lambda}_{p}(\Omega)\,.
\]
Hence, $M^{\lambda}_{p}(\Omega)$ is continuously embedded into ${\mathcal{M}}^{r^{-\lambda}}_{p }(\Omega)\cap L^{p}(\Omega)$ and 
\[
\|f\|_{
{\mathcal{M}}^{r^{-\lambda}}_{p }(\Omega)\cap L^{p}(\Omega)
}
=
\max\{\|f\|_{{\mathcal{M}}^{r^{-\lambda}}_{p }(\Omega)},\|f\|_{L^{p}(\Omega)}
\}\leq
\|f\|_{M^{\lambda}_{p}(\Omega)}
\qquad\forall f\in M^{\lambda}_{p}(\Omega)\,.
\]
On the other hand, if $f\in {\mathcal{M}}^{r^{-\lambda}}_{p }(\Omega)\cap L^{p}(\Omega)$, then we have 
\begin{eqnarray*}
\lefteqn{
w_{\lambda}(r)\|f\|_{L^{p}(
{\mathbb{B}}_{n}(x,r)\cap\Omega)}
}
\\ \nonumber
&&\qquad
=r^{-\lambda}\|f\|_{L^{p}(
{\mathbb{B}}_{n}(x,r)\cap\Omega)}
\leq\|f\|_{{\mathcal{M}}^{r^{-\lambda}}_{p }(\Omega)}\qquad {\mathrm{if}}\ (x,r)\in \Omega\times]0,1]\,,
\\
\lefteqn{
w_{\lambda}(r)\|f\|_{L^{p}(
{\mathbb{B}}_{n}(x,r)\cap\Omega)}
}
\\ \nonumber
&&\qquad
= \|f\|_{L^{p}(
{\mathbb{B}}_{n}(x,r)\cap\Omega)}
\leq\|f\|_{L^{p}(\Omega)}\qquad {\mathrm{if}}\ (x,r)\in \Omega\times [1,+\infty[\,.
\end{eqnarray*}
Hence, $f\in M^{\lambda}_{p}(\Omega)$, 
\[
\|f\|_{  M^{\lambda}_{p}(\Omega)  }
\leq
\max\{
\|f\|_{{\mathcal{M}}^{r^{-\lambda}}_{p }(\Omega)}
,\|f\|_{L^{p}(\Omega)}
\}
=\|f\|_{{\mathcal{M}}^{r^{-\lambda}}_{p }(\Omega)
\cap L^{p}(\Omega)} 
\]
and statement (i) holds true.

We now consider statement (ii). If $\Omega $ is empty, then statement (ii) is obviously true. If $\Omega\neq\emptyset$, then there exists $\xi\in \Omega$ and we have $\Omega\subseteq 
{\mathbb{B}}_{n}(\xi,{\mathrm{diam}}\,(\Omega))$ and 
\begin{eqnarray*}
\lefteqn{
\|f\|_{L^{p}(\Omega)}
=
\|f\|_{  
L^{p}({\mathbb{B}}_{n}(\xi,{\mathrm{diam}}\,(\Omega))\cap\Omega)}
}
\\ \nonumber
&&\qquad
\leq
{\mathrm{diam}}^{\lambda}(\Omega)
{\mathrm{diam}}^{-\lambda}(\Omega)
\|f\|_{  
L^{p}({\mathbb{B}}_{n}(\xi,{\mathrm{diam}}\,(\Omega))\cap\Omega)}
\\ \nonumber
&&\qquad
\leq
{\mathrm{diam}}^{\lambda}(\Omega)
\|f\|_{{\mathcal{M}}^{r^{-\lambda}}_{p }(\Omega)}
\qquad
\forall f\in
{\mathcal{M}}^{r^{-\lambda}}_{p }(\Omega)\,.
\end{eqnarray*}
Hence, statement (ii) holds true.\par

We now consider statement (iii). Let $f\in M^{\lambda}_{p}(\Omega)$. Then statement (i) implies that $f\in {\mathcal{M}}^{r^{-\lambda}}_{p }(\Omega)$ and that
\[
\|f\|_{ {\mathcal{M}}^{r^{-\lambda}}_{p }(\Omega) }
\leq
\|f\|_{ M^{\lambda}_{p}(\Omega) }
\,.
\]
Conversely, if $f\in {\mathcal{M}}^{r^{-\lambda}}_{p }(\Omega)$, then statement (ii) implies that $f\in 
L^{p}(\Omega)$ and that
\[
\|f\|_{ L^{p}(\Omega)  }
\leq {\mathrm{diam}}\,^{\lambda}(\Omega)
\|f\|_{{\mathcal{M}}^{r^{-\lambda}}_{p }(\Omega) }\,.
\]
Then we have $f\in {\mathcal{M}}^{r^{-\lambda}}_{p }(\Omega)\cap L^{p}(\Omega)  $ and 
\[
\max\{
\|f\|_{{\mathcal{M}}^{r^{-\lambda}}_{p }(\Omega) }
,\|f\|_{ L^{p}(\Omega)  }
\}
\leq 
\max\{1,{\mathrm{diam}}\,^{\lambda}(\Omega)\}
\|f\|_{{\mathcal{M}}^{r^{-\lambda}}_{p }(\Omega) }\,.
\]
By statement (i), $M^{\lambda}_{p}(\Omega)=
{\mathcal{M}}^{r^{-\lambda}}_{p }(\Omega)\cap L^{p}(\Omega)$, and the corresponding norms coincide. Hence, statement (iii) follows. 
\hfill  $\Box$ 

\vspace{\baselineskip}

\begin{rem}{\em 
 By Corollaries \ref{corol:prembd}, \ref{mlpwp}, if $\Omega$ is a bounded open subset of 
${\mathbb{R}}^{n}$, $p\in [1,+\infty]$, $\lambda\in [0,+\infty[$, $\rho\in [0,+\infty[$, then
\[
{\mathcal{M}}_{p}^{r^{-\lambda},\rho }(\Omega) 
={\mathcal{M}}_{p}^{w_\lambda }(\Omega)(\equiv M^{\lambda}_{p}(\Omega))={\mathcal{M}}^{r^{-\lambda}}_{p }(\Omega)\,.
\]}
\end{rem}

\section{An embedding theorem for generalized Morrey spaces in case the domain has finite volume}
 Next we see that in case $m_n(\Omega)$ is finite, we can  weaken the assumptions of Theorem \ref{thm:emgm}. 
\begin{thm}\label{thm:emgmf}
 Let $\Omega$ be an open subset of ${\mathbb{R}}^{n}$, $m_n(\Omega)<+\infty$. Let $p,q\in [1,+\infty]$, $p\geq q$.  Let $v_{1}$, $v_{2}$ be   functions from $]0,+\infty[$ to $[0,+\infty[$. 
Assume that there exist  $r_1$, $r_2\in]0,+\infty[$ such that $v_1(r_1)\neq 0\neq v_2(r_2)$. Let
\begin{eqnarray*}
\lefteqn{
J (p,q,\rho)\equiv\bigl\{\bigr.
c\in [0,+\infty[:\,   v_2(r)r^{n/q}
\leq c v_1(r)r^{n/p}  \ \forall r\in]0,1[, r<\rho\,,
}
\\
&&\qquad\qquad\qquad\qquad\qquad\qquad\qquad
\ v_2(r) \leq c v_1(r)\ \forall r\in [1,+\infty[, \ r<\rho
\bigl.\bigr\}  \,,
\\
\lefteqn{
\eta  (p,q,\rho)\equiv \inf J (p,q,\rho)\qquad\forall \rho\in]0,+\infty]\,,}
\end{eqnarray*}
(see Theorem \ref{thm:emgm}.) Then we have
\begin{equation}\label{thm:emgmf1a}
J (p,q,\rho)=I(p,q,\rho)\,,\quad\eta  (p,q,\rho)=\iota  (p,q,\rho)\qquad\forall \rho\in]0,1[\,.
\end{equation}
If  
\begin{equation}\label{thm:emgmf1}
J (p,q,+\infty)\neq\emptyset\,,
\end{equation}
then ${\mathcal{M}}^{v_{1}}_{p }(\Omega)$ is continuously embedded into
${\mathcal{M}}^{v_{2}}_{q }(\Omega)$ and ${\mathcal{M}}^{v_{1},0}_{p }(\Omega)$ is continuously embedded into
${\mathcal{M}}^{v_{2},0}_{q }(\Omega)$. If we further assume that (\ref{thm:emgm2}) holds, then ${\mathcal{M}}^{v_{1}}_{p }(\Omega)$ is continuously embedded into
${\mathcal{M}}^{v_{2},0}_{q }(\Omega)$. 
\end{thm}
{\bf Proof.} We first note that the equalities in (\ref{thm:emgmf1a}) are an immediate consequence of the definition of $J (p,q,\rho)$ and of $\eta (p,q,\rho)$. Now
let $f\in {\mathcal{M}}^{v_{1}}_{p }(\Omega)$. Since $p\geq q$, then the H\"{o}lder inequality implies that
\begin{equation}\label{thm:emgmf4}
v_{2}(r)\|f\|_{L^{q}({\mathbb{B}}_{n}(x,r)\cap\Omega)}
\leq 
 \|f\|_{L^{p}({\mathbb{B}}_{n}(x,r)\cap\Omega)}    v_{2}(r) m_{n}({\mathbb{B}}_{n}(x,r)\cap\Omega)^{\frac{1}{q}-\frac{1}{p}}
\end{equation}
for all $(x,r)\in\Omega\times]0,+\infty[$. Next we assume that $J (p,q,+\infty)\neq\emptyset$ and we prove the first embedding of the last part of the statement. If $c\in J (p,q,+\infty)$,   we have
\begin{eqnarray}
\label{thm:emgmf6}
\lefteqn{
v_{2}(r)\|f\|_{L^{q}({\mathbb{B}}_{n}(x,r)\cap\Omega)}
}
\\
\nonumber
&&\qquad
\leq v_{2}(r)\omega_{n}^{\frac{1}{q}-\frac{1}{p}}r^{\frac{n}{q}-\frac{n}{p}}\|f\|_{L^{p}({\mathbb{B}}_{n}(x,r)\cap\Omega)} 
\leq
\omega_{n}^{\frac{1}{q}-\frac{1}{p}}cv_1(r)
\|f\|_{L^{p}({\mathbb{B}}_{n}(x,r)\cap\Omega)}   \,,
\end{eqnarray}
for all $(x,r)\in \Omega\times]0,1[$. Since $m_{n}(\Omega)<+\infty$, we have 
\[
m_{n}({\mathbb{B}}_{n}(x,r)\cap\Omega)^{\frac{1}{q}-\frac{1}{p}}
\leq
m_{n}(\Omega)^{\frac{1}{q}-\frac{1}{p}}
\qquad\forall r\in[1,+\infty[\,,
\]
Then inequality (\ref{thm:emgmf4}) implies that
\begin{eqnarray}
\label{thm:emgmf7}
\lefteqn{
v_{2}(r)\|f\|_{L^{q}({\mathbb{B}}_{n}(x,r)\cap\Omega)}
}
\\
\nonumber
&&\qquad
\leq
m_{n}( \Omega)^{\frac{1}{q}-\frac{1}{p}}
c
v_{1}(r)\|f\|_{L^{p}({\mathbb{B}}_{n}(x,r)\cap\Omega)}   \,,
\end{eqnarray}
for all $(x,r)\in \Omega\times[1,+\infty[$. By combining inequality
(\ref{thm:emgmf6})  and inequality (\ref{thm:emgmf7}), we deduce that
\[
\|f\|_{  {\mathcal{M}}^{v_{2}}_{q }(\Omega)  }
\leq
\max\left\{
\omega_{n}^{\frac{1}{q}-\frac{1}{p}}
 ,
m_{n}(\Omega)^{\frac{1}{q}-\frac{1}{p}} 
\right\}c
\|f\|_{  {\mathcal{M}}^{v_{1}}_{p }(\Omega)  }\,,
\]
an inequality which implies that $ {\mathcal{M}}^{v_{1}}_{p }(\Omega) $ is continuously embedded into the space ${\mathcal{M}}^{v_{2}}_{q }(\Omega)$. Since $J (p,q,+\infty)\neq\emptyset$, we have $I (p,q,\rho)=J (p,q,\rho)\neq\emptyset$ for all $\rho\in ]0,1[$ (cf.~(\ref{thm:emgmf1a})). Then Theorem   
\ref{thm:emgm} implies the validity of inequality (\ref{thm:emgm1a}). Hence, $f\in {\mathcal{M}}^{v_{2},0}_{q }(\Omega)$ whenever
$f\in {\mathcal{M}}^{v_{1},0}_{p }(\Omega) $. Moreover, if further assume that (\ref{thm:emgm2}) holds, then ${\mathcal{M}}^{v_{1}}_{p }(\Omega)$ is contained in 
${\mathcal{M}}^{v_{2},0}_{q }(\Omega)$. \hfill  $\Box$ 

\vspace{\baselineskip}

\begin{rem}\label{rem:emgmf} 
{\em Under the assumptions of the previous statement, if we further assume that
\[
v_1(r)\neq 0\qquad\forall r\in]0,+\infty[\,,
\]
 then condition
 \begin{equation}
 \label{emgmf1}
\sup_{r\in]0,1[}\frac{v_{2}(r)r^{n/q}}{v_{1}(r)r^{n/p}}<+\infty\,,
\qquad
\sup_{r\in [1,+\infty[}\frac{v_{2}(r)}{v_{1}(r)}<+\infty\,.
\,,
\end{equation}
 implies that $J(p,q,+\infty)\neq\emptyset$.}\end{rem}

Then we have the following corollary for the generalized Morrey spaces with the weights $r^{-\lambda}$ and  $w_\lambda$.
 \begin{corol}
\label{corol:emm}
Let $\Omega$ be an open subset of ${\mathbb{R}}^{n}$, $m_n(\Omega)<+\infty$. Let $p,q\in [1,+\infty]$.  Let $\lambda,\nu\in [0,+\infty[$. Then the following statements hold.
 \begin{enumerate}
\item[(i)] If
\begin{equation}\label{corol:emgmf2}
p
\geq q\,,\qquad \left(\lambda- (n/p)\right)
\geq\left(\nu- (n/q)\right)\,,\qquad\lambda\leq\nu\,,
\end{equation}
then $ {\mathcal{M}}^{r^{-\lambda}}_{p }(\Omega)$ is continuously embedded into
${\mathcal{M}}^{r^{-\nu} }_{q }(\Omega)$ and $ {\mathcal{M}}^{r^{-\lambda} ,0}_{p }(\Omega)$ is continuously embedded into
${\mathcal{M}}^{r^{-\nu} ,0}_{q }(\Omega)$. If we further assume that the second inequality of (\ref{corol:emgmf2}) is strict, then 
$ {\mathcal{M}}^{r^{-\lambda} }_{p }(\Omega)$ is continuously embedded into
${\mathcal{M}}^{r^{-\nu} ,0}_{q }(\Omega)$.

\item[(ii)] If
\begin{equation}
\label{emm3}
p\geq q\,,\qquad \left(\lambda- (n/p)\right)
\geq\left(\nu- (n/q)\right)\,,
\end{equation}
then $M^{\lambda}_{p }(\Omega)$ is continuously embedded into
$M^{\nu}_{q }(\Omega)$.  If we further assume that the second inequality of (\ref{emm3}) is strict, then 
$M^{\lambda}_{p }(\Omega)$ is continuously embedded into
$M^{\nu,0}_{q }(\Omega)$. 
\end{enumerate}
\end{corol}
 {\bf Proof.} We first consider statement (i) and we note that
 \[
 \frac{r^{-\nu}r^{n/q}}{r^{-\lambda} r^{n/p}}=r^{ \left(\lambda- (n/p)\right)
-\left(\nu- (n/q)\right)}\quad\forall r\in ]0,1]\,,
\qquad
\frac{r^{-\nu} }{r^{-\lambda} }=r^{\lambda-\nu}
\quad\forall r\in ]1,+\infty[\,.
 \]
 Since $p\geq q$,  $\left(\lambda- (n/p)\right)
-\left(\nu- (n/q)\right)\geq 0$, $\lambda\leq\nu$, Theorem \ref{thm:emgmf} and Remark \ref{rem:emgmf} imply that (i) holds true.\par

 (ii) By definition of $w_{\lambda}$ and $w_{\nu}$, we have
\[
\frac{w_{\nu}(r)r^{n/q}}{w_{\lambda}(r)r^{n/p}}=
r^{ \left(\lambda- (n/p)\right)
-\left(\nu- (n/q)\right)}\quad\forall r\in ]0,1]\,,
\qquad
\frac{w_{\nu}(r) }{w_{\lambda}(r) }=1
\quad\forall r\in ]1,+\infty[\,.
\]
Since $p\geq q$ and $\left(\lambda- (n/p)\right)
-\left(\nu- (n/q)\right)\geq 0$, Theorem \ref{thm:emgmf} and Remark \ref{rem:emgmf} imply that (ii) holds true.
 \hfill  $\Box$

\vspace{\baselineskip}

 By exploiting   Theorem \ref{thm:emgmf}, we can also improve the conditions of Proposition \ref{prop:bdmo} that ensure that
a generalized Morrey space    contains the space of essentially bounded functions in case $\Omega$ has a finite measure. Indeed, we have the following statement. 
\begin{corol}
\label{bfgm}
Let $\Omega$ be an open subset of ${\mathbb{R}}^{n}$ of finite measure. Let $p\in [1,+\infty]$. Let $w$ be a   function from $]0,+\infty[$ to itself which satisfies the following conditions
\begin{equation}
\label{bfgm1}
\sup_{r\in]0,1]} w(r)r^{n/p} <+\infty\,,\qquad \sup_{r\in [1,+\infty[} w (r) <+\infty\,.
\end{equation}
Then $L^{\infty}(\Omega)$ is continuously embedded into ${\mathcal{M}}^{w}_{p}(\Omega)$. 

If we further assume that
\begin{equation}
\label{bfgm2}
\lim_{\rho\to 0}\sup_{r\in]0,\rho]} w(r)r^{n/p} =0\,,
\end{equation}
then $L^{\infty}(\Omega)$ is continuously embedded into ${\mathcal{M}}^{w,0}_{p}(\Omega)$. 
\end{corol}
{\bf Proof.}  Let $v(r)=1$ for all $r\in]0,+\infty[$. 
Then  
Corollary \ref{corol:mo=lp}  implies that
 $L^{\infty}(\Omega)={\mathcal{M}}^{v}_{\infty}(\Omega)$.
 Then inequalities
(\ref{bfgm1}) imply that the weights $v_{1}\equiv v$, 
$v_{2}\equiv w$ with $p$, $q$ replaced by $\infty$, $p$, respectively satisfy inequalities
(\ref{emgmf1}). Then Theorem \ref{thm:emgmf}  and Remark \ref{rem:emgmf} imply that
$L^{\infty}(\Omega)={\mathcal{M}}^{v}_{\infty}(\Omega)$ is continuously embedded into ${\mathcal{M}}^{w}_{p}(\Omega)$. If we further assume that (\ref{bfgm2}) holds, then Theorem \ref{thm:emgmf}  implies that
$L^{\infty}(\Omega)={\mathcal{M}}^{v}_{\infty}(\Omega)$ is continuously embedded into ${\mathcal{M}}^{w,0}_{p}(\Omega)$. 
\hfill  $\Box$

\vspace{\baselineskip}

By applying the previous corollary to the weights $w_{\lambda,\rho}$, with $\rho\in]0,+\infty[$, $r^{-\lambda}$, $w_\lambda$, we immediately deduce  the validity of the following statement. 
\begin{corol}
\label{bfm}
Let $\Omega$ be an open subset of ${\mathbb{R}}^{n}$ of finite measure. Let $p\in [1,+\infty]$. 
If $\lambda\in[0,n/p[$, then 
\begin{enumerate}
\item[(i)] $L^{\infty}(\Omega)$ is continuously embedded into 
${\mathcal{M}}^{r^{-\lambda},\rho,0}_{p}(\Omega)$ for all $\rho\in]0,+\infty[$. 
\item[(ii)] $L^{\infty}(\Omega)$ is continuously embedded into 
${\mathcal{M}}^{r^{-\lambda}, 0}_{p}(\Omega)$. 
\item[(iii)] $L^{\infty}(\Omega)$ is continuously embedded into 
$M^{\lambda,0}_{p}(\Omega)$. 
\end{enumerate}
\end{corol}

\vspace{\baselineskip}

\section{Mollifiers and  functions in generalized Morrey spaces}

If $\phi\in {\mathbb{R}}^{{\mathbb{R}}^{ n}}$ and $t\in]0,+\infty[$, we denote by $\phi_{t}(\cdot)$ the function from ${\mathbb{R}}^{ n}$ to ${\mathbb{R}}$ defined by
\begin{equation}
\label{phd}
\phi_{t}(x)\equiv t^{-n}\phi (x/t)\qquad\forall x\in {\mathbb{R}}^{ n}\,.
\end{equation}
By the formula of change of variables in integrals, we conclude that
\[
\int_{ {\mathbb{R}}^{ n} }\phi_{t}(x)\,dx=\int_{ {\mathbb{R}}^{ n} }\phi (x)\,dx\qquad\forall t\in]0,+\infty[\,,
\]
whenever $\phi\in L^{1}({\mathbb{R}}^{ n})$. We also note that if ${\mathrm{supp}}\,\phi\subseteq\overline{{\mathbb{B}}_n(0,1)}$, then
\[
{\mathrm{supp}}\,\phi_t\subseteq\overline{{\mathbb{B}}_n(0,t)}\qquad\forall t\in]0,+\infty[\,.
\]
Then we have the following result of approximation by convolution.
\begin{thm}
\label{apgm}
Let $\phi\in C_{c}^{\infty}({\mathbb{R}}^{ n})$, $\int_{ {\mathbb{R}}^{ n} }\phi\,dx=1$, ${\mathrm{supp}}\,\phi\subseteq\overline{{\mathbb{B}}_n(0,1)}$.  Let $w$ be a function from $]0,+\infty[$ to $[0,+\infty[$.  Assume that there exists $r_0\in]0,+\infty[$ such that $w(r_0)\neq 0$. Then the following statements hold.
\begin{enumerate}
\item[(i)] Let $p\in[1,+\infty]$. If $f\in {\mathcal{M}}^{w}_{p}({\mathbb{R}}^{n})$ and $\epsilon>0$, then the function $f*\phi_{\epsilon}$ from ${\mathbb{R}}^{n}$ to ${\mathbb{R}}$ defined by 
\[
f*\phi_{\epsilon}(x)\equiv\int_{{\mathbb{R}}^{n}}f(x-y)\phi_{\epsilon}(y)\,dy\qquad\forall x\in {\mathbb{R}}^{n}\,,
\]
belongs to ${\mathcal{M}}^{w}_{p}({\mathbb{R}}^{n})\cap C^{\infty}({\mathbb{R}}^{n}) 
$ and
\begin{eqnarray}
 \nonumber
 |f*\phi_{\epsilon}|_{\rho,w,p,  {\mathbb{R}}^{n} }&\leq&
 \int_{{\mathbb{R}}^{n}}|\phi (y)|
\,dy
|f|_{\rho,w,p,{\mathbb{R}}^{n}} \qquad
\forall \rho\in]0,+\infty]\,,
\\  \label{apgm1}
 \|f*\phi_{\epsilon}\|_{  {\mathcal{M}}^{w}_{p}({\mathbb{R}}^{n})  }&\leq &
\int_{ {\mathbb{R}}^{ n} }|\phi|\,dx
\|f\|_{
{\mathcal{M}}^{w}_{p}({\mathbb{R}}^{n})  }\qquad\forall f\in {\mathcal{M}}^{w}_{p}({\mathbb{R}}^{n}).
\end{eqnarray}

\item[(ii)] Let $p\in[1,+\infty]$.  If $f\in {\mathcal{M}}^{w,0}_{p}({\mathbb{R}}^{n})$ and $\epsilon>0$, then $f*\phi_{\epsilon}$ belongs to ${\mathcal{M}}^{w,0}_{p}({\mathbb{R}}^{n})\cap C^{\infty}({\mathbb{R}}^{n}) 
$.
\end{enumerate}
\end{thm}
{\bf Proof.}  Let $f\in {\mathcal{M}}^{w}_{p}({\mathbb{R}}^{n})$, $\epsilon\in]0,+\infty[$.  Since 
$f\in {\mathcal{M}}^{w}_{p}({\mathbb{R}}^{n})\subseteq L^{p}_{ {\mathrm{loc}} }({\mathbb{R}}^{n})$  and  $\phi_{\epsilon}\in C^{\infty}_{c}({\mathbb{R}}^{n}) $,  standard differentiability properties of integrals depending on a parameter imply that
\[
f*\phi_{\epsilon}\in 
 C^{\infty}({\mathbb{R}}^{n})
 \subseteq L^{p}_{ {\mathrm{loc}} }({\mathbb{R}}^{n})
 \,,
\]
(cf.~\textit{e.g.},  Folland~\cite[Prop.~8.10]{Fo99}.) By the Minkowski inequality for integrals, we have
\begin{eqnarray*}
\lefteqn{w(r)\|f*\phi_{\epsilon}\|_{ L^{p}(  {\mathbb{B}}_{n}(x,r) )   }
=w(r)\left\|
\int_{{\mathbb{R}}^{n}}f(\cdot -y)\phi_{\epsilon}(y)\,dy
\right\|_{ L^{p}(  {\mathbb{B}}_{n}(x,r) )   }
}
\\
\nonumber
&& 
\leq w(r)\int_{{\mathbb{R}}^{n}}\left\|
f(\cdot -y)\phi_{\epsilon}(y)
\right\|_{ L^{p}(  {\mathbb{B}}_{n}(x,r) )   }
\,dy
\\
\nonumber
&& 
= w(r)\int_{{\mathbb{R}}^{n}}\left\|
f(\cdot -y)
\right\|_{ L^{p}(  {\mathbb{B}}_{n}(x,r) )   }|\phi_{\epsilon}(y)|
\,dy
\\
\nonumber
&& 
\leq w(r)\int_{{\mathbb{R}}^{n}}|\phi (y)|
\,dy
\sup_{y\in {\mathbb{B}}_{n}(0,\epsilon)}\left\|
f(\cdot -y)\right\|_{ L^{p}(  {\mathbb{B}}_{n}(x,r) )   }
\ \ \forall (x,r)\in {\mathbb{R}}^{n}\times]0,+\infty[\,,
\end{eqnarray*}
(cf.~\textit{e.g.},  Folland~\cite[6.19]{Fo99}.) 
Now by translation invariance of the $L^{p}$-norm, we have
\[
\left\|
f(\cdot -y)\right\|_{ L^{p}(  {\mathbb{B}}_{n}(x,r) )   }
=
\left\|
f \right\|_{ L^{p}(  {\mathbb{B}}_{n}(x-y,r) )   }\,,
\]
for all $x,y\in {\mathbb{R}}^{n}$ and $r\in ]0,+\infty[$. Then the above inequality implies that
\[
w(r)\|f*\phi_{\epsilon}\|_{ L^{p}(  {\mathbb{B}}_{n}(x,r) )   }
\leq 
w(r)\int_{{\mathbb{R}}^{n}}|\phi (y)|
\,dy
\sup_{y\in {\mathbb{B}}_{n}(0,\epsilon)}
\left\|
f \right\|_{ L^{p}(  {\mathbb{B}}_{n}(x-y,r) )   }
\]
for all $(x,r)\in {\mathbb{R}}^{n}\times]0,+\infty[$\,.
Now we have
\[
w(r)\left\|
f \right\|_{ L^{p}(  {\mathbb{B}}_{n}(x-y,r) )   }
\leq\sup_{(z,r)\in {\mathbb{R}}^{n}\times]0,\rho[ }
w(r)\left\|
f \right\|_{ L^{p}(  {\mathbb{B}}_{n}(z,r) )   }
=|f|_{\rho,w,p,{\mathbb{R}}^{n}}
\]
for all $x,y\in  {\mathbb{R}}^{n}$ and $r\in]0,\rho[$ and $\rho\in]0,+\infty]$.  Hence,
\begin{eqnarray}\label{apgmlp1}
\lefteqn{|f*\phi_{\epsilon}|_{\rho,w,p,  {\mathbb{R}}^{n} }=\sup_{(x,r)\in  {\mathbb{R}}^{n}\times]0,\rho[ }
w(r)\|f*\phi_{\epsilon}\|_{ L^{p}(  {\mathbb{B}}_{n}(x,r) )   }
}
\\
\nonumber
&&\qquad\qquad
\leq
\sup_{(x,r)\in  {\mathbb{R}}^{n}\times]0,\rho[ }
 w(r)\int_{{\mathbb{R}}^{n}}|\phi (y)|
\,dy
\sup_{y\in {\mathbb{B}}_{n}(0,\epsilon)}\left\|
f \right\|_{ L^{p}(  {\mathbb{B}}_{n}(x-y,r) )   }
\\
\nonumber
&&\qquad\qquad
\leq
 \int_{{\mathbb{R}}^{n}}|\phi (y)|
\,dy
\sup_{(x,r)\in  {\mathbb{R}}^{n}\times]0,\rho[ }
\left\{
\sup_{(z,r)\in {\mathbb{R}}^{n}\times]0,\rho[ }
w(r)\left\|
f \right\|_{ L^{p}(  {\mathbb{B}}_{n}(z,r) )   }
\right\}
\\
\nonumber
&&\qquad\qquad
 \leq  \int_{{\mathbb{R}}^{n}}|\phi (y)|
\,dy
\sup_{(x,r)\in  {\mathbb{R}}^{n}\times]0,\rho[ }
|f|_{\rho,w,p,{\mathbb{R}}^{n}}
\\
\nonumber
&&\qquad\qquad
 \leq  \int_{{\mathbb{R}}^{n}}|\phi (y)|
\,dy
|f|_{\rho,w,p,{\mathbb{R}}^{n}} \qquad
\forall \rho\in]0,+\infty]\,,
\end{eqnarray}
an inequality which implies the validity of both statement (i) and (ii). \hfill  $\Box$ 

\vspace{\baselineskip}

In the next two sections we introduce two important subspaces of ${\mathcal{M}}^{w }_{p}({\mathbb{R}}^{n})$ of functions that can be approximated by mollifiers. For a characterization of those functions in ${\mathcal{M}}^{w }_{p}({\mathbb{R}}^{n})$ that can be approximated by functions of $C^\infty_c({\mathbb{R}}^n)$, we  refer to Almeida and Samko \cite{AlSa17}, where one can find a more extensive treatment of the problem. See also Sawano, Di Fazio, Hakim \cite[Chap.~8, Vol.~1]{SaDiHa20}.

\section{Zorko's subspace of a generalized Morrey space}

If $p\in[1,+\infty[$, we know that each element of $f\in L^p({\mathbb{R}}^n)$ can be approximated by the family $\{f\ast\phi_\epsilon\}_{\epsilon\in]0,+\infty[}$. 

To prove it, one can estimate the norm $\|f-f\ast\phi_\epsilon\|_{ L^p({\mathbb{R}}^n) }$ in terms of
  $\sup_{y\in{\mathbb{B}}_n(0,\epsilon)}\|f-\tau_y f\|_{ L^p({\mathbb{R}}^n) }$ by exploiting the Minkowski inequality for integrals 
   and remember that the
  translation   map $\tau_{(\cdot)} f $ from ${\mathbb{R}}^n$ to $L^p({\mathbb{R}}^n) $ that takes 
    $y\in{\mathbb{R}}^n$ to $\tau_y f$ is uniformly continuous and that accordingly
 \[
    \lim_{y\to 0}\tau_y f=f\qquad\text{in}\ L^p({\mathbb{R}}^n)\,.
 \]
 We now follow the idea of Zorko \cite{Zo86}   and adopt the same type of argument, but only for  those  functions $f\in {\mathcal{M}}^{w }_{p}({\mathbb{R}}^{n})$ such that
 \[
\lim_{y\to 0}\tau_y f=f\qquad\text{in}\ {\mathcal{M}}^{w }_{p}({\mathbb{R}}^{n})\,.
 \]   
 To do so, we introduce the  Zorko's subspace of ${\mathcal{M}}^{w }_{p}({\mathbb{R}}^{n})$ by means of the following.
 \begin{defn}\label{defn:zorko} Let $p\in[1,+\infty]$. 
  Let $w$ be a function from $]0,+\infty[$ to $[0,+\infty[$.  Assume that there exists $r_0\in]0,+\infty[$ such that $w(r_0)\neq 0$. Then we define as  Zorko's subspace  of  ${\mathcal{M}}^{w }_{p}({\mathbb{R}}^{n})$, the subspace
  \begin{equation}\label{eq:zorko}
{\mathcal{M}}^{w,Z }_{p}({\mathbb{R}}^{n})
\equiv
\left\{
f\in {\mathcal{M}}^{w }_{p}({\mathbb{R}}^{n}):\,
\lim_{y\to 0}\tau_y f=f
\ \text{in}\ {\mathcal{M}}^{w }_{p}({\mathbb{R}}^{n})
\right\}
\end{equation}
of  ${\mathcal{M}}^{w }_{p}({\mathbb{R}}^{n})$. 
\end{defn}
Then we prove the following approximation theorem for the functions in the Zorko's subspace.
\begin{thm}\label{thm:zorkoap}
Let $\phi\in C_{c}^{\infty}({\mathbb{R}}^{ n})$, $\int_{ {\mathbb{R}}^{ n} }\phi\,dx=1$, ${\mathrm{supp}}\,\phi\subseteq\overline{{\mathbb{B}}_n(0,1)}$.  
Let $p\in[1,+\infty]$.   Let $w$ be a function from $]0,+\infty[$ to $[0,+\infty[$.  Assume that there exists $r_0\in]0,+\infty[$ such that $w(r_0)\neq 0$. If $f\in {\mathcal{M}}^{w,Z}_{p}({\mathbb{R}}^{n})
$, then 
\[
f*\phi_{\epsilon}\in 
{\mathcal{M}}^{w}_{p}({\mathbb{R}}^{n})
\cap
 C^{\infty}({\mathbb{R}}^{n})
\]
 for all $\epsilon\in ]0,+\infty[$ and 
\begin{equation}
\label{thm:zorkoap2}
\lim_{\epsilon\to 0}f*\phi_{\epsilon}=f\qquad{\text{in}}\ 
{\mathcal{M}}^{w }_{p}({\mathbb{R}}^{n})
\,.
\end{equation}
\end{thm}
{\bf Proof.} By Theorem \ref{apgm}, we already know that $f*\phi_{\epsilon}\in 
{\mathcal{M}}^{w}_{p}({\mathbb{R}}^{n})
\cap
 C^{\infty}({\mathbb{R}}^{n})
$
 for all $\epsilon\in ]0,+\infty[$. We now prove the limiting relation (\ref{thm:zorkoap2}). If $(x,r)\in {\mathbb{R}}^n\times]0,+\infty[$, then the Minkowski Inequality for integrals implies that
 \begin{eqnarray*}
\lefteqn{
w(r)\|f-f\ast\phi_\epsilon\|_{ L^p( {\mathbb{B}}_n(x,r))  }
}
\\ \nonumber
&&\qquad
=w(r)\left\|
f(z)\int_{{\mathbb{R}}^{n}}\phi_{\epsilon}(y)\,dy-\int_{{\mathbb{R}}^{n}}f(z-y)\phi_{\epsilon}(y)\,dy
\right\|_{ L^p( {\mathbb{B}}_n(x,r)_z)  }
\\ \nonumber
&&\qquad
=w(r)\left\|
\int_{{\mathbb{R}}^{n}}(f(z)-f(z-y))\phi_{\epsilon}(y)\,dy
\right\|_{ L^p( {\mathbb{B}}_n(x,r)_z)  }
\\ \nonumber
&&\qquad
=w(r)\left\|
\int_{{\mathbb{B}}_{n}(0,\epsilon)}(f(z)-\tau_yf(z))\phi_{\epsilon}(y)\,dy
\right\|_{ L^p( {\mathbb{B}}_n(x,r)_z)  }
\\ \nonumber
&&\qquad
\leq w(r) 
\int_{{\mathbb{B}}_{n}(0,\epsilon)}\left\|(f(z)-\tau_yf(z))\phi_{\epsilon}(y)
\right\|_{ L^p( {\mathbb{B}}_n(x,r)_z)  }
\,dy
\\ \nonumber
&&\qquad
\leq 
\sup_{y\in {\mathbb{B}}_{n}(0,\epsilon)} w(r) \left\|f -\tau_yf \right\|_{ L^p( {\mathbb{B}}_n(x,r)) }
\int_{{\mathbb{B}}_{n}(0,\epsilon)}|\phi_{\epsilon}(y)|\,dy
\\ \nonumber
&&\qquad
\leq 
\sup_{y\in {\mathbb{B}}_{n}(0,\epsilon)}   \left\|f -\tau_yf \right\|_{ 
{\mathcal{M}}^{w }_{p}({\mathbb{R}}^{n})
}
\int_{{\mathbb{R}}^{n} }|\phi (y)|\,dy
 \end{eqnarray*}
 where the subscript $z$ of $ {\mathbb{B}}_n(x,r)_z$ is to remind that the $L^p$ norm is to be taken with respect to the variable $z$. Since $f\in {\mathcal{M}}^{w,Z}_{p}({\mathbb{R}}^{n})
$, we have 
\[
\lim_{\epsilon\to 0}\sup_{y\in {\mathbb{B}}_{n}(0,\epsilon)}   \left\|f -\tau_yf \right\|_{ 
{\mathcal{M}}^{w }_{p}({\mathbb{R}}^{n})
}
\]
and thus (\ref{thm:zorkoap2}) holds true.\hfill  $\Box$ 

\vspace{\baselineskip}

\begin{rem}{\em
 By Kato \cite[Cor.~3.3]{Ka92}, 
${\mathcal{M}}^{r^{-\lambda},Z}_{p}({\mathbb{R}}^{n})\subseteq  {\mathcal{M}}^{r^{-\lambda},0}_{p}({\mathbb{R}}^{n})$ for all 
$p\in[1,+\infty[$ and $\lambda\in [0,+\infty[$. By Chiarenza and Franciosi  \cite[Lem.~1.2]{ChFr92}, equality holds if we further assume that $m_n(\Omega)$ is finite. 
}\end{rem}

\section{Approximation properties of $L^p$  functions in generalized Morrey spaces}
In the present section we consider the approximation by smooth functions of those functions of ${\mathcal{M}}^{w }_{p}({\mathbb{R}}^{n})$ that are also $p$-summable.

If $A$ is a subset of a normed space $X$, then ${\mathrm{cl}}_XA$ denotes the closure of $A$ in $X$.
\begin{thm}
\label{apgmlp}
Let $\phi\in C_{c}^{\infty}({\mathbb{R}}^{ n})$, $\int_{ {\mathbb{R}}^{ n} }\phi\,dx=1$, ${\mathrm{supp}}\,\phi\subseteq\overline{{\mathbb{B}}_n(0,1)}$.  Let $w$ be a function from $]0,+\infty[$ to $[0,+\infty[$.  Assume that there exists $r_0\in]0,+\infty[$ such that $w(r_0)\neq 0$. Then the following statements hold.
\begin{enumerate}
\item[(i)] Let $p\in[1,+\infty[$. Let 
\begin{equation}
\label{apgm1a}
\sup_{r\in [\rho,+\infty[}w(r)<+\infty
\qquad\forall\rho\in
]0,+\infty[\,.
\end{equation}
If $f\in {\mathcal{M}}^{w,0}_{p}({\mathbb{R}}^{n})
\cap
L^{p}({\mathbb{R}}^{n})
$, then 
\[
f*\phi_{\epsilon}\in 
{\mathcal{M}}^{w,0}_{p}({\mathbb{R}}^{n})
\cap
L^{p}({\mathbb{R}}^{n})\cap C^{\infty}({\mathbb{R}}^{n})\cap
C^{0}_{ub}({\mathbb{R}}^{n})
\]
 for all $\epsilon\in ]0,+\infty[$ and 
\begin{equation}
\label{apgm2}
\lim_{\epsilon\to 0}f*\phi_{\epsilon}=f\qquad{\text{in}}\ 
{\mathcal{M}}^{w,0}_{p}({\mathbb{R}}^{n})
\cap L^{p}
({\mathbb{R}}^{n})\,.
\end{equation}
\item[(ii)] Let $p\in[1,+\infty[$. Let condition (\ref{apgm1a}) hold. Then
\begin{eqnarray}
\label{apgm3}
\lefteqn{
{\mathrm{cl}}_{{\mathcal{M}}^{w}_{p}({\mathbb{R}}^{n})\cap L^{p}
({\mathbb{R}}^{n})}\left(
{\mathcal{M}}^{w,0}_{p}({\mathbb{R}}^{n})
\cap L^{p}
({\mathbb{R}}^{n})
\cap C^{\infty}({\mathbb{R}}^{n})\cap
C^{0}_{ub}({\mathbb{R}}^{n})
\right)
}
\\
\nonumber
&&\qquad\qquad\qquad\qquad\qquad\qquad\qquad\qquad={\mathcal{M}}^{w,0}_{p}({\mathbb{R}}^{n})\cap L^{p}
({\mathbb{R}}^{n})\,.
\end{eqnarray}
If we further assume that condition 
\begin{equation}
\label{apgmlp3}
\lim_{\rho\to 0}\sup_{r\in]0,\rho]} w(r)r^{n/p} =0\,,
\end{equation}
 holds, then we have
\begin{eqnarray}
\label{apgm3a}
\lefteqn{
{\mathrm{cl}}_{{\mathcal{M}}^{w}_{p}({\mathbb{R}}^{n})\cap L^{p}
({\mathbb{R}}^{n})}\left(
{\mathcal{M}}^{w,0}_{p}({\mathbb{R}}^{n})
\cap L^{p}
({\mathbb{R}}^{n})
\cap C^{\infty}({\mathbb{R}}^{n})\cap
C^{0}_{ub}({\mathbb{R}}^{n})
\right)
}
\\
\nonumber
&&\qquad
=
{\mathrm{cl}}_{{\mathcal{M}}^{w}_{p}({\mathbb{R}}^{n})\cap L^{p}
({\mathbb{R}}^{n})}\left(
{\mathcal{M}}^{w}_{p}({\mathbb{R}}^{n})\cap L^{p}
({\mathbb{R}}^{n})\cap C^{\infty}({\mathbb{R}}^{n})\cap
C^{0}_{ub}({\mathbb{R}}^{n})
\right)
\\
\nonumber
&&\qquad
=
{\mathrm{cl}}_{{\mathcal{M}}^{w}_{p}({\mathbb{R}}^{n})\cap L^{p}
({\mathbb{R}}^{n})}\left(
{\mathcal{M}}^{w}_{p}({\mathbb{R}}^{n})\cap L^{p}
({\mathbb{R}}^{n})\cap L^{\infty}({\mathbb{R}}^{n})
\right)
\\
\nonumber
&&\qquad
={\mathcal{M}}^{w,0}_{p}({\mathbb{R}}^{n})\cap L^{p}
({\mathbb{R}}^{n})\,.
\end{eqnarray}
 \end{enumerate}
 \end{thm}
 {\bf Proof.}  (i) Since $f\in L^{p}({\mathbb{R}}^{n})$ and $\phi_{\epsilon}\in C^{\infty}_{c}({\mathbb{R}}^{n})\subseteq L^{1}({\mathbb{R}}^{n})\cap L^{\infty}({\mathbb{R}}^{n})\subseteq L^{p'}({\mathbb{R}}^n)$, the H\"{o}lder inequality and the Young inequality together with standard differentiability properties of integrals depending on a parameter imply that
\[
f*\phi_{\epsilon}\in 
L^{p}({\mathbb{R}}^{n})\cap C^{\infty}({\mathbb{R}}^{n})\cap
C^{0}_{ub}({\mathbb{R}}^{n})\,,
\]
for all $\epsilon>0$ 
(cf.~\textit{e.g.},  Folland~\cite[Prop.~8.7, 8.8, 8.10]{Fo99}.) By Theorem \ref{apgm} (ii), we also know that
$f*\phi_{\epsilon}\in {\mathcal{M}}^{w,0}_{p}({\mathbb{R}}^{n})$. 

We now turn to show the limiting relation (\ref{apgm2}). Let $\eta>0$. Since $f\in  {\mathcal{M}}^{w,0}_{p}({\mathbb{R}}^{n})$, there exists $\rho_{\eta}>0$ such that 
\[
|f|_{\rho,w,p,{\mathbb{R}}^{n}}\leq\eta
\left(
1+\int_{{\mathbb{R}}^{n}}|\phi (y)|
\,dy
\right)^{-1}
\qquad\forall \rho\in]0,\rho_{\eta}]\,.
\]
Then the triangular inequality and the first inequality of (\ref{apgm1}) imply that
\begin{eqnarray}
\label{apgm4}
\lefteqn{
\sup_{(x,r)\in  {\mathbb{R}}^{n}\times]0,\rho_{\eta}[ }
w(r)\|f-f*\phi_{\epsilon}\|_{ L^{p}(  {\mathbb{B}}_{n}(x,r) )   }
 }
\\
\nonumber
&&\qquad
\leq
\sup_{(x,r)\in  {\mathbb{R}}^{n}\times]0,\rho_{\eta}[ }
w(r)\|f \|_{ L^{p}(  {\mathbb{B}}_{n}(x,r) )   }
+
\sup_{(x,r)\in  {\mathbb{R}}^{n}\times]0,\rho_{\eta}[ }
w(r)\| f*\phi_{\epsilon}\|_{ L^{p}(  {\mathbb{B}}_{n}(x,r) )   }
\\
\nonumber
&&\qquad
\leq |f|_{\rho_{\eta},w,p,{\mathbb{R}}^{n}}
+
\int_{{\mathbb{R}}^{n}}|\phi (y)|
\,dy
|f|_{\rho_{\eta},w,p,{\mathbb{R}}^{n}}
\\
\nonumber
&&\qquad
\leq
|f|_{\rho_{\eta},w,p,{\mathbb{R}}^{n}}  \left(
1+\int_{{\mathbb{R}}^{n}}|\phi (y)|
\,dy
\right)\leq \eta\,,
\end{eqnarray}
for all $\epsilon\in]0,+\infty[$.  
By assumption (\ref{apgm1a}), we have 
\[
\sup_{(x,r)\in  {\mathbb{R}}^{n}\times [\rho_{\eta},+\infty[ }
w(r)\|f-f*\phi_{\epsilon}\|_{ L^{p}(  {\mathbb{B}}_{n}(x,r) )   }
 \leq
 \left(\sup_{r\in [\rho _{\eta},+\infty[}w(r)\right)\|f-f*\phi_{\epsilon}\|_{ L^{p}(  {\mathbb{R}}^{n} )   }\,,
\]
for all $\epsilon>0$.
Since $f\in  L^{p}(  {\mathbb{R}}^{n} )$ and $p\in[1,+\infty[$ and $\int_{ {\mathbb{R}}^{n} }\phi\,dy=1$, standard properties of approximate identities of convolution imply that
\[
\lim_{\epsilon\to 0}\|f-f*\phi_{\epsilon}\|_{ L^{p}(  {\mathbb{R}}^{n} )   }=0\,,
\]
(cf.~\textit{e.g.},  Folland~\cite[Thm.~8.14]{Fo99}.) Then there exists $\epsilon_{\eta}>0$ such that
\[
\left(\sup_{r\in [\rho _{\eta},+\infty[}w(r)\right)\|f-f*\phi_{\epsilon}\|_{ L^{p}(  {\mathbb{R}}^{n} )   }
\leq\eta
\qquad\forall\epsilon\in]0,\epsilon_{\eta} ]\,.
\]
Then we have
\begin{equation}
\label{apgm5}
\sup_{(x,r)\in  {\mathbb{R}}^{n}\times [\rho_{\eta},+\infty[ }
w(r)\|f-f*\phi_{\epsilon}\|_{ L^{p}(  {\mathbb{B}}_{n}(x,r) )   }
 \leq\eta
\qquad\forall\epsilon\in]0,\epsilon_{\eta}]\,.
\end{equation}
By combining inequalities (\ref{apgm4}) and (\ref{apgm5}), we deduce that
\[
\|  f-f*\phi_{\epsilon}\|_{{\mathcal{M}}^{w}_{p}(  {\mathbb{R}}^{n}  )}=
\sup_{(x,r)\in  {\mathbb{R}}^{n}\times ]0,+\infty[ }
w(r)\|f-f*\phi_{\epsilon}\|_{ L^{p}(  {\mathbb{B}}_{n}(x,r) )   }
 \leq\eta
 \]
 for all  $\epsilon\in]0,\epsilon_{\eta}]$. Hence, statement (i) holds true.
 
 Next we consider statement (ii). By Proposition \ref{intml} (iii), 
 ${\mathcal{M}}^{w,0}_{p}(  {\mathbb{R}}^{n}  )\cap L^{p}(  {\mathbb{R}}^{n}  )$ is closed in ${\mathcal{M}}^{w}_{p}(  {\mathbb{R}}^{n}  )\cap L^{p}(  {\mathbb{R}}^{n}  )$. Then we have
 \begin{eqnarray*}
\lefteqn{
{\mathrm{cl}}_{{\mathcal{M}}^{w}_{p}({\mathbb{R}}^{n})\cap
L^{p}({\mathbb{R}}^{n})}\left(
{\mathcal{M}}^{w,0}_{p}({\mathbb{R}}^{n})\cap
L^{p}({\mathbb{R}}^{n})\cap C^{\infty}({\mathbb{R}}^{n})\cap
C^{0}_{ub}({\mathbb{R}}^{n})
\right)
}
\\
\nonumber
&&\qquad\qquad\qquad\qquad\qquad\qquad
\subseteq
{\mathcal{M}}^{w,0}_{p}({\mathbb{R}}^{n})\cap
L^{p}({\mathbb{R}}^{n})\,.
\end{eqnarray*}
On the other hand, if $f\in {\mathcal{M}}^{w,0}_{p}({\mathbb{R}}^{n})\cap
L^{p}({\mathbb{R}}^{n})$, then statement  (i)  implies that 
\[
f\in 
{\mathrm{cl}}_{{\mathcal{M}}^{w}_{p}({\mathbb{R}}^{n})\cap
L^{p}({\mathbb{R}}^{n})}\left(
{\mathcal{M}}^{w,0}_{p}({\mathbb{R}}^{n})\cap
L^{p}({\mathbb{R}}^{n}) \cap C^{\infty}({\mathbb{R}}^{n})\cap
C^{0}_{ub}({\mathbb{R}}^{n})
\right)\,.
\]
If we further assume that condition (\ref{apgmlp3}) holds, then
 Proposition \ref{bgm} implies that
 ${\mathcal{M}}^{w}_{p}({\mathbb{R}}^{n})\cap
L^{\infty}({\mathbb{R}}^{n})\subseteq
{\mathcal{M}}^{w,0}_{p}({\mathbb{R}}^{n})
$, and accordingly we have
\begin{eqnarray*}
\lefteqn{
{\mathrm{cl}}_{{\mathcal{M}}^{w}_{p}({\mathbb{R}}^{n})\cap
L^{p}({\mathbb{R}}^{n})}\left(
{\mathcal{M}}^{w,0}_{p}({\mathbb{R}}^{n})\cap
L^{p}({\mathbb{R}}^{n})\cap C^{\infty}({\mathbb{R}}^{n})\cap
C^{0}_{ub}({\mathbb{R}}^{n})
\right)
}
\\
\nonumber
&&\qquad
\subseteq
{\mathrm{cl}}_{{\mathcal{M}}^{w}_{p}({\mathbb{R}}^{n})\cap
L^{p}({\mathbb{R}}^{n})}\left(
{\mathcal{M}}^{w}_{p}({\mathbb{R}}^{n})\cap
L^{p}({\mathbb{R}}^{n})\cap C^{\infty}({\mathbb{R}}^{n})\cap
C^{0}_{ub}({\mathbb{R}}^{n})
\right)
\\
\nonumber
&&\qquad
\subseteq
{\mathrm{cl}}_{{\mathcal{M}}^{w}_{p}({\mathbb{R}}^{n})\cap
L^{p}({\mathbb{R}}^{n})}\left(
{\mathcal{M}}^{w}_{p}({\mathbb{R}}^{n})\cap
L^{p}({\mathbb{R}}^{n})\cap L^{\infty}({\mathbb{R}}^{n})
\right)
\\
\nonumber
&&\qquad
\subseteq
{\mathcal{M}}^{w,0}_{p}({\mathbb{R}}^{n})\cap
L^{p}({\mathbb{R}}^{n})\,.
\end{eqnarray*}
Then the equality follows as above.
  \hfill  $\Box$

\vspace{\baselineskip}

Then we have the following special case in which we consider weights $w$ such that $\eta_{w}>0$ (cf.~(\ref{prelprgm1_})).
\begin{corol}
\label{capgm}
Let $\phi\in C_{c}^{\infty}({\mathbb{R}}^{ n})$, $\int_{ {\mathbb{R}}^{ n} }\phi\,dx=1$, ${\mathrm{supp}}\,\phi\subseteq\overline{{\mathbb{B}}_n(0,1)}$.  Let $w$ be a function from $]0,+\infty[$ to itself. Let $\eta_{w}>0$.  Then the following statements hold.
\begin{enumerate}
\item[(i)] Let $p\in[1,+\infty]$.  If $f\in {\mathcal{M}}^{w,0}_{p}({\mathbb{R}}^{n})$ and $\epsilon>0$, then $f*\phi_{\epsilon}$ belongs to ${\mathcal{M}}^{w,0}_{p}({\mathbb{R}}^{n})\cap C^{\infty}({\mathbb{R}}^{n})\cap
C^{0}_{ub}({\mathbb{R}}^{n})
$.
\item[(ii)] Let $p\in[1,+\infty[$. Let condition
(\ref{apgm1a}) hold.
If $f\in {\mathcal{M}}^{w,0}_{p}({\mathbb{R}}^{n})$, then
\begin{equation}
\label{capgm2}
\lim_{\epsilon\to 0}f*\phi_{\epsilon}=f\qquad{\text{in}}\ 
{\mathcal{M}}^{w,0}_{p}({\mathbb{R}}^{n})\,.
\end{equation}
\item[(iii)] Let $p\in[1,+\infty[$. Let condition (\ref{apgm1a}) hold. Then
\begin{equation}
\label{capgm3}
{\mathrm{cl}}_{{\mathcal{M}}^{w}_{p}({\mathbb{R}}^{n})}\left(
{\mathcal{M}}^{w,0}_{p}({\mathbb{R}}^{n})\cap C^{\infty}({\mathbb{R}}^{n})\cap
C^{0}_{ub}({\mathbb{R}}^{n})
\right)
={\mathcal{M}}^{w,0}_{p}({\mathbb{R}}^{n})\,.
\end{equation}
If we further assume that condition (\ref{apgmlp3}) holds, then we have
\begin{eqnarray}
\label{capgm3a}
\lefteqn{
{\mathrm{cl}}_{{\mathcal{M}}^{w}_{p}({\mathbb{R}}^{n})}\left(
{\mathcal{M}}^{w,0}_{p}({\mathbb{R}}^{n})\cap C^{\infty}({\mathbb{R}}^{n})\cap
C^{0}_{ub}({\mathbb{R}}^{n})
\right)
}
\\
\nonumber
&&\qquad
=
{\mathrm{cl}}_{{\mathcal{M}}^{w}_{p}({\mathbb{R}}^{n})}\left(
{\mathcal{M}}^{w}_{p}({\mathbb{R}}^{n})\cap C^{\infty}({\mathbb{R}}^{n})\cap
C^{0}_{ub}({\mathbb{R}}^{n})
\right)
\\
\nonumber
&&\qquad
=
{\mathrm{cl}}_{{\mathcal{M}}^{w}_{p}({\mathbb{R}}^{n})}\left(
{\mathcal{M}}^{w}_{p}({\mathbb{R}}^{n})\cap L^{\infty}({\mathbb{R}}^{n})
\right)
={\mathcal{M}}^{w,0}_{p}({\mathbb{R}}^{n})\,.
\end{eqnarray}
\end{enumerate}
\end{corol}
{\bf Proof.} Since $\eta_{w}>0$, we have $w(r_0)>0$ for all $r_0\in]0,+\infty[$ and Proposition \ref{intml}
 (iv) implies that
 \[
 {\mathcal{M}}_{p}^{w}({\mathbb{R}}^{n})
\cap
L^{p}({\mathbb{R}}^{n})={\mathcal{M}}_{p}^{w}({\mathbb{R}}^{n})
\]
 and ${\mathcal{M}}_{p}^{w,0}({\mathbb{R}}^{n})
\cap
L^{p}({\mathbb{R}}^{n})={\mathcal{M}}_{p}^{w,0}({\mathbb{R}}^{n})$ both algebraically and topologically. 
We now consider statement (i). Since $f\in L^{p}({\mathbb{R}}^{n})$ and $\phi_{\epsilon}\in L^{1}({\mathbb{R}}^{n})
\cap L^{\infty}({\mathbb{R}}^{n})$, standard properties of the convolution imply that $f *\phi_{\epsilon}\in C^{0}_{ub}({\mathbb{R}}^{n})$ (cf.~\textit{e.g.},  Folland~\cite[Prop.~8.8]{Fo99}.)  Then statement (i) follows by statement (ii) of Theorem \ref{apgm}. 
Statements  (ii), (iii)  follow immediately from Theorem \ref{apgmlp}.
\hfill  $\Box$

\vspace{\baselineskip}

Then we have the following characterization of the vanishing generalized Morrey spaces on domain. We exploit the notation with $\tilde{C}$, that has been introduced in (\ref{eq:ctildom}).
\begin{thm}
\label{apgmo}
Let $\Omega$ be an open subset of ${\mathbb{R}}^{n}$. 
 Let $p\in[1,+\infty[$.  Let $w$ be a  function from $]0,+\infty[$ to itself. Let $\eta_{w}>0$ (cf.~(\ref{prelprgm1_})).
 Let $\sigma_{w}<+\infty$ (cf.~(\ref{prelprgm2})). Then the following statements hold
 \begin{enumerate}
\item[(i)]If $w$ satisfies condition   (\ref{apgm1a}), then
\begin{eqnarray*}
\lefteqn{
{\mathrm{cl}}_{{\mathcal{M}}^{w}_{p}(\Omega)}\left(
{\mathcal{M}}^{w,0}_{p}(\Omega)\cap \tilde{C}^{\infty}(\overline{\Omega})\cap
C^{0}_{ub}(\overline{\Omega})
\right)
}
\\
\nonumber
&&\qquad\qquad
=
{\mathrm{cl}}_{{\mathcal{M}}^{w}_{p}(\Omega)}\left(
{\mathcal{M}}^{w,0}_{p}(\Omega)\cap C^{\infty}(\overline{\Omega})\cap
C^{0}_{ub}(\overline{\Omega})
\right)
={\mathcal{M}}^{w,0}_{p}(\Omega)
\,.
\end{eqnarray*}
\item[(ii)] If $m_{n}(\Omega)<+\infty$ and if $w$ satisfies conditions (\ref{bfgm1}) and (\ref{bfgm2}), then 
$C^{0}_{ub}(\overline{\Omega})\subseteq L^{\infty}(\Omega)
\subseteq {\mathcal{M}}^{w,0}_{p}(\Omega)$ and
\[
{\mathrm{cl}}_{{\mathcal{M}}^{w}_{p}(\Omega)}\left(
 \tilde{C}^{\infty}(\overline{\Omega})\cap
C^{0}_{ub}(\overline{\Omega})
\right)
=
{\mathrm{cl}}_{{\mathcal{M}}^{w}_{p}(\Omega)}\left(
C^{\infty}(\overline{\Omega})\cap
C^{0}_{ub}(\overline{\Omega})
\right)
={\mathcal{M}}^{w,0}_{p}(\Omega)\,.
\] 
\item[(iii)]  If $\Omega $ is bounded and if $w$ satisfies conditions   (\ref{bfgm1}) and (\ref{bfgm2}), then 
$\tilde{C}^{\infty}(\overline{\Omega})\subseteq
C^{\infty}(\overline{\Omega})\subseteq C^{0}_{ub}(\overline{\Omega})\subseteq L^{\infty}(\Omega)
\subseteq {\mathcal{M}}^{w,0}_{p}(\Omega)$ and
\[
{\mathrm{cl}}_{{\mathcal{M}}^{w}_{p}(\Omega)} 
 \tilde{C}^{\infty}(\overline{\Omega}) 
 =
 {\mathrm{cl}}_{{\mathcal{M}}^{w}_{p}(\Omega)} 
C^{\infty}(\overline{\Omega}) 
={\mathcal{M}}^{w,0}_{p}(\Omega)\,.
\] 
\end{enumerate}
\end{thm}
{\bf Proof.} By Theorem \ref{lmc}, ${\mathcal{M}}^{w,0}_{p}(\Omega)$ is closed in ${\mathcal{M}}^{w}_{p}(\Omega)$. Then we have 
\begin{eqnarray*}
\lefteqn{
{\mathrm{cl}}_{{\mathcal{M}}^{w}_{p}(\Omega)}\left(
{\mathcal{M}}^{w,0}_{p}(\Omega)\cap \tilde{C}^{\infty}(\overline{\Omega})\cap
C^{0}_{ub}(\overline{\Omega})
\right)
}
\\
\nonumber
&&\qquad\qquad
\subseteq
{\mathrm{cl}}_{{\mathcal{M}}^{w}_{p}(\Omega)}\left(
{\mathcal{M}}^{w,0}_{p}(\Omega)\cap C^{\infty}(\overline{\Omega})\cap
C^{0}_{ub}(\overline{\Omega})
\right)
\subseteq
{\mathcal{M}}^{w,0}_{p}(\Omega)\,.
\end{eqnarray*}
Now let $f\in {\mathcal{M}}^{w,0}_{p}(\Omega)$. By Proposition \ref{prelprgm} (ii), we have $E_{\Omega}f\in {\mathcal{M}}^{w,0}_{p}({\mathbb{R}}^{n})$. By Corollary \ref{capgm} (iii), there exists a sequence $\{g_{j}\}_{j\in {\mathbb{N}} }$ in 
\[
{\mathcal{M}}^{w,0}_{p}({\mathbb{R}}^{n})\cap C^{\infty}({\mathbb{R}}^{n})\cap
C^{0}_{ub}({\mathbb{R}}^{n})
\]
 such that
\[
\lim_{j\to\infty} g_{j}=E_{\Omega}f\qquad{\mathrm{in}}\ 
{\mathcal{M}}^{w}_{p}({\mathbb{R}}^{n})\,.
\]
Then Proposition \ref{prelprgm} (i)  implies that
\[
f_{j}\equiv (g_{j})_{|\Omega}\in
{\mathcal{M}}^{w,0}_{p}(\Omega)\cap \tilde{C}^{\infty}(\overline{\Omega})\cap
C^{0}_{ub}(\overline{\Omega})\qquad\forall j\in{\mathbb{N}}\,.
\]
Hence, $f\in {\mathrm{cl}}_{{\mathcal{M}}^{w}_{p}(\Omega)}\left(
{\mathcal{M}}^{w,0}_{p}(\Omega)\cap \tilde{C}^{\infty}(\overline{\Omega})\cap
C^{0}_{ub}(\overline{\Omega})
\right)
$ and statement (i) follows. 

If $m_{n}(\Omega)<+\infty$ and if $w$ satisfies conditions (\ref{bfgm1}) and (\ref{bfgm2}), then Corollary \ref{bfgm} implies that
$L^{\infty}(\Omega)\subseteq {\mathcal{M}}^{w,0}_{p}(\Omega)$. Then we have
\[
C^{0}_{ub}(\overline{\Omega})\subseteq L^{\infty}(\Omega)
\subseteq {\mathcal{M}}^{w,0}_{p}(\Omega)
\]
 and
\[
{\mathcal{M}}^{w,0}_{p}(\Omega)\cap \tilde{C}^{\infty}(\overline{\Omega})\cap
C^{0}_{ub}(\overline{\Omega})
=  \tilde{C}^{\infty}(\overline{\Omega})\cap  C^{0}_{ub}(\overline{\Omega})
\]
and
\[
{\mathcal{M}}^{w,0}_{p}(\Omega)\cap C^{\infty}(\overline{\Omega})\cap
C^{0}_{ub}(\overline{\Omega})
=  C^{\infty}(\overline{\Omega})\cap  C^{0}_{ub}(\overline{\Omega})\,.
\]
Since conditions (\ref{bfgm1}) and (\ref{bfgm2}) imply condition  (\ref{apgm1a}), 
statement (i) implies the validity of statement (ii).

If we further assume that $\Omega$  is bounded, then $\overline{\Omega}$ is compact and
\[
\tilde{C}^{\infty}(\overline{\Omega})\subseteq 
C^{\infty}(\overline{\Omega})
\subseteq
C^{0}_{ub}(\overline{\Omega})\subseteq L^{\infty}(\Omega)
\subseteq {\mathcal{M}}^{w,0}_{p}(\Omega)\,,
\]
and
$ \tilde{C}^{\infty}(\overline{\Omega})\cap
C^{0}_{ub}(\overline{\Omega})
=  \tilde{C}^{\infty}(\overline{\Omega})$
and $ C^{\infty}(\overline{\Omega})\cap
C^{0}_{ub}(\overline{\Omega})
=  C^{\infty}(\overline{\Omega})$. Hence, statement (iii) follows by statement (ii). \hfill  $\Box$

\vspace{\baselineskip}

Next we introduce the corresponding statement for the classical Morrey spaces. 

\begin{corol}
\label{apgml}
Let $\Omega$ be an open subset of ${\mathbb{R}}^{n}$. 
 Let $p\in[1,+\infty[$. Let $\lambda\in [0,n/p[$. Then the following statements hold.
  \begin{enumerate}
\item[(i)]
\begin{eqnarray*}
\lefteqn{
{\mathrm{cl}}_{M^{\lambda}_{p}(\Omega)}\left(
M^{\lambda,0}_{p}(\Omega)\cap \tilde{C}^{\infty}(\overline{\Omega})\cap
C^{0}_{ub}(\overline{\Omega})
\right)
}
\\
\nonumber
&&\qquad\qquad
={\mathrm{cl}}_{M^{\lambda}_{p}(\Omega)}\left(
M^{\lambda,0}_{p}(\Omega)\cap C^{\infty}(\overline{\Omega})\cap
C^{0}_{ub}(\overline{\Omega})
\right)
=M^{\lambda,0}_{p}(\Omega)
\,.
\end{eqnarray*}
\item[(ii)] If $m_{n}(\Omega)<+\infty$, then
$C^{0}_{ub}(\overline{\Omega})\subseteq L^{\infty}(\Omega)
\subseteq M^{\lambda,0}_{p}(\Omega)$ and
\[
{\mathrm{cl}}_{M^{\lambda}_{p}(\Omega)}\left(
 \tilde{C}^{\infty}(\overline{\Omega})\cap
C^{0}_{ub}(\overline{\Omega})
\right)
=
{\mathrm{cl}}_{M^{\lambda}_{p}(\Omega)}\left(
C^{\infty}(\overline{\Omega})\cap
C^{0}_{ub}(\overline{\Omega})
\right)
=M^{\lambda,0}_{p}(\Omega)\,.
\] 
\item[(iii)]  If $\Omega $ is bounded, then 
$\tilde{C}^{\infty}(\overline{\Omega})\subseteq
C^{\infty}(\overline{\Omega})\subseteq C^{0}_{ub}(\overline{\Omega})\subseteq L^{\infty}(\Omega)
\subseteq M^{\lambda,0}_{p}(\Omega)$ and
\[
{\mathrm{cl}}_{M^{\lambda}_{p}(\Omega)} 
 \tilde{C}^{\infty}(\overline{\Omega}) 
 =
 {\mathrm{cl}}_{M^{\lambda}_{p}(\Omega)} 
 C^{\infty}(\overline{\Omega}) 
=M^{\lambda,0}_{p}(\Omega)\,.
\] 
\end{enumerate}
\end{corol}
{\bf Proof.} It suffices to note that if $\lambda\in [0,n/p[$, then $w_{\lambda}$ satisfies (\ref{bfgm1}), (\ref{bfgm2}) and (\ref{apgm1a}) and to apply Theorem \ref{apgmo}. \hfill  $\Box$

\vspace{\baselineskip}

By Corollary \ref{apgml} (iii) the vanishing Morrey space on a bounded domain $M^{\lambda,0}_{p}(\Omega)$ coincides with the completion of $ C^{\infty}(\overline{\Omega}) $ with respect to the norm 
$\|\cdot\|_{ M^{\lambda}_{p}(\Omega) }$. Thus the Morrey space on a bounded domain of authors as Campanato~\cite{Ca63} coincides
 with our vanishing Morrey space $M^{\lambda,0}_{p}(\Omega)$.
 
 Instead, our definition of Morrey space $M^{\lambda}_{p}(\Omega)$ 
 on a bounded domain is equivalent to that of Kufner, John and Fu\v cik~\cite{KuJoFu77}.

\appendix

\section{Appendix on measure theory}
 
Let $\Omega$ be an open subset of ${\mathbb{R}}^{n}$. Let $p\in [1,+\infty]$. Let $L^{p}_{{\mathrm{loc}} }(\Omega )$ be the vector space   of the equivalence classes of the measurable functions $f$ from $\Omega$ to ${\mathbb{R}}$ such that 
$f_{|K}\in L^{p}(K)$ for all compact subsets $K$ of $\Omega$. 
 Let $K$ be a compact subset of $\Omega$. Then we set
\[
p_{K}\equiv\|f_{|K}\|_{L^{p}(K)}
\qquad\forall f\in L^{p}_{{\mathrm{loc}} }(\Omega )
\]
As one can easily verify, $p_{K}$ is a seminorm, but not a norm. The family
\[
{\mathcal{P}}\equiv \left\{
p_{K}:\,K\ {\mathrm{is\ compact}},\ K\subseteq \Omega
\right\}
\]
induces a locally convex  topology on $L^{p}_{{\mathrm{loc}}}(\Omega )$. Then one can prove that the topology of $L^{p}_{{\mathrm{loc}}}(\Omega )$ is metrizable and that $L^{p}_{{\mathrm{loc}}}(\Omega )$ is complete.

A sequence $\{ f_{j} \}_{ j\in {\mathbb{N}} }$   in the space $L^{p}_{{\mathrm{loc}}}(\Omega )$ is a Cauchy sequence  if and only if $\{ (f_{j})_{|K} \}_{ j\in {\mathbb{N}} }$ is a Cauchy sequence in $L^{p}(K)$ for all compact subsets $K$ of 
$\Omega$.

A sequence $\{ f_{j} \}_{ j\in {\mathbb{N}} }$ converges to $f$ in the space $L^{p}_{{\mathrm{loc}}}(\Omega )$ if and only if
\[
\lim_{j\to\infty}(f_{j})_{|K}=f_{|K}\qquad{\mathrm{in}}\ L^{p}(K)
\]
for all compact subsets $K$ of $\Omega$.

As a consequence, one can easily prove that if a sequence $\{ f_{j} \}_{ j\in {\mathbb{N}} }$ converges to $f$ in the space $L^{p}_{{\mathrm{loc}}}(\Omega )$, then there exists a subsequence $\{ f_{j_{k}} \}_{k\in {\mathbb{N}} }$ that converges to $f$ pointwise almost everywhere in $\Omega$.

Finally, in $X$ is a normed space and if $T$ is a linear map from $X$ to $L^{p}_{{\mathrm{loc}}}(\Omega )$, then $T$ is continuous if and only if for each compact subset $K$ of $\Omega$ there exists a constant $c_{K}>0$ such that
\[
\|T[x]\|_{
L^{p} (K )
}\leq c_{K}\|x\|_{X}\qquad\forall x\in X\,.
\]

Next we introduce the following classical result (cf.~\textit{e.g.}, Folland~\cite[Ex.~15, p.~187]{Fo99}. We also mention the paper of Fichera~\cite{Fi43}.)
\begin{thm}[Vitali Convergence]
\label{vitali}
Let $(X,{\mathcal{M}},\mu)$ be a measured sp\-ace with finite measure. Let $\{f_{j}\}_{j\in{\mathbb{N}}}$ be a sequence of functions in $L^{1}_{\mu}(X)$ which converges pointwise to a measurable function $f\in {\mathbb{R}}^{X}$. If  for each $\epsilon>0$ there exists $\delta>0$ such that
\[
\int_{E}\left|f_{j}\right|\,d\mu
\leq \epsilon
\qquad \forall E\in {\mathcal{M}}, \mu(E)\leq\delta, \forall j\in {\mathbb{N}}\,,
\]
then  $f\in L^{1}_{\mu}(X)$ and $\lim_{j\to\infty}f_{j}=f$ in $L^{1}_{\mu}(X)$.
\end{thm}

Then we have the following  important result (cf.~\textit{e.g.}, Folland~\cite[Th.~6.19, p.~194]{Fo99}.)
\begin{thm}[Minkowski inequality for integrals]\label{minint}
Let $(X,{\mathcal{M}},\mu)$ and $(Y,{\mathcal{N}},\nu)$ be $\sigma$-finite  measured spaces. Let $f$ be a measurable function from $(X\times Y,{\mathcal{M}}\otimes {\mathcal{N}})$ to ${\mathbb{R}}$. Let $p\in [1,+\infty]$. Then the following statements hold.  
\begin{enumerate}
\item[(i)] If $f$ is a measurable function from $(X\times Y,{\mathcal{M}}\otimes {\mathcal{N}})$ to $[0,+\infty[$, then there exist $h\in {\mathcal{L}}^{+}(X)$ and $k_{p}\in 
{\mathcal{L}}^{+}(Y)$ such that 
\begin{eqnarray}
\label{minint0}
&&
h(x)=\int_{Y}f(x,y)\,d\nu (y)\qquad{\mathrm{a.a}}\ x\in X\,,
\\  \nonumber
&&
k_{p}(y)\equiv
\left\{
\begin{array}{ll}
\left(
\int_{X}f^{p}(x,y)\,d\mu(x)
\right)^{1/p} & {\mathrm{a.a.}}\ y\in Y\qquad {\mathrm{if}}\ p\in[1,+\infty[\,,
\\ \nonumber
{\mathrm{ess\,sup}}_{x\in X}f(x,y) & {\mathrm{a.a.}}\ y\in Y\qquad {\mathrm{if}}\ p=+\infty\,,
\end{array}
\right.
\end{eqnarray}
and
\begin{eqnarray}\nonumber
&&\left(
\int_{X}h^{p}(x)\,d\mu(x)
\right)^{1/p}
\leq\int_{Y}k_{p}(y)\,d\nu(y)
\qquad {\mathrm{if}}\ p\in[1,+\infty[\,,
\\ \label{minint1}
&&
{\mathrm{ess\,sup}}_{x\in X} h(x)\leq
 \int_{Y}k_{\infty}(y)\,d\nu(y)
\qquad {\mathrm{if}}\ p=+\infty\,,
\end{eqnarray}
 which we rewrite in the form
 \begin{eqnarray*}
&&\left(
\int_{X}\left(\int_{Y}f(x,y)\,d\nu (y)\right)^{p}\,d\mu(x)
\right)^{1/p}
\\
&&\qquad\leq\int_{Y}\left(
\int_{X}f^{p}(x,y)\,d\mu(x)
\right)^{1/p}\,d\nu(y) \qquad 
\qquad {\mathrm{if}}\ p\in[1,+\infty[\,,
\\
&&
{\mathrm{ess\,sup}}_{x\in X}\int_{Y}f(x,y)\,d\nu (y)\leq
 \int_{Y}{\mathrm{ess\,sup}}_{x\in X}f(x,y)\,d\nu(y)\\
&&\qquad\qquad\qquad\qquad\qquad\qquad\qquad\qquad\qquad\qquad
\qquad {\mathrm{if}}\ p=+\infty\,.
\end{eqnarray*}
\item[(ii)] If $f$ is a measurable function from $(X\times Y,{\mathcal{M}}\otimes {\mathcal{N}})$ to ${\mathbb{R}}$
such that
\begin{enumerate}
\item[(j)] $f(\cdot,y)\in {\mathcal{L}}^{p}_{\mu}(X)$ for almost all $y\in Y$.
\item[(jj)] There exists $k_{p}\in {\mathcal{L}}^{1}_{\nu}(Y)$ such that
\[
k_{p}(y)=\|f(\cdot,y)\|_{   {\mathcal{L}}^{p}_{\mu}(X)  }
\qquad {\mathrm{a.a.}}\ y\in Y\,.
\]
\end{enumerate}
Then $f(x,\cdot)\in {\mathcal{L}}^{1}_{\nu}(Y)$  for almost all $x\in X$ and there exists $h\in  {\mathcal{L}}^{p}_{\mu}(X)$ such that
\[
h(x)=\int_{Y}f(x,y)\,d\nu(y)\qquad{\mathrm{a.a.}}\ x\in X\,,
\]
and
\begin{equation}
\label{minint2}
\|h\|_{  {\mathcal{L}}^{p}_{\mu}(X)  }
\leq \int_{Y}k_{p}(y)\,d\nu(y)\,,
\end{equation}
an inequality which we also write in the form
\[
\left\|
\int_{Y}f(\cdot,y)\,d\nu(y)
\right\|_{  {\mathcal{L}}^{p}_{\mu}(X)  }
\leq \int_{Y}
\|
f(\cdot,y)
\|_{  {\mathcal{L}}^{p}_{\mu}(X)  }
\,d\nu(y)\,.
\]
\end{enumerate}
\end{thm}

\section{Appendix on Functional Analysis}
As is well known, if $({\mathcal{X}},\|\cdot\|_{{\mathcal{X}}})$
and $({\mathcal{Y}},\|\cdot\|_{{\mathcal{Y}}})$ are Banach spaces, then the function $\|\cdot\|_{ {\mathcal{X}}\cap {\mathcal{Y}} }$
from ${\mathcal{X}}\cap {\mathcal{Y}}$ to $[0,+\infty[$ defined by
\[
\|\xi\|_{ {\mathcal{X}}\cap {\mathcal{Y}} }
\equiv
\max
\left\{
\|\xi\|_{ {\mathcal{X}} }
,
\|\xi\|_{ {\mathcal{Y}} }
\right\}
\qquad\forall  
\xi\in {\mathcal{X}}\cap {\mathcal{Y}}\,,
\]
is a norm on ${\mathcal{X}}\cap {\mathcal{Y}}$ and $({\mathcal{X}}\cap {\mathcal{Y}},\|\cdot\|_{ {\mathcal{X}}\cap {\mathcal{Y}} } )$ is a Banach space. Moreover, $({\mathcal{X}}\cap {\mathcal{Y}},\|\cdot\|_{ {\mathcal{X}}\cap {\mathcal{Y}} } )$ is continuously embedded  both into ${\mathcal{X}}$ and into ${\mathcal{Y}}$.

\vspace{\baselineskip}

\noindent
{\bf Acknowledgement} M. Lanza de Cristoforis wishes to express his gratitude to
Prof.~Alexey Karapetyants and to the Doctoral School in Mathematics of the
Southern Federal State University of Rostov-on-Don for the  warm hospitality.

\end{document}